\documentclass[11pt]{article}

\usepackage{titlesec}
\titleformat{\subsection}[hang]{\normalfont\bfseries}{\thesubsection}{1em}{}
\usepackage{titlesec}

\titlespacing\section{0pt}{3.5ex plus 0.5ex minus .2ex}{0.3ex plus .2ex}
\titlespacing\subsection{0pt}{2.5ex plus 0.5ex minus .2ex}{0.3ex plus .2ex}
\titlespacing\subsubsection{0pt}{2.5ex plus 0.5ex minus .2ex}{0.3ex plus .2ex}
\usepackage[latin1]{inputenc}
\usepackage{amsmath}
\usepackage{amsfonts}
\usepackage{amssymb}
\usepackage{amsthm}

\usepackage{mathrsfs} % for \mathscr (script letters)

\usepackage{multind}
\makeindex{notation}
\makeindex{definition}

\usepackage[alphabetic,lite]{amsrefs} % for different reference style
\usepackage{fullpage}
\usepackage{setspace}
\usepackage{hyperref}
\usepackage{color}
\usepackage[shortlabels]{enumitem} % to change bullet points for itemize, seem to includes enumerate (after adding shortlabels)
\usepackage{ulem} %um text durchzustreichen

\usepackage{comment} %to comment library out
\excludecomment{Jessi}

\usepackage[disable, color=orange!50]{todonotes} % to insert notes via \todo[options]{todo-text} e.g. \todo[inline, color=yellow!20]{Add update}
\newcommand{\todor}{\todo[size=\scriptsize, color=red!20]}

\usepackage{marginnote}

\usepackage[all,cmtip]{xy} % crazy diagrams
 \usepackage{fancyhdr}
\newtheorem{Thm}{Theorem}[subsection]
\newtheorem{Lemma}[Thm]{Lemma}
\newtheorem{Cor}[Thm]{Corollary}

\newtheorem{Prop}[Thm]{Proposition}
\newtheorem{Assumption}[Thm]{Assumption}
\newtheorem{DefProp}[Thm]{Definition / Proposition}

\theoremstyle{definition}
\newtheorem{Def}[Thm]{Definition}
\newtheorem{Rem}[Thm]{Remark}%helpful commands

\newcommand{\Tr}{\text{Tr}}
\newcommand{\<}{\left\langle}
\renewcommand{\>}{\right\rangle}
\newcommand{\al}{\alpha}

\newcommand{\bF}{\mathbb{F}}
\newcommand{\bG}{\mathbb{G}}

\newcommand{\bN}{\mathbb{N}}
\newcommand{\bP}{\mathbb{P}}
\newcommand{\bQ}{\mathbb{Q}}
\newcommand{\bR}{\mathbb{R}}

\newcommand{\bZ}{\mathbb{Z}}

\newcommand{\cB}{\mathcal{B}}

\newcommand{\cF}{\mathcal{F}}

\newcommand{\cO}{\mathcal{O}}

\newcommand{\cX}{\mathcal{X}}

\newcommand{\ff}{\mathfrak{f}}
\newcommand{\fg}{\mathfrak{g}}

\newcommand{\fp}{\mathfrak{p}}

\newcommand{\ft}{\mathfrak{t}}
\newcommand{\fu}{\mathfrak{u}}

\newcommand{\sA}{\mathscr{A}}
\newcommand{\sB}{\mathscr{B}}

\newcommand{\sG}{\mathscr{G}}
\newcommand{\sH}{\mathscr{H}}

\newcommand{\sS}{\mathscr{S}}
\newcommand{\sT}{\mathscr{T}}
\newcommand{\sU}{\mathscr{U}}
\newcommand{\sV}{\mathscr{V}}

\DeclareMathAlphabet{\mathpzc}{OT1}{pzc}{m}{it}

\newcommand{\cx}{\mathpzc{x}}

\newcommand{\ra}{\rightarrow}

\newcommand{\Ra}{\Rightarrow}

\newcommand{\xra}{\xrightarrow}

\newcommand{\wt}{\widetilde}

\newcommand{\eps}{\epsilon}

\newcommand{\ov}{\overline}
\newcommand{\un}{\underline}

\providecommand{\abs}[1]{\left\lvert#1\right\rvert}

\DeclareMathOperator{\Image}{Im}		% Image
\DeclareMathOperator{\Spec}{Spec}		% Spectrum of a ring
\DeclareMathOperator{\Hom}{Hom}			% Set of arrows between two object
\DeclareMathOperator{\Id}{Id}		  	% Identity
\DeclareMathOperator{\SL}{SL}		  	% SL
\DeclareMathOperator{\SU}{SU} 			% SU, special unitary group
\DeclareMathOperator{\Dyn}{Dyn}			% Dynkin diagram
\DeclareMathOperator{\Gal}{Gal}			% Galois group
\DeclareMathOperator{\Aut}{Aut}			% automorhpism group
\DeclareMathOperator{\Ad}{Ad}	  		% Adjoint
\DeclareMathOperator{\diag}{diag}  	% diagonal matrix
\DeclareMathOperator{\Lie}{Lie}  	  % Lie algebra
\DeclareMathOperator{\ind}{ind}  	  % compact induction
\DeclareMathOperator{\Res}{Res}  	  % restriction
\DeclareMathOperator{\Sym}{Sym}  	  % Sym
\newcommand{\Waff}{W_{\mathrm{aff}}}  	  % affine Weil group
\DeclareMathOperator{\Pic}{Pic}  	  % Picard group 
\DeclareMathOperator{\Stab}{Stab}  	  % stabilizer
\newcommand{\fpb}{\overline{{\mathbb F}}_p}  	  % F-p-bar
\newcommand{\ba}{\begin{aligned}}
	\newcommand{\ea}{\end{aligned}}
\newcommand{\beqn}{\begin{eqnarray}}
	\newcommand{\eeqn}{\end{eqnarray}}
\newcommand{\beqns}{\begin{eqnarray*}}
	\newcommand{\eeqns}{\end{eqnarray*}}
\newcommand{\benum}{\begin{enumerate}}
	\newcommand{\eenum}{\end{enumerate}}

\newcommand{\ZG}{\sH} % for the integral model, maybe I want to change the notation?
\newcommand{\ZV}{\sV}
\newcommand{\ZVold}{\wt \sV}
\newcommand{\ZVnm}{\ZV_{\mathrm{nm}}}
\newcommand{\ZVmul}{\ZV_{\mathrm{mul}}}
\newcommand{\ZzV}{\sV'}%{\sV^{E'}}

\newcommand{\ZS}{\sS}

\newcommand{\ZZU}{\sU}
\newcommand{\ZzU}{\sU}%{{\sU^{E'}}}
\newcommand{\ZzB}{\sB}%{{\sB^{E'}}}
\newcommand{\ZB}{\sB_H}%{{\sB^{E'}}}
\newcommand{\ZZG}{\sG}
\newcommand{\ZzG}{{\sH'}}%{{\sH ^{E'}}}
\newcommand{\ZZT}{\sT}
\newcommand{\Zx}{{\cx_H}}
\newcommand{\ZZx}{\cx}
\newcommand{\Rx}{\ov x}
\newcommand{\RX}{\ov X}
\newcommand{\Rxq}{\ov x_q}

\newcommand{\ZZX}{\cX}

\newcommand{\PP}{\bP} % parahoric groupscheme
\newcommand{\RP}{\mathbf{G}} %reductive quotient of parahoric subgroup
\newcommand{\RPF}{\mathbf{G}^{F}} %reductive quotient of parahoric subgroup for G_F
\newcommand{\RPFq}{\mathbf{G}^{F_q}} %reductive quotient of parahoric subgroup for G_F
\newcommand{\PV}{\mathbf{V}} % higher Moy--Prasad quotients
\newcommand{\PVF}{\mathbf{V}^{F}} % higher Moy--Prasad quotients for G_E' 
\newcommand{\PVFq}{\mathbf{V}^{F_q}} % higher Moy--Prasad quotients for G_E' 
\newcommand{\RU}{\mathbf{U}} % root groups of reductive quotient
\newcommand{\RUq}{\mathbf{U_q}} % root groups of reductive quotient
\newcommand{\RUF}{\mathbf{U}^{\Fp}} % root groups of reductive quotient over F
\newcommand{\RUFq}{\mathbf{U}^{\Fq}} % root groups of reductive quotient over Fq
\newcommand{\RT}{\mathbf{T}} % split maximal torus of reductive quotient
\newcommand{\RTF}{\mathbf{T}^\Fp} % split maximal torus of reductive quotient
\newcommand{\RTFq}{\mathbf{T}^\Fq} % split maximal torus of reductive quotient

\newcommand{\mmu}{\boldsymbol\mu}
\DeclareMathOperator{\Conjug}{Inn} 
\newcommand{\me}{{e'}} % m=\me \cdot \mf 
\newcommand{\mf}{f} % m=\me \cdot \mf 
\newcommand{\tl}{\theta_\lambda} % map from mu_m into Aut(G)
\newcommand{\out}{\vartheta} % outer endomorphism
\newcommand{\UE}{U^E} % root group over $E$
\newcommand{\xE}{x^E} % root group over $E$
\newcommand{\xF}{x^{\Fp}} % root group over $E$
\newcommand{\xFq}{x^{\Fq}} % root group over $E$
\newcommand{\GE}{G^E} % Moy--Prasad filtration group over $E$
\newcommand{\ftE}{\ft^{E_t}} % Moy--Prasad filtration over $E$
\newcommand{\GF}{G^F} % Moy--Prasad filtration group over $F$
\newcommand{\GFq}{G^{F_q}} % Moy--Prasad filtration group over $F$
\newcommand{\fgF}{\fg^F} % Moy--Prasad filtration over $F$
\newcommand{\fgFq}{\fg^{F_q}} % Moy--Prasad filtration over $F$
\newcommand{\va}{\mathrm{v}} % valuation
\newcommand{\Root}{\Psi_K} % affine roots
\newcommand{\Rootq}{\Psi_{K_q}} % affine roots
\newcommand{\good}{good} % good, property of groups, to get correct spacing, write \good{}
\newcommand{\Fp}{F} % field to define groups over other primes, maybe change letter b/c F was used generic in earlier sections
\newcommand{\Fq}{{F_q}} %
\newcommand{\Eq}{{E_q}} %
\newcommand{\Gammapr}{\Gamma}
\newcommand{\gammapr}{\gamma}
\newcommand{\zetaG}{\zeta_{G}}
\newcommand{\zetaGq}{\zeta_{G_q}}
\newcommand{\iotaQ}{\iota_{\ov \bQ}}
\newcommand{\Et}{E \cap K^{\rm{tame}}}

\newcommand{\OH}{\cO_H} %ring of integers in Hilbert class field
\newcommand{\Ze}{\OEpr[1/e]}
\newcommand{\OEpr}{\cO_{e'}} %ring of integers i
\newcommand{\Phixr}{\Phi_{x,r}^{\rm{max}}}
\newcommand{\fZV}{f_{\ZV}}
\newcommand{\Nprime}{N'}

\newcommand{\Kal}{{\overline{K}}} % algebraic closure of K
\DeclareMathOperator{\fqb}{\ov \bF_q}  	  % F-p-bar

\AtEndDocument{\bigskip{\footnotesize%
		\textsc{Jessica Fintzen, Department of Mathematics, Harvard University, One Oxford Street, Cambridge, MA 02138, USA} \\
		\textit{Current address}: Trinity College, Trinity Street, Cambridge, CB2 1TQ, UK\\
		\textit{E-mail address}: \texttt{fintzen@umich.edu} 
	}}

\begin{document}
	\author{Jessica Fintzen} 
	\title{On the Moy--Prasad filtration}
	\date{}
	\maketitle
	
	\begin{abstract}Let $K$ be a maximal unramified extension of a nonarchimedean local field with arbitrary residual characteristic $p$. Let $G$ be a reductive group over $K$ which splits over a tamely ramified extension of $K$. We show that the associated Moy--Prasad filtration representations are in a certain sense independent of $p$. We also establish descriptions of these representations in terms of explicit Weyl modules and as representations occurring in a generalized Vinberg--Levy theory.

		As an application, we use these results to provide necessary and sufficient conditions for the existence of stable vectors in Moy--Prasad filtration representations, which extend earlier results by Reeder and Yu (which required $p$ to be large) and by Romano and the author (which required $G$ to be absolutely simple and split). This yields new supercuspidal representations.

		We also treat reductive groups $G$ that are not necessarily split over a tamely ramified field extension.		

	\end{abstract}

	{
		\renewcommand{\thefootnote}{}  % to delete the footnote number 	
		\footnotetext{MSC: Primary 20G25, 20G07, 22E50, 14L15; Secondary 11S37, 14L24 }
		\footnotetext{Keywords: Moy--Prasad filtration, reductive group schemes, stable vectors, supercuspidal representations, Weyl modules}
		\footnotetext{The author was partially supported by the Studienstiftung des deutschen Volkes.}
	}

	\newpage
	\tableofcontents
	\newpage

\section{Introduction} \label{Sec-Intro}

The introduction of Moy--Prasad filtrations in the 1990s revolutionized the study of the representation theory of $p$-adic groups. As one example, their introduction enabled a construction of supercuspidal representations -- the building blocks in the representation theory of $p$-adic groups -- that is exhaustive for large primes $p$ under certain tameness assumptions. However, while this and similar advances are remarkable, 
the restrictions on the prime $p$ are unsatisfying.
Given their critical role, we expect that a better understanding of the Moy--Prasad filtrations will be a key ingredient for future progress.
To that end, we introduce a ``global'' model for the Moy--Prasad filtration quotients. This allows us to compare the Moy--Prasad filtrations for different primes $p$ and to deduce results for \textit{all} primes $p$ that were previously only known for large primes. Our global model also enables us to express the Moy--Prasad filtration quotients in terms of more traditional, well studied concepts, e.g. as explicit Weyl modules or in terms of a generalized Vinberg--Levy theory. As an application, we exhibit new supercuspidal representations for non-split $p$-adic groups, including non-tame groups.

To explain the content and background of the paper in more detail, let us introduce some notation.
Let $k$ be a nonarchimedean local field 
 with residual characteristic $p>0$. Let $K$ be a maximal unramified extension of $k$ and identify its residue field with $\fpb$. Let $G$ be a reductive group over $K$. In \cite{BT1, BT2}, Bruhat and Tits defined a building $\sB(G,K)$ associated to $G$. For each point $x$ in $\sB(G,K)$, they constructed a bounded subgroup $G_x$ of $G(K)$, called a parahoric subgroup. In \cite{MP1, MP2}, Moy and Prasad defined a filtration of these parahoric subgroups by smaller subgroups
\begin{eqnarray*}
	G_x=G_{x,0} \triangleright G_{x, r_1} \triangleright G_{x, r_2} \triangleright \hdots  ,
\end{eqnarray*}
where $0< r_1 < r_2 < \hdots$ are real numbers depending on $x$. For simplicity, we assume that $r_1, r_2, \hdots$ are rational numbers. The quotient $G_{x,0}/G_{x,r_1}$ can be identified with the $\ov \bF_p$-points of a reductive group $\RP_x$, and $G_{x,r_i}/G_{x,r_{i+1}}$ ($i>0$) can be identified with an $\fpb$-vector space $\PV_{x,r_i}$ on which $\RP_x$ acts.

\textbf{Results about Moy--Prasad filtrations.}
We  show for a large class of reductive groups $G$, 
 which we call \good{} groups (see Definition \ref{Def-good}), that Moy--Prasad filtrations are in a certain sense (made precise below) independent of the residue field characteristic $p$. The class of \good{} groups contains reductive groups that split over a tamely ramified field extension (which is the class that many authors restrict to), as well as simply connected and adjoint semisimple groups, and products and restriction of scalars along finite separable (not necessarily tamely ramified) field extensions of any of these. The restriction to this (large) subclass of reductive groups is necessary as the main result (Theorem \ref{Thm-general}) fails in general, see Remark \ref{Remark-failure}. Given a \good{} reductive group $G$ over $K$, where $K$ is a maximal unramified extension of $k$ as above, a point $x$ of the Bruhat--Tits building $\sB(G,K)$ as above, and an arbitrary prime $q$ coprime to a certain integer $N$  that depends on the splitting field of $G$ (for details see Definition \ref{Def-good}), we construct a finite extension $K_q$ of $\bQ_q^{ur}$, a reductive group $G_q$ over $K_q$ and a point $x_q$ in $\sB(G_q,K_q)$. To these data, one can attach a Moy--Prasad filtration as above. The corresponding reductive quotient $\RP_{x_q}$ is a reductive group over $\ov \bF_q$ that acts on the quotients $\PV_{x_q,r_i}$, which are identified with $\ov \bF_q$-vector spaces. For a given positive integer $i$, we show in Theorem \ref{Thm-general} that there exists a split reductive group scheme $\ZG$ over $\ov \bZ[1/N]$ acting on a free $\ov \bZ[1/N]$-module $\ZV$ such that the special fiber of this representation over $\ov \bF_q$ is the above constructed Moy--Prasad filtration representations of $\RP_{x_q}$ on $\PV_{x_q, r_i}$ for all $q$ coprime to $N$, and the special fiber over $\ov \bF_p$ is the Moy--Prasad filtration representations of $\RP_{x}$ on $\PV_{x, r_i}$. This allows us to compare the Moy--Prasad filtration representations for different primes.

We also give a new description of the Moy--Prasad filtration representations, i.e. of $\RP_x$ acting on $\PV_{x,r_i}$, for reductive groups that split over a tamely ramified field extension of $K$. Let $m$ be the order of $x$ (see page \pageref{page-order} for the definition of ``order''). We define an action of the group scheme $\mmu_m$ of $m$-th roots of unity on a reductive group $\ZZG_{\fpb}$ over $\fpb$, and denote by $\ZZG_{\fpb}^{\mmu_m,0}$ the identity component of the fixed-point group scheme. In addition, we define a related action of $\mmu_m$ on the Lie algebra $\Lie(\ZZG_{\fpb})$, which yields a decomposition $\Lie(\ZZG_{\fpb})(\fpb)= \bigoplus_{i=1}^{m}\Lie(\ZZG_{\fpb})_i(\fpb)$. Then we prove that the action of $\RP_x$ on $\PV_{x,r_i}$ corresponds to the action of $\ZZG_{\fpb}^{\mmu_m,0}$ on one of the graded pieces $\Lie(\ZZG)_j(\fpb)$ of the Lie algebra of $\ZZG_{\ov \bF_p}$. This was previously known by \cite{ReederYu} for sufficiently large primes $p$, and representations of the latter kind have been studied by Vinberg \cite{Vinberg} in characteristic zero and generalized to positive characteristic coprime to $m$ by Levy \cite{Levy}. To be precise, in this paper we even prove a global version of the above mentioned result. See Theorem \ref{Thm-tame} for details. We also show that the same statement holds true for all \good{} reductive groups after base change of $\ZG$ and $\ZV$ to $\ov \bQ$, see Corollary \ref{Cor-Vinberg}.

Moreover, the global version of the Moy--Prasad filtration representations given by Theorem \ref{Thm-general} allows us to describe the representations occurring in the Moy--Prasad filtrations of good reductive groups explicitly in terms of Weyl modules, see Section \ref{section-reps} for precise formulas.

\textbf{An application to supercuspidal representations.} 
Suppose $G$ is defined over $k$.
In 1998, Adler (\cite{Adler}) used the Moy--Prasad filtrations to construct supercuspidal representations of $G(k)$, and Yu (\cite{Yu}) generalized his construction three years later, both assuming that $G$ splits over a tamely ramified field extension of $k$. If $p$ does not divide the order of the Weyl group of $G$, then Yu's construction yields all supercuspidal representations (\cites{Fi-exhaustion, Kim}). However, it is known that the construction does not give rise to all supercuspidal representations for small primes $p$.

In 2014, Reeder and Yu (\cite{ReederYu}) gave a new construction of supercuspidal representations of smallest positive depth, which they called epipelagic representations. A vector in the dual  $\check \PV_{x,r_1}=\left(G_{x,r_1}/G_{x,r_2}\right)^\vee$ of the first Moy--Prasad filtration quotient is called stable (in the sense of geometric invariant theory) if its orbit under $\RP_x$ is closed and its stabilizer in $\RP_x$ is finite. The only input for the new construction of supercuspidal representations in \cite{ReederYu} is such a stable vector. Assuming that $G$ is a semisimple group that splits over a tamely ramified field extension, Reeder and Yu gave a necessary and sufficient criterion for the existence of stable vectors for sufficiently large primes $p$. In \cite{Splitcase}, Romano and the author removed the assumption on the prime $p$ for absolutely simple split reductive groups $G$, which yielded new supercuspidal representations for split groups.

One application of our results on Moy--Prasad filtrations is a criterion for the existence of stable vectors for all primes $p$ for a much larger class of semisimple groups, see Corollary \ref{Cor-stable}. 
 As a consequence we obtain new supercuspidal representations for a class of non-split $p$-adic reductive groups, including non-tame groups.
Similarly, we prove in Theorem \ref{Thm-semistable} that the existence of semistable vectors is independent of the residue field characteristic. Semistable vectors play an important role when moving from epipelagic representations to representations of higher depth.

\textbf{Structure of the paper.} In Section, \ref{section-basics} we first recall the Moy--Prasad filtration of $G$, and then in Section \ref{Section-reductive-quotient} we introduce a Chevalley system for the reductive quotient that will be used for the construction of the reductive group scheme $\ZG$ that appears in Theorem \ref{Thm-general}. In Section \ref{section-iota}, we construct an inclusion of the Moy--Prasad filtration representation of $G$ into that of $G_F$ for a sufficiently large field extension $F$ of $K$ that will allow us to define the action of $\ZG$ on $\ZV$ in Theorem \ref{Thm-general}.
Afterwards, in Section \ref{section-globalMP}, we move from a previously fixed residue field characteristic $p$ to other residue field characteristics $q$. More precisely, we first introduce the notion of a \good{} group and define $K_q / \bQ_q^{ur}$, $G_q$ over $K_q$, and $x_q \in \sB(G_q,K_q)$. In Section \ref{section-global-MP-filtration}, we prove our first main theorem, Theorem \ref{Thm-general}.
Section \ref{section-Vinberg} is devoted to giving a different description of the Moy--Prasad filtration representations and their global version as generalized Vinberg--Levy representations (Theorem \ref{Thm-tame}).
In Section \ref{section-semistablestablevectors}, we use the results of the previous sections to show that the existence of (semi)stable vectors is independent of the residue characteristic. This leads to new supercuspidal representations.
We conclude the paper by giving a description of the Moy--Prasad filtration representations in term of Weyl modules in Section \ref{section-reps}.

\textbf{Conventions and notation.} If $M$ is a free module over some ring $A$, and if there is no danger of confusion, then we denote the associated scheme whose functor of points is $B \mapsto M \otimes_A B$ for any $A$-algebra $B$ also by $M$.  In addition, if $G$ and $T$ are schemes over a  scheme $S$, then we may abbreviate the base change $G \times_S T$ by $G_T$; and, if $T=\Spec A$ for some ring $A$, then we may also write $G_A$ instead of $G_T$.

When we talk about the identity component of a smooth group scheme $G$ of finite presentation, we mean the unique open subgroup scheme whose fibers are the connected components of the respective fibers of the original scheme that contains the identity. The identity component of $G$ will be denoted by $G^0$.
If $G$ is a group scheme defined over a ring $R$, then $\Lie(G)$ denotes the corresponding Lie algebra functor over $R$, and, if $f: G \ra H$ is a map between group schemes over $R$, then we write $\Lie(f)$ for the corresponding induced map $\Lie(G) \ra \Lie(H)$.

Throughout the paper, we require reductive groups to be connected. 

For each prime number $q$, we fix an algebraic closure $\ov \bQ_q$ of $\bQ_q$ and an algebraic closure $\ov{\bF_{q}((t))}$ of $\bF_{q}((t))$. All algebraic field extensions of $\bQ_q$ and $\bF_{q}((t))$ are assumed to be contained in $\ov \bQ_q$ and $\ov{\bF_{q}((t))}$, respectively. We then denote by $\bQ_q^{ur}$ the maximal unramified extension of $\bQ_q$ (inside $\ov \bQ_q$), and by $\bF_q((t))^{ur}$ the maximal unramified extension of $\bF_q((t))$. For any field extension $F$ of $\bQ_q$ (or of $\bF_q((t))$), we denote by $F^{\rm{tame}}$ its maximal tamely ramified field extension. Similarly, we fix an algebraic closure $\ov \bQ$ of $\bQ$, and we denote by $\ov \bZ$ the integral closure of $\bZ$ in $\ov \bQ$ and by $\ov \bZ_q$ the integral closure of $\bZ_q$ in $\ov \bQ_q$.

In addition, we will use the following notation throughout the paper: $p$ denotes a fixed prime number, $k$ is a nonarchimedean local field (of arbitrary characteristic)
with residual characteristic $p$, and $K$ \index{notation}{$K$} is the maximal unramified extension of $k$\index{notation}{$k$} contained in the fixed algebraic closure above. We write $\cO$\index{notation}{$\cO$} for the ring of integers of $K$, $\va:K \ra \bZ \cup \{\infty\}$\index{notation}{$\va$} for a valuation on $K$ with image $\bZ \cup \{\infty\}$, and $\varpi$\index{notation}{$\varpi$} for a uniformizer.
$G$ is a reductive group over $K$, and $E$ denotes a splitting field of $G$, i.e., $E$ is a minimal field extension of $K$ such that $G_E$ is split. Note that all reductive groups over $K$ are quasi-split and hence $E$ is unique up to conjugation. Let $e$ be the degree of $E$ over $K$, $\cO_E$ the ring of integers of $E$, and $\varpi_E$ a uniformizer of $E$. Without loss of generality, we assume that $\varpi$ is chosen to equal $\varpi_E^e$ modulo $\varpi_E^{e+1} \cO_E$. 
We denote the (absolute) root datum of $G$ by $R(G)$,\index{notation}{$R(G)$} and its root system by $\Phi=\Phi(G)$\index{notation}{$\Phi$}. We fix a point $x$ in the (reduced) Bruhat--Tits building $\sB(G,K)$ of $G$, denote by $S$\index{notation}{$S$} a maximal split torus of $G$ such that $x$ is contained in the apartment $\sA(S,K)$ associated to $S$, and let $T$\index{notation}{$T$} be the centralizer of $S$, which is a maximal torus of $G$. Moreover, we fix a Borel subgroup $B$ \index{notation}{$B$} of $G$ containing $T$, which yields a choice of simple roots $\Delta$ and positive roots $\Phi^+$ in $\Phi$. In addition, we denote by $\Phi_K=\Phi_K(G)$\index{notation}{$\Phi_K$} the restricted root system of $G$, i.e., the restrictions of the roots in $\Phi$ from $T$ to $S$. 
For $a \in \Phi_K$, we denote its preimage in $\Phi$ by $\Phi_a$.

Moreover, to help the reader, we will adhere to the convention of labeling roots in $\Phi$ by Greek letters: $\alpha, \beta, \hdots$, and roots in $\Phi_K$ by Latin letters: $a, b, \hdots$.

\textbf{Acknowledgment.} The author thanks her advisor 
Benedict Gross for his support. In addition, she thanks Mark Reeder for his lectures on epipelagic representations and inspiring discussions during the spring school and conference on representation theory and geometry of reductive groups.
She also thanks Jeffrey Adams, Jeffrey Adler, Stephen DeBacker, Brian Conrad, Wee Teck Gan, Thomas Haines, Tasho Kaletha, Gopal Prasad, Beth Romano, Cheng-Chiang Tsai and Zhiwei Yun for interesting discussions related to this paper.
Moreover, the author is particularly grateful to Jeffrey Hakim for comments on an initial draft, to Loren Spice and Siu-Fung Wong for carefully reading part of an earlier version of this paper, and to Stephen DeBacker for his feedback on several versions of the introduction.	
The author also thanks an anonymous referee for some comments and suggestions.

\section{Parahoric subgroups and Moy--Prasad filtration} \label{section-basics}

In order to talk about the Moy--Prasad filtration, we will first recall the structure of the root groups following \cite[Section~4]{BT2}. For more details and proofs we refer to \textit{loc.~cit.} 

\subsection{Chevalley--Steinberg system}

For $\alpha \in \Phi$, we denote by $\UE_\alpha$ \index{notation}{$\UE_\alpha$} the root subgroup of $G_E$ corresponding to $\alpha$. Note that $\Gal(E/K)$ acts on $\Phi$. We denote by $E_\alpha$ \index{notation}{$E_\alpha$} the fixed subfield of $E$ of the stabilizer $\Stab_{\Gal(E/K)}(\alpha)$ of $\alpha$ in $\Gal(E/K)$. In order to parameterize the root groups of $G$ over $K$, we fix a \textit{Chevalley--Steinberg system} $\{ \xE_\alpha: \bG_a \ra \UE_\alpha\}_{\alpha \in \Phi}$ of $G$ with respect to $T$, i.e. a Chevalley system \index{definition}{Chevalley--Steinberg system} $\{ \xE_\alpha: \bG_a \ra \UE_\alpha\}_{\alpha \in \Phi}$ \index{notation}{$\xE_\alpha$} of $G_E$ (see Remark \ref{Rem-Chevalley}) satisfying the following additional properties for all roots $\alpha \in \Phi$:
\begin{enumerate}[(i)]
	\item The isomorphism $\xE_\alpha:\bG_a \ra \UE_\alpha$ is defined over $E_\alpha$. 
	\item If the restriction $a \in \Phi_K$ of $\alpha$ to $S$ is not divisible, i.e. $\frac{a}{2} \notin \Phi_K$, then $\xE_{\gamma(\alpha)}=\gamma \circ \xE_\alpha \circ \gamma^{-1}$ for all $\gamma \in \Gal(E/K)$.
	\item  If the restriction $a \in \Phi_K$ of $\alpha$ to $S$ is divisible, then there exist $\beta, \beta' \in \Phi$ restricting to $\frac{a}{2}$ such that 
	 $E_\beta=E_{\beta'}$ is a quadratic extension of $E_\alpha$, and $\xE_{\gamma(\alpha)} =\gamma \circ \xE_\alpha \circ \gamma^{-1} \circ \eps$ for all $\gamma \in \Gal(E/E_\alpha)$, where $\eps \in \{\pm1 \}$ is 1 if and only if $\gamma$ induces the identity on $E_\beta$.
\end{enumerate} 

According to \cite[4.1.3]{BT2} such a Chevalley--Steinberg system does exist. It is a generalization of a Chevalley system to non-split groups and it will allow us to define a valuation of root groups in Section \ref{Sec-rootgroups} even if the group $G$ is non-split.
\begin{Rem} \label{Rem-Chevalley}
	We follow the conventions resulting from \cite[XXIII Définition~6.1]{SGA3III}, so we do not add the requirement of Bruhat and Tits that for each root $\alpha$, $\xE_\alpha$ and $\xE_{-\alpha}$ are associated, i.e. $\xE_\al(1)\xE_{-\alpha}(1)\xE_\alpha(1)$ is contained in the normalizer of $T$. However, there exists $\eps_{\al,\al} \in \{1, -1\}$ such that $$m_\al:=\xE_\al(1)\xE_{-\alpha}(\eps_{\al,\al})\xE_\alpha(1)$$ \index{notation}{$m_\al$}
	 is contained in the normalizer of $T$. Moreover, 
	  $\Ad(m_\al)(\Lie(\xE_\al)(1))=\eps_{\al,\al}\Lie(\xE_{-\al})(1)$.
\end{Rem}
	
\begin{Def} \label{Def-signs} For $\al, \beta \in \Phi$, we define $\eps_{\alpha,\beta} \in \{\pm 1\}$ \index{notation}{$\eps_{\al,\beta}$} by 
	$$ \Ad(m_\al)(\Lie(x_\beta)(1))=\eps_{\alpha,\beta}\Lie(x_{s_\alpha(\beta)})(1), $$ 
\end{Def}	
where $s_\al$\index{notation}{$s_\al$} denotes the reflection in the Weyl group $W$ of $\Phi(G)$ corresponding to $\al$.
The integers $\eps_{\al,\beta}$ for $\al$ and $\beta$ in $\Phi$ are called the \textit{signs} \index{definition}{signs of a Chevalley--Steinberg system} of the Chevalley--Steinberg system $\{x_\al^E\}_{\al \in \Phi}$.
	
\subsection{Parametrization and valuation of root groups}	\label{Sec-rootgroups}
In this section, we associate a parametrization and a valuation to each root group of $G$.

Let $a\in \Phi_K=\Phi_K(G)$, and let $U_a$ \index{notation}{$U_a$} be the corresponding root subgroup of $G$, i.e., the connected unipotent (closed) subgroup of $G$ normalized by $S$ whose Lie algebra is the sum of the root spaces corresponding to the roots that are a positive integral multiple of $a$. 

Let $G_a$ be the subgroup of $G$ generated by $U_a$ and $U_{-a}$, and let $\pi:G^a \ra G_a$ \index{notation}{$\pi$} \index{notation}{$G^a$} \index{notation}{$G_a$} be a simply connected cover. Note that $\pi$ induces an isomorphism between a %(negative)
 root group $U_+$ %($U_-$)
  of $G^a$ and $U_a$. We call $U_+$ the positive root group of $G^a$. 
   In order to describe 
   the root group $U_a$, we distinguish two cases.

\textbf{Case 1:} The root $a \in \Phi_K$ is neither divisible nor multipliable, i.e. $\frac{a}{2}$ and $2a$ are both not in $\Phi_K$.\\
Let $\alpha \in \Phi_a$ be a root that equals $a$ when restricted to $S$. Then $G^a$ is isomorphic to the Weil restriction $\Res_{E_\alpha/K}\SL_2$ of $\SL_2$ over $E_\alpha$ to $K$, and $U_a \simeq \Res_{E_\alpha/K}\UE_\alpha$, where $\UE_\alpha$ is the root group of $G_E$ corresponding to $\alpha$ as above. Note that $(U_a)_E$ is the product $\prod_{\beta \in \Phi_a}\UE_\beta$. Using the $E_\alpha$-isomorphism $\xE_\alpha:\bG_a \ra \UE_\alpha$, we obtain a $K$-isomorphism 
\begin{equation*}
	x_a:=\Res_{E_\alpha/K}\xE_\alpha: \Res_{E_\alpha/K}\bG_a \ra \Res_{E_\alpha/K}\UE_\alpha \xrightarrow{\simeq} U_a,
\end{equation*}
which we call a \textit{parametrization} \index{definition}{parametrization of $U_a$} of $U_a$. Note that for $u\in  \Res_{E_\alpha/K}\bG_a (K)=E_\alpha$, we have 
\begin{equation*}
	x_a(u)=\prod_{\beta \in \Phi_a}\xE_\beta(u_\beta), \text{ with } u_{\gamma(\alpha)}=\gamma(u) \index{notation}{$u_\beta$} \text{ for } \gamma \in \Gal(E/K).
\end{equation*}
This allows us to define the \textit{valuation} $\varphi_a:U_a(K) \ra \frac{1}{[E_\alpha:K]} \bZ \cup \{\infty \}$ \index{notation}{$\varphi_a$} of $U_a(K)$ \index{definition}{valuation of $U_a(K)$} by
\begin{equation*}
	\varphi_a(x_a(u))=\va(u).
\end{equation*}

\textbf{Case 2:} The root $a \in \Phi_K$ is divisible or multipliable, i.e. $\frac{a}{2}$ or $2a \in \Phi_K$.\\
We assume that $a$ is multipliable and describe $U_a$ and $U_{2a}$.
Let $\alpha, \wt \alpha \in \Phi_a$ \index{notation}{$\wt \alpha$} be such that $\alpha+\wt \alpha$ is a root in $\Phi$. Then $G^a$ is isomorphic to $\Res_{E_{\alpha+\wt \alpha}/K}\SU_3$, where $\SU_3$ \index{notation}{$\SU_3$} is the special unitary group over $E_{\alpha+\wt \alpha}$ defined by the hermitian form $(x,y,z) \mapsto \sigma(x) z+\sigma(y) y+\sigma(z) x$ on $E_\alpha^3$ with $\sigma$ the nontrivial element in $\Gal(E_\alpha/E_{\alpha+\wt \alpha})$. Hence, in order to parametrize $U_a$, we first parametrize the positive root group $U_+$ of $SU_3$. To simplify notation, write $L=E_{\alpha}=E_{\wt \alpha}$ \index{notation}{$L$} and $L_2=E_{\alpha + \wt \alpha}$. \index{notation}{$L_2$} Following \cite{BT2}, we define the subset $H_0(L,L_2)$ \index{notation}{$H_0(L,L_2)$} of $L\times L$ by
\begin{equation*}
	H_0(L,L_2)=\{(u,v) \in L \times L \, | \, v + \sigma(v)=\sigma(u)u \}.
\end{equation*}
Viewing $L \times L$ as a four dimensional vector space over $L_2$, and considering the corresponding scheme over $L_2$ (as described in ``Conventions and notation'' in Section \ref{Sec-Intro}), we can view $H_0(L, L_2)$ as a closed subscheme of $L \times L$ over $L_2$, which we will again denote by $H_0(L, L_2)$. Then there exists an $L_2$-isomorphism $\mu: H_0(L,L_2) \ra U_+$ \index{notation}{$\mu$} given by
\begin{equation*}
	(u,v) \mapsto \left( 
	\begin{matrix} 1 & -\sigma(u) & -v \\
					0 & 1 & u \\
					0 & 0 & 1 \end{matrix}	\right) , 
\end{equation*}
where $\sigma$ is induced by the nontrivial element in $\Gal(L/L_2)$. Using this isomorphism, we can transfer the group structure of $U_+$ to $H_0(L,L_2)$ and thereby turn the latter into a group scheme over $L_2$. Let us denote the  restriction of scalars $\Res_{L_2/K}H_0(L,L_2)$ of $H_0(L,L_2)$ from $E_{\alpha + \wt \alpha}=L_2$ to $K$ by $H(L, L_2)$\index{notation}{$H(L, L_2)$}. Then, by identifying $G^a$ with $\Res_{E_{\alpha+\wt \alpha}/K}\SU_3$, we obtain an isomorphism
\begin{equation*}
	x_a:=\pi \circ \Res_{E_{\alpha + \wt  \alpha}/K} \mu: H(L,L_2) \xra{\simeq} U_a,
\end{equation*}	
which we call the \textit{parametrization} \index{definition}{parametrization of $U_a$} of $U_a$.
We can describe the isomorphism $x_a$ on $K$-points as follows. Let $[\Phi_a]$ \index{notation}{$[\Phi_a]$} be a set of representatives in $\Phi_a$ of the orbits of the action of $\Gal(E_\alpha/E_{\alpha + \wt \alpha})=\<\sigma\>$ on $\Phi_a$. We will choose the sets of representatives for $\Phi_a$ and $\Phi_{-a}$ such that $[\Phi_a]$ and $-[\Phi_{-a}]$ are disjoint. For $\beta \in [\Phi_a]$, choose $\gamma \in \Gal(E/K)$ such that $\beta=\gamma(\alpha)$ and set $\wt \beta = \gamma(\wt \alpha)$ \index{notation}{$\wt \beta$} and $u_\beta=\gamma(u)$ \index{notation}{$u_\beta$} for every $u \in L$.
 By replacing some $\xE_{\beta + \wt \beta}$ by $\xE_{\beta + \wt \beta} \circ (-1)$ if necessary, we ensure that $\xE_{\beta + \wt \beta} = \Conjug( m_{\wt \beta}^{-1}) \circ \xE_\beta$ (where $m_{\wt \beta}$ is defined as in Remark \ref{Rem-Chevalley})\footnote{Note that our choice of $\xE_\beta$ or $\xE_{\beta+\wt \beta}$ for negative roots $\beta, \wt \beta$ deviates from Bruhat and Tits. It allows us a more uniform construction of the root group parameterizations that does not require us to distinguish between positive and negative roots, but that coincides with the ones defined by Bruhat and Tits in \cite{BT2}.}. Moreover, we choose the identification of $G^a$ with $\Res_{E_{\alpha+\wt \alpha}/K}\SU_3$ so that its restriction to the positive root group  arises from the restriction of scalars of the identification 
  that satisfies 
 \begin{equation*}
 	\pi \left( \left( 
 	\begin{matrix} 1 & -w & v \\
 		0 & 1 & u \\
 		0 & 0 & 1 \end{matrix}	\right) \right) = \xE_\al(u)\xE_{\al+\wt\al}(v)\xE_{\wt \al}(w). 
 \end{equation*}
 
Then we have for $(u,v)\in H_0(L,L_2)=H(L,L_2)(K) \subset L \times L$ that
\begin{equation} \label{equation-xa}
	x_a(u,v)=\prod_{\beta \in [\Phi_a]}	\xE_\beta(u_\beta)\xE_{\beta + \wt \beta}(-v_\beta)\xE_{\wt \beta}(\sigma(u)_\beta).
\end{equation}

The root group $U_{2a}$ corresponding to $2a$ is the subgroup of $U_{a}$ given by the image of $x_a(0,v)$. Hence $U_{2a}(K)$ is identified with the group of elements in $E_\alpha$ of trace zero with respect to the quadratic extension $E_\alpha/E_{\alpha + \wt \alpha}$, which we denote by $E_\alpha^0$. \index{notation}{$E_\alpha^0$}

Using the parametrization $x_a$, we define the \textit{valuation} \index{definition}{valuation of $U_a(K)$} $\varphi_a$ \index{notation}{$\varphi_a$} of $U_a(K)$ and $\varphi_{2a}$ \index{notation}{$\varphi_{2a}$} of $U_{2a}(K)$ by
\begin{eqnarray*}
	\varphi_a(x_a(u,v))&=&\frac{1}{2}\va(v) \\
	\varphi_{2a}(x_a(0,v))&=&\va(v) \, .
\end{eqnarray*}

\begin{Rem} \label{Rmk-val-u}
	\begin{enumerate}[(i)]
	\item Note that $v + \sigma(v) =\sigma(u)u$ implies that $\frac{1}{2}\va(v) \leq \va(u).$
	\item	The valuation of the root groups $U_a$ can alternatively be defined for all roots $a\in\Phi_K$ as follows. Let $u \in U_a(K)$, and write $u=\prod\limits_{\alpha \in \Phi_a \cup \Phi_{2a}}  u_\alpha$ with $u_\alpha \in U_\alpha(E)$. Then 
	\begin{equation*}
		\varphi_a(u)=\inf\left(\inf_{\alpha \in \Phi_a} \varphi^E_\alpha(u_\alpha), \inf_{\alpha \in \Phi_{2a}} \frac{1}{2} \varphi^E_\alpha(u_\alpha) \right),	
	\end{equation*}
	where $\varphi_\alpha^E(x_\alpha(v))=\va(v)$. The equivalence of the definitions is an easy exercise, see also \cite[4.2.2]{BT2}.
	\end{enumerate}
\end{Rem}

\subsection{Affine roots} \label{section-affine-roots}
Recall that the apartment $\sA=\sA(S,K)$ corresponding to the maximal split torus $S$ of $G$ is an affine space under the $\bR$-subspace of $X_*(S) \otimes_{\bZ} \bR$ spanned by the coroots of $G$, where $X_*(S)=\Hom_{K}(\bG_m,S)$.\index{notation}{$X_*(S)$} The apartment $\sA$ can be defined as corresponding to all valuations of $(T(K),$ $(U_a(K))_{a \in \Phi_K})$ in the sense of \cite[Section~6.2]{BT1} that are equipolent to the one constructed in Section \ref{Sec-rootgroups}, i.e., families of maps $(\wt \varphi_a: U_a(K) \ra \bR \cup \{\infty\})_{a \in \Phi_K}$ such that there exists $v \in X_*(S) \otimes_{\bZ} \bR$ satisfying $\wt \varphi_a(u)=\varphi_a(u)+a(v)$ for all $u \in U_a(K)$, for all $a \in \Phi_K$. In particular, the valuation defined in Section \ref{Sec-rootgroups} corresponds to a (special) point in $\sA$ that we denote by $x_0$. \index{notation}{$x_0$} Then the set of affine roots $\Psi_K$ \index{notation}{$\Psi_K$} on $\sA$ consists of the affine functions on $\sA$ given by
\begin{equation*}
	\Root = \Root(\sA)= \left\{y \mapsto a(y-x_0)+{\gamma'} \, | \, a \in \Phi_K, {\gamma'} \in \Gamma_a' \right\} ,
\end{equation*}
where \index{notation}{$\Gamma_a'$}
\begin{equation*}
	\Gamma_a' =\left\{ \varphi_a(u) \, | \, u \in U_a-\{1\}, \varphi_a(u)=\sup \varphi_a(uU_{2a}) \right\} .
\end{equation*} 
It will turn out to be handy to introduce a more explicit description of $\Gamma_a'$. In order to do so, consider a multipliable root $a$ and $\alpha \in \Phi_a$, and define \index{notation}{$(E_\alpha)^0$} \index{notation}{$(E_\alpha)^1$} \index{notation}{$(E_\alpha)^1_{\text{max}}$} 
\begin{eqnarray*}
	(E_\alpha)^0&=&\{u \in E_\alpha \, | \, \Tr_{E_\alpha/E_{\alpha+\wt \alpha}}(u)=0\},\\
	(E_\alpha)^1&=&\{u \in E_\alpha \, | \, \Tr_{E_\alpha/E_{\alpha+\wt \alpha}}(u)=1\},\\
	(E_\alpha)^1_{\text{max}}&=&\left\{u \in (E_\alpha)^1 \, | \, \va(u)= \sup\{\va(v) \, | \, v \in (E_\alpha)^1\}  \right\}.
\end{eqnarray*} 
Then, by \cite[4.2.20, 4.2.21]{BT2}, the set $(E_\alpha)^1_{\text{max}}$ is nonempty, and, with $\lambda$ any element of $(E_\alpha)^1_{\text{max}}$ and $a$ still being multipliable, we have
\begin{eqnarray} 
		\Gamma_a'&=&\tfrac{1}{2}\va(\lambda)+\va(E_\alpha-\{0\}) \label{equation-Ealpha-mult} \\
		\Gamma_{2a}'&=&\va((E_\alpha)^0-\{0\}) = \va(E_\alpha-\{0\}) -2 \cdot \Gamma_a'. \label{equation-Ealpha-div}
\end{eqnarray}	
For $a$ being neither multipliable nor divisible and $\alpha \in \Phi_a$, we have
\begin{equation} \label{equation-Ealpha-nonmult}
	\Gamma_a'=\va(E_\alpha-\{0\}).
\end{equation}
\begin{Rem} \label{Rem-lambda}
	Note that if the residue field characteristic $p$ is not 2, then $\frac{1}{2} \in (E_\alpha)^1_{\rm{max}}$ for $a$ a multipliable root and $\alpha \in \Phi_a $, and hence $\Gamma_a'=\va(E_\alpha-\{0\})$. If the residue field characteristic is $p=2$, then $\va(\lambda)<0$ for $\lambda \in (E_\alpha)^1_{\rm{max}}$.
\end{Rem}

\subsection{Moy--Prasad filtration} \label{section-Moy--Prasad}
Bruhat and Tits (\cite{BT1, BT2}) associated to each point $x$ in the (reduced) Bruhat--Tits building $\sB(G,K)$ \index{notation}{$\sB(G,K)$} a parahoric group scheme over $\cO$, which we denote by $\PP_x$\index{notation}{$\PP_x$}, whose generic fiber is isomorphic to $G$. We will quickly recall the filtration of $G_{x}:=\PP_x(\cO)$ \index{notation}{$G_x$} introduced by Moy and Prasad in \cite{MP1,MP2} and thereby specify our convention for the involved parameter.

Define $T_0=T(K)\cap\PP_x(\cO)$.\index{notation}{$T_0$} Then $T_0$ is a subgroup of finite index in the maximal bounded subgroup $\{ t\in T(K) \, | \, \va(\chi(t))=0 \, \forall \, \chi \in X^*(T)=\Hom_{\Kal}(T,\bG_m) \}$\index{notation}{$X^*(T)$} of $T(K)$.  
Note that this index equals one if $G$ is split.

For every positive real number $r$, we define 
\begin{equation*}
	T_r=\{t \in T_0 \, | \, \va(\chi(t)-1) \geq r  \text{ for all } \chi \in X^*(T)=\Hom_{\Kal}(T,\bG_m) \}. \index{notation}{$T_r$}
\end{equation*}

For every affine root $\psi \in \Root$, we denote by $\dot \psi$ \index{notation}{$\dot \psi$} its gradient and define the subgroup $U_\psi$ of $U_{\dot \psi}(K)$ by
\begin{equation*}
	U_\psi=\{u \in U_{\dot \psi}(K) \, | \, u=1 \, \text{ or } \, \varphi_{\dot \psi}(u)\geq \psi(x_0)\}. \index{notation}{$U_\psi$}
\end{equation*}
Then the Moy--Prasad filtration subgroups of $G_x$ are given by
\begin{equation*}
	G_{x,r}=\<T_r, U_\psi \, | \, \psi \in \Root, \psi(x) \geq r \> \, \text{ for } \, r \geq 0, \index{notation}{$G_{x,r}$}
\end{equation*}
and we set 
\begin{equation*}
	G_{x,r+}=\bigcup_{s>r} G_{x,s}. \index{notation}{$G_{x,r+}$}
\end{equation*}
The quotient $G_x/G_{x,0+}$ can be identified with the $\fpb$-points of the reductive quotient of the special fiber $\PP_x \times_\cO \fpb$ of the parahoric group scheme $\PP_x$, which we denote by $\RP_x$.\index{notation}{$\RP_x$} From \cite[Corollaire~4.6.12]{BT2} we deduce the following lemma.
\begin{Lemma}[\cite{BT2}] \label{lemma-reductive-quotient}
	Let $R_K(G)=(X_K=X^*(S),\Phi_K,\check X_K=X_*(S), \check \Phi_K)$ be the restricted root datum of $G$. Then the root datum $R(\RP_x)$ of $\RP_x$ is canonically identified with $(X_K,\Phi',\check X_K, \check \Phi')$ where 
	\begin{equation*}
		\Phi'
		=\{a \in \Phi \, | \,  a(x-x_0) \in \Gamma_a'\} \quad \text{ and } \quad		\check\Phi'
			=\{\check a \in \check\Phi \, | \,  a(x-x_0) \in \Gamma_a'\} .
	\end{equation*}
\end{Lemma}

We can define a filtration of the Lie algebra $\fg=\Lie(G)(K)$ similar to the filtration of $G_x$. In order to do so, we denote the $\cO$-lattice $\Lie(\PP_x)$ of $\fg$ by $\fp$. Define $\fp_a = \fp \cap \fg_a$ for $a \in \Phi_K$ and $\ft=\Lie(T)(K)$, where $\fg_a$ \index{notation}{$\fg_a$} is the subspace of $\fg$ on which $\ft$ acts via $\Lie(a)$.

We define the Moy--Prasad filtration of the Lie algebra $\ft$ for $r \in \bR$ to be
\begin{equation}
	\ft_r = \left\{ X \in \ft \, | \, \va(\Lie(\chi)(X)) \geq r \text{ for all } \chi \in X^*(T) \right\} \index{notation}{$\ft_r$}
\end{equation}

For every root $a \in \Phi_K$, we define the Moy--Prasad filtration of $\fg_a$ as follows. Let $\psi_a$ be the smallest affine root with gradient $a$ such that $\psi_a(x)\geq 0$. For every $\psi \in \Root$ with gradient $a$, we let $n_{\psi}=e_{\alpha}(\psi-\psi_a)$, where $e_\alpha=[E_\alpha:K]$ \index{notation}{$\varpi_\alpha$} for some root $\alpha \in \Phi_a$ that restricts to $a$. Note that $n_\psi$ is an integer. Choosing a uniformizer $\varpi_\alpha \in E_\alpha$ \index{notation}{$\varpi_\alpha$} and viewing $\fp_a$ inside $\Lie(G)(E_\al)$ we set\footnote{Note that $\fu_{\psi}$ does not depend on the choice of $x$ inside $\sA$.}
\begin{equation*} 
	\fu_{\psi}=(\varpi_{\alpha}^{n_\psi}\cO_{E_\al}\fp_a) \cap \fg. \index{notation}{$\fu_\psi$}
\end{equation*}

Then the Moy--Prasad filtration of the Lie algebra $\fg$ is given by	
\begin{equation*}
	\fg_{x,r}=\<\ft_r, \fu_\psi \, | \, \psi(x) \geq r \> \text{ for } r \in \bR. \index{notation}{$\fg_{x,r}$}
\end{equation*}

In general, the quotient $G_{x,r}/G_{x,r+}$ is not isomorphic to $\fg_{x,r}/\fg_{x,r+}$ for $r>0$. However, it turns out that we can identify them (as $\fpb$-vector spaces) under the following assumption. 

\begin{Assumption} \label{assumption}
	The maximal (maximally split) torus $T$ of $G$ becomes an induced torus over a tamely ramified extension.
\end{Assumption}

Recall that the torus $T$ is called \textit{induced} \index{definition}{induced torus} if it is a product of separable Weil restrictions of $\bG_m$, i.e. $T\simeq\prod\limits_{i=1}^N\Res_{K_i/K} \bG_m$ for some integer $N$ and finite separable field extensions $K_i/K$, $1 \leq i \leq N$.

For the rest of Section \ref{section-basics}, we impose Assumption \ref{assumption}. 
\begin{Rem}
	Assumption \ref{assumption} holds, for example, if $G$ is either adjoint or simply connected semisimple, or if $G$ splits over a tamely ramified extension. 
\end{Rem}

For $r \in \bR$, we denote the quotient $\fg_{x,r}/\fg_{x,r+}$ ($\simeq G_{x,r}/G_{x,r+}$ for $r>0$) by \index{notation}{$\PV_{x,r}$} $\PV_{x,r}$. The adjoint action of $G_{x,0}$ on $\fg_{x,r}$ (or, equivalently, the conjugation action of $G_{x,0}$ on $G_{x,r}$ for $r>0$) induces an action of the algebraic group $\RP_x$ on the quotients $\PV_{x,r}$.

\subsection{Chevalley system for the reductive quotient} \label{Section-reductive-quotient}
In this section we construct a Chevalley system for the reductive quotient $\RP_x$ by reduction of the root group parameterizations given in Section \ref{Sec-rootgroups}. Let $\RU_a$ \index{notation}{$\RU_a$} denote the root group of $\RP_x$ corresponding to the root $a \in \Phi(\RP_x) \subset \Phi_K(G)$. We denote by $\cO_{\bQ_p^{ur}}$ the ring of integers in ${\bQ_p^{ur}}$.
If $K$ is an extension of $\bQ_p^{ur}$, we let $\chi: \fpb \ra \cO_{\bQ_p^{ur}}$ \index{notation}{$\chi$} be the Teichmüller lift, i.e. the unique multiplicative section of the surjection $\cO_{\bQ_p^{ur}} \twoheadrightarrow \fpb$. If $K$ is an extension of $\bF_p((t))^{ur}=\varinjlim_{n \in \bN} \bF_{p^n}((t))$, we let $\chi: \fpb=\varinjlim_{n \in \bN} \bF_{p^n} \ra \varinjlim_{n \in \bN} \bF_{p^n}[[t]]$ be the usual inclusion.

\begin{Lemma} \label{Lemma-Chevalley} 
	Let $\lambda=\lambda_a \in   (E_{\alpha})^1_{\text{max}}$ for some $\alpha \in \Phi_a$, and write $\lambda = \lambda_0 \cdot \varpi_{E}^{\va(\lambda)e}\cdot \eps_0$ with $\lambda_0 \in \chi(\ov \bF_p)$ and $\eps_0 \in 1 + \varpi_E\cO_{E}$;\index{notation}{$\lambda_0$} e.g., take $\lambda_0\eps_0=\lambda=\frac{1}{2}$ if $p \neq 2$. Consider the map 
	\begin{eqnarray*} 
		\fpb & \ra & G_{x,0} \\
		u   & \mapsto  & \left\{ \begin{array}{ll} 
			{x_a\left(\sqrt{\frac1{\lambda_0}}\chi(u) \varpi_E^{s} \eps_1, \, \chi(u)\varpi_E^s\eps_1\sigma(\chi(u)\varpi_E^s\eps_1)\cdot \varpi_E^{\va(\lambda)e}\eps_0\right)} & \mbox{if $a$ is multipliable} \\ 
			x_a(0,\chi(u) \cdot \varpi_E^{-2a(x-x_0)\cdot e}\eps_2) & \mbox{if $a$ is divisible} \\
			{x_a(\chi(u) \cdot \varpi_E^{-a(x-x_0)\cdot e}\eps_3)} & \mbox{otherwise} \, ,
		\end{array} \right. 
	\end{eqnarray*}
	where $s=-(a(x-x_0)+\va(\lambda)/2)\cdot e$, and $\eps_1, \eps_2,\eps_3 \in 1+\varpi_E \cO_E$ such that $\sqrt{\frac1{\lambda_0}}\chi(u) \varpi_E^{s} \eps_1$, 
	 $\chi(u) \varpi_E^{-2a(x-x_0)\cdot e}\eps_2$ and $\chi(u) \varpi_E^{-a(x-x_0)\cdot e}\eps_3$ are contained in $E_\al$, and $\sqrt{\frac{1}{\lambda_0}} \in \chi(\ov \bF_p)$ with $\sqrt{\frac{1}{\lambda_0}}^2=\frac{1}{\lambda_0}$.

	Then the composition of this map with the quotient map $G_{x,0} \twoheadrightarrow G_{x,0}/G_{x,0+}$ yields a root group parameterization $\Rx_a: \bG_a \ra \RU_a \subset \RP_x$. \index{notation}{$\Rx_a:$}
	
	Moreover, the root group parameterizations $\{ \Rx_a\}_{a \in \Phi(\RP_x)}$ form a Chevalley system for $\RP_x$. \todor{chevalley system for ... add torus here}
\end{Lemma}

We remark that Gopal Prasad pointed out to us that a similar Chevalley system construction can be found in \cite[2.19,~2.20]{PrasadRaghunathan}. 

\textbf{Proof.} Note first that since $a \in \Phi(\RP_x)$, we have $a(x-x_0)\in  \Gamma'_a$ by Lemma \ref{lemma-reductive-quotient}. 
Suppose $a$ is multipliable.
Then $\RU_a(\fpb)$ is the image of 
$$\Image:= \left\{ x_a(U,V) \, | \, (U,V) \in H_0(E_\alpha,E_{\alpha+\wt \alpha}), \frac{1}{2}\va(V) = -a(x-x_0)\right\}.$$
in $ G_{x,0}/G_{x,0+}$. Set 
$$U(u)=\sqrt{\tfrac1{\lambda_0}}\chi(u) \cdot \varpi_E^{-(a(x-x_0)+\va(\lambda)/2)\cdot e} \eps_1$$ and $$V(u)=\chi(u)\varpi_E^s\eps_1\sigma(\chi(u)\varpi_E^s\eps_1)\cdot \varpi_E^{\va(\lambda)e}\eps_0.$$
Then $V(u)+\sigma(V(u))=U(u)\sigma(U(u))$, i.e. $(U(u),V(u))$ is in $H_0(E_\alpha, E_{\alpha+\wt \alpha})$, and $\va(V(u))= -2a(x-x_0)$. Moreover, every element in $\Image$ is of the form $(U(u),V(u)+v_0)$ for $u \in \fpb$ and some element $v_0 \in (E_{\alpha})^0$ with $\va(v_0)> -2a(x-x_0)$, because $2a(x-x) \notin \va((E_{\alpha})^0)$ (by Equation (\ref{equation-Ealpha-div}), page \pageref{equation-Ealpha-div}). Note that the images of $x_a(U(u),V(u)+v_0)$ and $x_a(U(u),V(u))$ in $ G_{x,0}/G_{x,0+}$ agree. Thus, by the definition of $x_a$, we obtain an isomorphism of group schemes $\Rx_a: \bG_a \ra \RU_a$. Similarly, one can check that $\Rx_a$ yields an isomorphism $\bG_a \ra \RU_a$ for $a$ not multipliable.

In order to show that $\{ \Rx_a\}_{a \in \Phi(\RP_x)}$ is a Chevalley system, suppose for the moment that $a$ and $b$ in $\Phi(\RP_x)$ are neither multipliable nor divisible, and $\Phi_a=\{\alpha\}$ and $\Phi_b=\{\beta\}$ each contain only one root. Let $\check\al$\index{notation}{$\check\al$} be the coroot of the root $\al$, and denote by $s_\al$\index{notation}{$s_\al$} the reflection in the Weyl group $W$ of $G$ corresponding to $\al$. Then, using \cite[Cor.~5.1.9.2]{ConradSGA3}, we obtain 
\begin{eqnarray*}
	& &\Ad\left(\xE_\alpha(\varpi_E^{-\alpha(x-x_0)e})\xE_{-\alpha}(\eps_{\alpha,\alpha}\varpi_E^{-(-\alpha)(x-x_0)e})\xE_\alpha(\varpi_E^{-\alpha(x-x_0)e})\right)\left(\Lie(\xE_\beta)(\varpi_E^{-\beta(x-x_0)e})\right) \\
	&=&\Ad\left(\check \alpha(\varpi_E^{-\alpha(x-x_0)e})\right)\Ad\left(\xE_\alpha(1)\xE_{-\alpha}(\eps_{\alpha,\alpha})\xE_\alpha(1)\right)\left(\varpi_E^{-\beta(x-x_0)e}\Lie(\xE_\beta)(1)\right) \\
	&=&\Ad\left(\check \alpha(\varpi_E^{-\alpha(x-x_0)e})\right)\left(\eps_{\alpha,\beta}\varpi_E^{-\beta(x-x_0)e}\Lie(\xE_{s_\alpha(\beta)})(1)\right)\\
	&=&(s_\alpha(\beta))(\check \alpha(\varpi_E^{-\alpha(x-x_0)e}))\eps_{\alpha,\beta}\varpi_E^{-\beta(x-x_0)e}\Lie(\xE_{s_\alpha(\beta)})(1)\\
	&=&\varpi_E^{\<\check \alpha,s_\alpha(\beta)\>(-\alpha(x-x_0))e}\eps_{\alpha,\beta}\varpi_E^{-\beta(x-x_0)e}\Lie(\xE_{s_\alpha(\beta)})(1)\\
	&=&\varpi_E^{\<\check \alpha,\beta\>\alpha(x-x_0)e-\beta(x-x_0)e}\eps_{\alpha,\beta}\Lie(\xE_{s_\alpha(\beta)})(1)\\
	&=&\eps_{\alpha,\beta}\Lie(\xE_{s_\alpha(\beta)})(\varpi_E^{-(s_\alpha(\beta))(x-x_0)e}) .
\end{eqnarray*}
This implies (assuming $\eps_3=1$, otherwise it is an easy exercise to add in the required constants) that for $\ov m_a:=\ov x_a(1) \ov x_{-a}(\eps_{a,a}) \ov x_a(1)$ with $\eps_{a,a}=\eps_{\alpha,\alpha}$ we have
\begin{equation*}
	\Ad(\ov m_a)(\Lie(\Rx_b)(1))=\Ad(\ov x_a(1) \ov x_{-a}(\eps_{a,a}) \ov x_a(1))(\Lie(\Rx_b)(1)) = \eps_{\alpha,\beta} \Lie(\Rx_{s_a(b)})(1).
\end{equation*}
We obtain a similar result even if $\Phi_a$ and $\Phi_b$ are not singletons by the requirement that $\{\xE_\alpha\}_{\alpha \in \Phi}$ is a Chevalley--Steinberg system, i.e. compatible with the Galois action as described in Section \ref{section-basics}. Similarly, we can extend the result that $\Ad(\ov m_a)(\Lie(\Rx_b)(1))=\pm \Lie(\Rx_{s_a(b)})(1)$ to all non-multipliable roots $a, b \in \Phi(\RP_x) \subset \Phi_K$.

Suppose now that $a \in \Phi(\RP_x) \subset \Phi_K$ is multipliable, and let $\alpha \in \Phi_a$ and $\wt \alpha=\sigma(\alpha) \in \Phi_a$ as above. 
Following \cite[4.1.11]{BT2}, we define for $(u,v)\in H_0(E_{\alpha},E_{\alpha+\wt\alpha})$
\begin{equation*} m_a(U,V) = x_a(UV^{-1},\sigma(V^{-1}))x_{-a}(\eps_{\alpha,\alpha}U,V)x_a(U\sigma(V^{-1}),\sigma(V^{-1})).
\end{equation*}
Then Bruhat and Tits  show in \textit{loc. cit.} that $m_a(U,V)$ is in the normalizer of the maximal torus $T$ and
\begin{equation} \label{eqn-ma} m_a(U,V) =m_{a,1} \wt a(V) \quad \text{and} \quad x_{-a}(\eps_{\alpha,\alpha}U,V)=m_{a,1}x_a(U,V)m_{a,1}^{-1},
\end{equation}
where 
\begin{multline} \label{eqn-ma2} m_{a,1} = \pi \circ \Res_{E_{\alpha+\wt \alpha/K}} 
\begin{pmatrix} 0 & 0 & -1 \\ 
				0 & -1 & 0\\
				-1 & 0 & 0 
 \end{pmatrix}
 \, \text{ and }
 \,
 \wt a(V) = \pi \circ \Res_{E_{\alpha+\wt \alpha/K}} 
 \begin{pmatrix} V & 0 & 0 \\ 
 0 & V^{-1}\sigma(V) & 0\\
 0 & 0 & \sigma(V^{-1}) 
 \end{pmatrix} .
\end{multline} 
Note that we have $$m_a(\sqrt{\tfrac1{\lambda_0}}(-\varpi_E)^{(a(x-x_0)-\va(\lambda)/2)e}\eps_1,\varpi_E^{(a(x-x_0)-\va(\lambda)/2)e}\eps_1\sigma(\varpi_E^{(a(x-x_0)-\va(\lambda)/2)e}\eps_1) \varpi_E^{\va(\lambda)e}\eps_0) \in G_{x,0},$$
 and denote its image in $G_{x,0}/G_{x,0+}$ by $\ov m_a$. Using that $\va(\lambda)=0$ if $p \neq 2$, and $\sigma(\varpi_E^{(a(x-x_0)-\va(\lambda)/2)e}\eps_1)\equiv \pm \varpi_E^{(a(x-x_0)-\va(\lambda)/2)e}\eps_1 \equiv\varpi_E^{(a(x-x_0)-\va(\lambda)/2)e}\eps_1 \mod \varpi_E^{(a(x-x_0)-\va(\lambda)/2)e+1}$ if $p=2$ as well as the compatibility with Galois action properties of a Chevalley--Steinberg system, we obtain 
$$ \ov  m_a= \Rx_{a}(1) \Rx_{-a}(\eps_{a,a})\Rx_a(1) \quad \text{ with } \quad \eps_{a,a}=\eps_{\alpha,\alpha}(-1)^{(a(x-x_0)-\va(\lambda)/2)e}.$$
Moreover, using Equation \eqref{eqn-ma} and \eqref{eqn-ma2}, an easy calculation shows that 
$$ \Rx_{-a}(\eps_{a,a}u)=\ov m_a x_a(u) \ov m_a^{-1}$$
for all $u \in \fpb$. In other words,
$$ \Ad(\ov m_a)(\Lie(\Rx_a)(1))=\eps_{a,a}\Lie(\Rx_{-a})(1),$$
as desired. We obtain analogous results for $\ov m_{-a}$ being defined as above by substituting ``$a$'' by ``$-a$''. Moreover, $\ov m_a=\ov m_{-a}$, and hence $ \Ad(\ov m_{-a})(\Lie(\Rx_a)(1))=\eps_{a,a}\Lie(\Rx_{-a})(1)$.

\begin{figure}[hbt]
	\centering
	\includegraphics[width=0.450\textwidth]{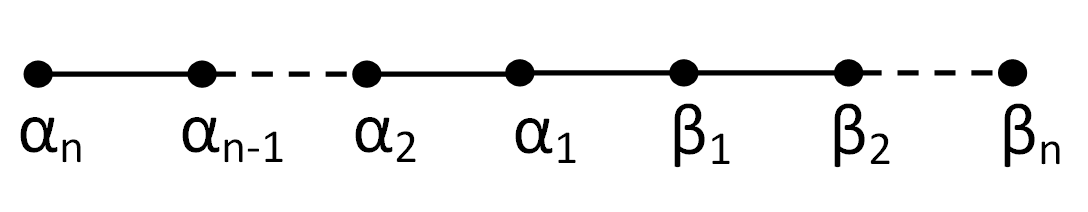}
	\caption{Dynkin diagram of type $A_{2n}$}
	\label{fig-A2n1}
\end{figure}

In order to show that $\{ \Rx_a \}_{a \in \Phi(\RP_x)}$ forms a Chevalley system, it is left to check that 
\begin{equation}\label{eqn-todo}
	 \Ad(\ov m_a)(\Lie(\Rx_b)(1))=\pm \Lie(\Rx_{s_a(b)})(1)
\end{equation}
holds for $a, b \in \Phi(\RP_x)$ with $a \neq \pm b$ and either $a$ or $b$ multipliable. Note that if $x_a$ and $x_{-a}$ commute with $x_b$, then the statement is trivial. 
Note also that if $b$ is multipliable and $\beta \in \Phi_b$, then $\beta$ lies in the span of the roots of a connected component of the Dynkin digram $\Dyn(G)$ of $\Phi(G)$ of type $A_{2n}$ for some positive integer $n$. Hence, for some $\al \in \Phi_a$, $\al$ and $\beta$ lie in the span of the roots of such a connected component. Moreover, by the compatibility of the Chevalley--Steinberg system $\{\xE_\alpha\}_{\alpha \in \Phi}$ with the Galois action, it suffices to restrict to the case where $\Dyn(G)$ is of type $A_{2n}$ with simple roots labeled by $\alpha_{n}, \alpha_{n-1}, \hdots, \alpha_1, \beta_1, \beta_2, \hdots, \beta_{n}$ as in Figure \ref{fig-A2n1},
 and the $K$-structure of $G$ arises from the unique outer automorphism of $A_{2n}$ of order two that sends $\alpha_i$ to $\beta_i$.
If a root in $\Phi_K(G)$ is multipliable, then it is the image of $\pm(\alpha_1 + \hdots + \alpha_s)$ in $\Phi_K$ for some $1 \leq s \leq n$.  In particular, the positive multipliable roots are orthogonal to each other, by which we mean that $\<\check a, b\>=0$ for two distinct positive multipliable roots $a$ and $b$.  Equation \eqref{eqn-todo} can now be verified by simple matrix calculations in $\SL_{2n+1}$. 
\qed

\subsection{Moy--Prasad filtration and field extensions} \label{section-iota}
Let $F$ be a field extension of $K$ of degree $d=[F:K]$ with ring of integers $\cO_F$, and denote by $\va:F \ra \frac{1}{d}\bZ \cup \{ \infty\}$ the extension of the valuation $\va:K \ra \bZ \cup\{\infty\}$ on $K$. Then there exists a $G(K)$-equivariant injection of the Bruhat--Tits building $\sB(G,K)$ of $G$ over $K$ into the Bruhat--Tits building $\sB(G_{F},F)$ of $G_F=G \times_K F$ over $F$. We denote the image of the point $x \in \sB(G,K)$ in $\sB(G_{F},F)$ by $x$ as well. Using the definitions introduced in Section \ref{section-Moy--Prasad}, but for notational convenience still with the valuation $\va$ (instead of replacing it by the normalized valuation $d \cdot \va$), we can define a Moy--Prasad filtration of $G(F)$ and $\fg_F$ at $x$, which we denote by $\GF_{x,r} (r \geq 0)$ \index{notation}{$\GF_{x,r}$} and $\fgF_{x,r} (r \in \bR)$\index{notation}{$\fgF_{x,r}$}, as well as its quotients $\PVF_{x,r} (r \in \bR)$ \index{notation}{$\PVF_{x,r}$} and the reductive quotient $\RPF_x$\index{notation}{$\RPF_x$}. 

Suppose now that $G_F$ is split, and that $\Gamma_a' \subset \va(F)$ for all restricted roots $a \in \Phi_K(G)$. This holds, for example, if $F$ is an even-degree extension of the splitting field $E$. Then, using Remark \ref{Rmk-val-u}(i) and the definition of the Moy--Prasad filtration, the inclusion $G(K) \hookrightarrow G(F)$ maps $G_{x,r}$ into $\GF_{x,r}$. Furthermore, recalling that for split tori $\wt T$ the subgroup $\wt T_0$ is the maximal bounded subgroup of the (rational points of) $\wt T$ and using the assumption that  $\Gamma_a' \subset \va(F)$ for all restricted roots $a \in \Phi_K(G)$, we observe that this map induces an injection
\begin{eqnarray} 
 \iota_{K,F}: & G_{x,0}/G_{x,0+} \hookrightarrow \GF_{x,0}/\GF_{x,0+} & , \label{equation-iota} 
\end{eqnarray}
which yields a map of algebraic groups $\RP_x \ra \RPF_x$, also denoted by $\iota_{K,F}$\index{notation}{$\iota_{K,F}$}. If $p \neq 2$ or $d$ is odd, then $\iota_{K,F}$ is a closed immersion.
 
To discuss a similar result for the higher-depth quotients, we denote by $\Phi_K^{\rm{mul}}$ \index{notation}{$\Phi_K^{\rm{mul}}$} the set of multipliable roots in $\Phi_K$ and by $\Phi_K^{\rm{nm}}$ \index{notation}{$\Phi_K^{\rm{nm}}$} the set of non-multipliable roots in $\Phi_K$.

\begin{Lemma} \label{Lemma-iota}
	For every $r \in \bR$, there exists an injection 
	\begin{equation*} \index{notation}{$\iota_{K,F,r}$}
		\iota_{K,F,r}: \PV_{x,r}=\fg_{x,r}/\fg_{x,r+} \hookrightarrow \fgF_{x,r}/\fgF_{x,r+}=\PVF_{x,r} 
	\end{equation*}
such that $\iota_{K,F}(\RP_{x})$ preserves $\iota_{K,F,r}(\PV_{x,r})$ under the action described in Section \ref{section-Moy--Prasad}. Moreover, we obtain a commutative diagram 
	\begin{eqnarray} \label{diagram-iota}   
		\xymatrix{
			\RP_x \times  \PV_{x,r} \ar[r] \ar[d]^{\iota_{K,F} \times \iota_{K,F,r}} &   \PV_{x,r} \ar[d]^{\iota_{K,F,r}} \\
			\RPF_x \times  \PVF_{x,r} \  \ar[r]    &  \PVF_{x,r}  \quad 
		} 
	\end{eqnarray} 
	unless $p=2$ and there exists $a \in \Phi_K^{\rm{mul}}$ with $a(x-x_0) \in \Gamma_a'$  such that $a(x-x_0)-r \in \Gamma_a'$ or such that there exists $b \in \Phi_K^{\rm{nm}}$ with $b(x-x_0)-r \in \Gamma_{b}'$ and $\<\check a, b\> \neq 0$.
\end{Lemma}

\textbf{Proof.}
For $p\neq 2$, let $\iota_{K,F,r}$ be induced by the inclusion $\fg \hookrightarrow \fg_F=\fg \otimes_K F$. This map is well defined, and it is easy to see that it is injective on $(\ft \cap \fg_{x,r})/\fg_{x,r+}$ and on $(\fg_a \cap \fg_{x,r})/\fg_{x,r+}$ for $a \in \Phi_K$ non-multipliable. Suppose $a$ is multipliable. If $r-a(x-x_0) \in \Gamma_a'$, i.e. there exists an affine root $\psi:y \mapsto a(y-x_0) + {\gamma'}$ with $\psi(x)=r$, and $\varphi_a(x_a(u,v))=\psi(x_0)=r-a(x-x_0) \in \Gamma_a'$, then $\va(u)=\frac{1}{2}\va(v)=r-a(x-x_0)$. This follows from the trace of $\frac{1}{2}$ being one, hence $v-\frac{1}{2}\sigma(u)u$ is traceless and therefore has valuation outside $2 \Gamma_a'$, while $\va(v) \in 2 \Gamma_a'$. Hence the image of $\fg_a \cap \fg_{x,r}$ in $\PVF_{x,r}$ is non-vanishing if it is non-trivial in $\PV_{x,r}$, i.e. if $r-a(x-x_0) \in \Gamma_a'$. Moreover, Diagram \eqref{diagram-iota} commutes.

In the case $p=2$, if $a \in \Phi_K$ is multipliable and $r-a(x-x_0) \in \Gamma_a'$ and $\varphi_a(x_a(u,v))=r-a(x-x_0)$, then $\va(u)=r-a(x-x_0)-\frac{1}{2}\va(\lambda_\al)$ for $\lambda_\al \in (E_\alpha)^1_{\rm{max}}$ by reasoning analogous to that above. However, recall from Remark \ref{Rem-lambda} that $\va(\lambda_\al)<0$ for $p=2$. Let $\varpi_{F}$ be a uniformizer of $F$ such that $\varpi_F^{\frac{d}{e}} =\varpi_F^{[F:E]}\equiv \varpi_E \mod \varpi_F^{[F:E]+1}$ and let $\varpi_\alpha$ be a uniformizer of $E_\alpha$ with $\varpi_\alpha \equiv \varpi_F^{[F:E_\alpha]}=\varpi_F^{\frac{d}{e_\al}} \mod \varpi_F^{[F:E_\alpha]+1}$. This allows us to define $\iota_{K,F,r}$ as follows. We define the linear morphism $i_{K,F,r}: \fg \hookrightarrow \fg_F$ to be the usual inclusion $\fg \hookrightarrow \fg_F=\fg \otimes_K F$ on $\ft \oplus \bigoplus\limits_{a \in \Phi_K^{\rm{nm}}} \fg_a$ and to be the linear map from $\bigoplus\limits_{a \in \Phi_K^{\rm{mul}}} \fg_a$ onto $\left(\bigoplus\limits_{a \in \Phi_K^{\rm{mul}}} \fg_a \otimes_{K} \varpi_F^{d \va(\lambda_\al)/2} K \right) \subset \fg_F$ on $\bigoplus\limits_{a \in \Phi_K^{\rm{mul}}} \fg_a$ such that 
$$i_{K,F,r}\left(\Lie(x_a)(\varpi_\al^{(r-a(x-x_0)-\va(\lambda_\al)/2)e_\alpha},0)\right)=\Lie(x_a)\left(\varpi_\al^{(r-a(x-x_0)-\va(\lambda_\al)/2)e_\alpha} \otimes \varpi_F^{d \va(\lambda_\al)/2},0\right),$$ 
where $\al \in \Phi_a$ for $a \in \Phi^{\rm{mul}}_K$.  
By restricting $i_{K,F,r}$ to $\fg_{x,r}$ and passing to the quotient, we obtain an injection $\iota_{K,F,r}$ of $\PV_{x,r}$ into $\PVF_{x,r}$.

In order to prove that $\iota_{K,F}(\RP_{x})$ preserves $\iota_{K,F,r}(\PV_{x,r})$ for $p=2$, it suffices to show that $\iota_{K,F}(\RP_{x})$ stabilizes the subspace $$V'=\iota_{K,F,r}\Bigl(\overline{\fg_{x,r} \cap \bigoplus\limits_{a \in \Phi_K^{\rm{mul}}} \fg_a}\Bigr),$$ where the overline denotes the image in $\PV_{x,r}$.

First suppose that the Dynkin diagram $\Dyn(G)$ of $\Phi(G)$ is of type $A_{2n}$ with simple roots labeled by $\alpha_{n}, \alpha_{n-1}, \hdots, \al_2, \alpha_1, \beta_1, \beta_2, \hdots, \beta_{n}$ as in Figure \ref{fig-A2n1} on page \pageref{fig-A2n1},
  and that the $K$-structure of $G$ arises from the unique outer automorphism of $A_{2n}$ of order two that sends $\alpha_i$ to $\beta_i$. If $a \in \Phi_K(G)$ is multipliable, then $a$ is the image of $\pm(\alpha_1 + \hdots + \alpha_s)$ for some $1 \leq s \leq n$. Suppose, without loss of generality, that $a$ is the image of $\alpha_1 + \hdots + \alpha_s$. 
 Consider the action of the image of $\Rx_b$ in $\RPF_{x}$ for $b$ the image of $-(\alpha_1+ \hdots + \alpha_{t})$ for some $1 \leq t \leq n$. Note that $\iota_{K,F}\left(\overline{x_b(H_0(E_{-(\alpha_1+ \hdots + \alpha_{t})},K))\cap G_{x,0}}\right)$ 
 is the image of $\xE_{-(\alpha_1 + \hdots +\alpha_t + \beta_1 + \hdots + \beta_t)}(E) \cap \GE_{x,0}$ in $\GE_{x,0}/\GE_{x,0+}$. Hence the orbit of  $\iota_{K,F}\left(\overline{x_b(H_0(E_{-(\alpha_1+ \hdots + \alpha_{t})},K))\cap G_{x,0}}\right)$ on $\iota_{K,F,r}\left(\overline{\fg_{x,r} \cap  \fg_a}\right)$ is contained in 
 $$\overline{(\fg \otimes_K \varpi_F^{d\va(\lambda_\alpha)/2}K) \cap \fgF_{x,r} \cap \left( \fgF_{\alpha_1+\hdots + \alpha_s} \oplus \fgF_{\beta_1+\hdots + \beta_s} \oplus \fgF_{-(\beta_1+ \hdots + \beta_t)} \oplus \fgF_{-(\alpha_1 + \hdots + \alpha_t)} \right)} \subset V'.$$
(Note that the last two summands can be deleted unless $s=t$.) Thus $V'$ is preserved under the action of the image of $\Rx_b$ in $\RPF_{x}$. Similarly (but more easily) one can check that the action of the image of $\Rx_b$ in $\RPF_{x}$ for all other $b \in \Phi(\RP_x)$ preserves $V'$, and the same is true for the image of $T \cap G_{x,0}$ in $\RPF_x$. Hence $\iota_{K,F}(\RP_x)$ stabilizes $V'$.

 The case of a general group $G$ follows using the observation that, if $a \in \Phi_K$ is multipliable, then each $\alpha \in \Phi_a$ is spanned by the roots of a connected component of the Dynkin diagram $\Dyn(G)$ of $\Phi(G)$ that is of type $A_{2n}$, together with the observation that the above explanation also works for $\Dyn(G)$ being a union of Dynkin diagrams of $A_{2n}$ that are permuted transitively by the action of the absolute Galois group of $K$. Thus $V'$ is preserved under the action of $\iota_{K,F}(\RP_x)$. 
 
 In order to show that $\iota_{K,F,r}$ is compatible with the action of $\RP_{x}$ as in Diagram \eqref{diagram-iota} for $p=2$, it remains to prove that $\RP_x$ preserves $\overline{\fg_{x,r} \cap (\ft \oplus \bigoplus\limits_{a \in \Phi_K^{\rm{nm}}} \fg_a)}$. We consider the action on $\overline{\fg_{x,r} \cap \ft }$ and $\overline{\fg_{x,r} \cap \bigoplus\limits_{a \in \Phi_K^{\rm{nm}}} \fg_a}$ separately. 
 
 We begin with the former, which is obviously preserved under the action of the image of $T \cap G_{x,0}$ in $\RP_x$. So consider the action of the image of $\Rx_b$ in $\RP_x$ for some $b \in \Phi(\bG_x) \subset \Phi_K$. If $b$ is non-multipliable in $\Phi_K$, then the image of the action lands in  $\overline{\fg_{x,r} \cap (\ft \oplus \fg_b)} \subset \overline{\fg_{x,r} \cap (\ft \oplus \bigoplus\limits_{a \in \Phi_K^{\rm{nm}}} \fg_a)}$. If $b \in \Phi_K^{\rm{mul}}$, then the image of the action is contained in  $\overline{\fg_{x,r} \cap (\ft \oplus \fg_b \oplus \fg_{2b})}$. However, by the assumption in our lemma, we have $b(x-x_0)-r \notin \Gamma'_b$ and hence $\overline{\fg_{x,r} \cap \fg_b} =\{0\}$. Therefore the image of the action of $\Rx_b$ on $\overline{\fg_{x,r} \cap \ft }$ is contained in $\overline{\fg_{x,r} \cap (\ft \oplus \fg_{2b})} \subset \overline{\fg_{x,r} \cap (\ft \oplus \bigoplus\limits_{a \in \Phi_K^{\rm{nm}}} \fg_a)}$.
 
 It remains to consider the action of $\RP_x$ on $\overline{\fg_{x,r} \cap \bigoplus\limits_{a \in \Phi_K^{\rm{nm}}} \fg_a}$. Note that the image of $T \cap G_{x,0}$ in $\RP_x$ preserves $\overline{\fg_{x,r} \cap \bigoplus\limits_{a \in \Phi_K^{\rm{nm}}} \fg_a}$. Thus it remains to consider the action of $\Rx_b(\bG_m)$ for some $b \in \Phi(\bG_x) \subset \Phi_K$, and we may restrict to the case that $\Dyn(G)$ is of type $A_{2n}$ with non-trivial Galois action as above. Let $b \in \Phi_K^{\rm{mul}}$, and assume without loss of generality that $b$ is the image of $\alpha_1 + \hdots + \alpha_s$ for some $1 \leq s \leq n$. Let $a \in \Phi_K^{\rm{nm}}$ with $\overline{\fg_{x,r} \cap \fg_a} \not\simeq \{0\}$, i.e. $a(x-x_0)-r \in \Gamma'_a$. 
 The assumption of the lemma implies that $\<\check b, a\> = 0$. Hence $a$ is the image of $\pm(\alpha_{s'} + \hdots + \alpha_{t'})$ for some $1 < s' < t' \leq n$ with $s' \neq s+1 \neq t'$, or of $\pm(\alpha_1 + \hdots + \alpha_{t'} + \beta_1 + \hdots + \beta_{s'})$ for some $1 \leq s',  t' \leq n$ with $s' \neq s \neq t'$ and $s' \neq t'$. In all cases, $\Rx_b(\bG_m)$ acts trivially on $\overline{\fg_{x,r} \cap \fg_a}$, and therefore $\Rx_b(\bG_m)$ preserves $\overline{\fg_{x,r} \cap \bigoplus\limits_{a \in \Phi_K^{\rm{nm}}} \fg_a}$. Similarly, i.e. using how non-multiple roots in the $A_{2n}$ case look like, we observe that if $b\in \Phi_K^{\rm{nm}}$, then $\Rx_b$ maps  $\overline{\fg_{x,r} \cap \bigoplus\limits_{a \in \Phi_K^{\rm{nm}}} \fg_a}$ to $\overline{\fg_{x,r} \cap (\ft \oplus \bigoplus\limits_{a \in \Phi_K^{\rm{nm}}} \fg_a)}$. 
 
Hence the Diagram \eqref{diagram-iota} commutes in the case $p=2$ if there does not exist $a \in \Phi_K^{\rm{mul}}$ with $a(x-x_0) \in \Gamma_a'$  such that $a(x-x_0)-r \in \Gamma_a'$ or such that there exists $b \in \Phi_K^{\rm{nm}}$ with $b(x-x_0)-r \in \Gamma_{b}'$ and $\<\check a, b\> \neq 0$. \qed

In the sequel we might abuse notation and identify $\PV_{x,r}$ with its image in $\PVF_{x,r}$ under $\iota_{K,F}$.

\section{Moy--Prasad filtration for different residual characteristics} \label{section-globalMP}
In this section we compare the Moy--Prasad filtration quotients for groups over nonarchimedean local fields of different residue field characteristics. In order to do so, we first introduce in Definition \ref{Def-good} the class of reductive groups that we are going to work with. 
We then show in Proposition \ref{Prop-examples} that this class contains reductive groups that split over a tamely ramified extension, i.e. those groups considered in \cite{ReederYu}, but also general simply connected and adjoint semisimple groups, among others. The restriction to this (large) class of reductive groups is necessary as the main result (Theorem \ref{Thm-general}) about the comparison of Moy--Prasad filtrations for different residue field characteristics does not hold true for some reductive groups that are not \good{} groups, see Remark \ref{Remark-failure}. 

\subsection{Definition and properties of good groups}

\begin{Def} \label{Def-good}
	We say that a reductive group $G$ over $K$, split over $E$, is \textit{\good{}} \index{definition}{good reductive group} if 
		there exist 
\begin{itemize}[label=$\cdot$]
			\item	an action of a finite cyclic group $\Gammapr=\<\gammapr\>$ \index{notation}{$\Gammapr$} \index{notation}{$\gammapr$}
			on the root datum $R(G)=(X,\Phi,\check X, \check \Phi)$ preserving the simple roots $\Delta$,
		 \item an element $u$ \index{notation}{$u$} generating the cyclic group $\Gal(\Et/K)$ and whose order $\abs{\Gal(\Et/K)}$ is divisible by $N$ where (throughout the remainder of the paper) we will write $\abs{\Gammapr}=p^s \cdot N$ \index{notation}{$N$} \index{notation}{$s$} for integers $s$ and $N$ with $(N,p)=1$
		 \end{itemize}
		 such that the following two conditions are satisfied. 
		\begin{enumerate}[(i)] 
			\item \label{item-root-action} The orbits of $\Gal(E/K)$ and $\Gammapr$ on $\Phi$ coincide, and, for every root $\alpha \in \Phi$, there exists $u_{1,\al} \in \Gal(E/K)$ \index{notation}{$u_{1,\al}$} such that 
			$$\gammapr(\alpha)=u_{1,\al}(\alpha) \quad \text{ and } \quad u_{1,\al} \equiv u \mod \Gal(E/\Et).$$
			\item	\label{item-Et}
			There exists a basis $\cB$ \index{notation}{$\cB$} of $X$ stabilized by $\Gal(E/\Et)$ and $\langle {\gammapr}^N\rangle$ on which the $\Gal(E/\Et)$-orbits and $\langle {\gammapr}^N\rangle$-orbits agree, and such that for any $B \in \cB$, there exists an element $v_{1,B} \in \Gal(E/K)$ \index{notation}{$v_{1,B}$}  satisfying
			$$ \gammapr(B)=v_{1,B}(B) \quad \text{ and } \quad v_{1,B} \equiv u \mod \Gal(E/\Et).$$				
		\end{enumerate}	
\end{Def}

\begin{Rem} \label{Rem-good}
	Note that condition \ref{item-root-action} of Definition \ref{Def-good} is equivalent to the condition 
			\begin{enumerate}[(i')] 
				\item \label{item-root-action-prime} The orbits of $\Gal(E/K)$ on $\Phi$ coincide with the orbits of $\Gammapr$ on $\Phi$, and there exist representatives $C_1, \hdots, C_n$\index{notation}{$C_i$} of the orbits of $\Gammapr$ on the connected components of the Dynkin diagram of $\Phi(G)$ satisfying the following. Denote by $\Phi_i$ \index{notation}{$\Phi_i$} the roots in $\Phi$ that are a linear combination of roots corresponding to $C_i$ ($1 \leq i \leq n$). Then for every root $\alpha \in \Phi_1 \cup \hdots \cup \Phi_n$ and $1 \leq t_1 \leq p^{s}N$, there exists $u_{t_1,\al} \in \Gal(E/K)$ such that
				$$ (\gammapr)^{t_1}(\alpha)=u_{t_1,\al}\alpha \quad \text{ and } \quad u_{t_1,\al} \equiv u^{t_1} \mod \Gal(E/\Et).$$
			\end{enumerate}			
		Condition \ref{item-Et} of Definition \ref{Def-good} is equivalent to the condition			
			\begin{enumerate}[(i')] 
				\setcounter{enumi}{1}
				\item	\label{item-Et-prime}
				There exists a basis $\cB$ of $X$ stabilized by $\Gal(E/\Et)$ and by $\langle {\gammapr}^N\rangle$ on which the $\Gal(E/\Et)$-orbits and $\langle {\gammapr}^N\rangle$-orbits agree, and such that there exist representatives $\{B_1, \hdots, B_{n'}\}$ for these orbits on $\cB$, and elements $v_{t_1,i} \in \Gal(E/K)$ for all  $1  \leq t_1 \leq p^sN$ and $1 \leq i \leq n'$   satisfying
				$$ (\gammapr)^{t_1}(B_i)=v_{t_1,i}(B_i) \quad \text{ and } \quad v_{t_1,i} \equiv u^{t_1} \mod \Gal(E/\Et).$$				
			\end{enumerate}	
\end{Rem}	

Before showing in Proposition \ref{Prop-examples} that a large class of reductive groups is \good{}, we prove a lemma that shows some more properties of \good{} groups. 

\begin{Lemma} \label{Lemma-good} 
	We assume that $G$ is a \good{} group, use the notation introduced in Definition \ref{Def-good} and Remark \ref{Rem-good}, and denote by $E_t$\index{notation}{$E_t$} the tamely ramified Galois extension of $K$ of degree $N$ contained in $E$. Then the following statements hold.
	\begin{enumerate}[(a)] 
		\item \label{item-old-ii} The basis $\cB$ of $X$ given in Property \ref{item-Et} is stabilized by $\Gal(E/E_t)$ and the $\Gal(E/E_t)$-orbits and $\langle{\gammapr}^N\rangle$-orbits on $\cB$ agree.
		\item \label{item-assumption} $G$ satisfies Assumption \ref{assumption}; more precisely, $T \times_K {E_t}$ is induced.
		\item \label{item-ii-2} We have $X^{{\gammapr}^{N}}=X^{\Gal(E/E_t)}$. Moreover, the action of $u$ on $X^{\Gal(E/E_t)}$ agrees with the action of $\gammapr$ on $X^{{\gammapr}^{N}}=X^{\Gal(E/E_t)}$, so $X^{\Gal(E/K)}=X^{\Gammapr}$.
	\end{enumerate}
\end{Lemma}
\textbf{Proof.} 
To show part \ref{item-old-ii}, consider a representative $B_i$ for a $\Gal(E/\Et)$-orbit on $\cB$ as in Remark \ref{Rem-good}. By Property \ref{item-Et-prime} there exists $v_{p^sN,i} \in \Gal(E/K)$ such that $v_{p^sN}(B_i)=({\gammapr})^{p^sN}(B_i)=B_i$ and $v_{p^sN,i} \equiv u^{p^sN} \mod \Gal(E/\Et)$. Choose $u_0 \in \Gal(E/K)$ such that $u_0 \equiv u \mod \Gal(E/\Et)$. Then we can write $v_{p^sN,i}=v\cdot u_0^{p^sN}$ for some $v \in \Gal(E/\Et)$ and $u_0^{p^sN}(B_i)=v^{-1}(B_i)$ is contained in the $\Gal(E/\Et)$-orbit of $B_i$. Note that the elements $u_0^{p^sNt_2}$ for $1 \leq t_2 \leq [(\Et):E_t]$ are in $\Gal(E/E_t)$ and form a set of representatives for $\Gal(E/E_t)/\Gal(E/\Et)$, and hence $\Gal(E/E_t)(B_i)=\Gal(E/\Et)(B_i)$. Thus $\cB$ is stabilized by $\Gal(E/E_t)$ and the $\Gal(E/E_t)$-orbits on $\cB$ coincide with the $\Gal(E/\Et)$-orbits, which coincide with the $\langle{\gammapr}^N\rangle$-orbits. This proves part \ref{item-old-ii}.

Part \ref{item-assumption} follows from part \ref{item-old-ii} by the definition of an induced torus. 

In order to show part \ref{item-ii-2}, note that $X^{\Gal(E/E_t)}$ is spanned (over $\bZ$) by
 $$\Bigl\{ \sum_{B \in \Gal(E/E_t)(B_i)} B\Bigr\}_{1 \leq i \leq n'}=\Bigl\{ \sum_{B \in \<{\gammapr}^N\>(B_i)} B\Bigr\}_{1 \leq i \leq n'}.$$ The $\bZ$-span of the latter equals $X^{{\gammapr}^N}$, which implies $X^{{\gammapr}^N}=X^{\Gal(E/E_t)}$. Using Definition \ref{Def-good}\ref{item-Et} and the observation that $u \mod \Gal(E\cap K^{\rm{tame}}/E_t)$ is a generator of $\Gal(E_t/K)$, we conclude that the action of $u$ on $X^{\Gal(E/E_t)}$ agrees with the action of $\gammapr$ on $X^{{\gammapr}^{N}}=X^{\Gal(E/E_t)}$ and that 
 \begin{align*}  & X^{\Gal(E/K)}=\left(X^{\Gal(E/E_t)}\right)^{\Gal(E_t/K)}=\left(X^{{\gammapr}^N}\right)^{\gammapr}=X^{\Gammapr} .  & & \qed
 \end{align*}

\begin{Prop} \label{Prop-examples}
Examples of \good{} groups include 
	\begin{enumerate}[(a)]
	\item \label{item-tame} reductive groups that split over a tamely ramified field extension of $K$,
	\item \label{item-adjoint} simply connected or adjoint (semisimple) groups,
	\item \label{item-products} products of \good{} groups,
	\item \label{item-restriction-of-scalars} groups that  
	 are the restriction of scalars of good groups along finite separable field extensions.
	\end{enumerate}
\end{Prop}
\textbf{Proof.}\\
\textbf{\ref{item-tame}} Part \ref{item-tame} follows by taking $\Gammapr=\Gal(E/K)$ and $u=\gammapr$.

\textbf{\ref{item-adjoint}} Part \ref{item-adjoint} can be deduced from \ref{item-products} and \ref{item-restriction-of-scalars} (whose proofs do not depend on \ref{item-adjoint}) as follows. If $G$ is a simply connected or adjoint group then $G$ is the direct product of restrictions of scalars of simply connected or adjoint absolutely simple groups. Hence by \ref{item-products} and \ref{item-restriction-of-scalars} it suffices to show that, if $G$ is a simply connected or adjoint absolutely simple group, then $G$ is \good{}. Recall that these groups are classified by choosing the attribute ``simply connected'' or ``adjoint'' and giving a connected finite Dynkin diagram together with an action of the absolute Galois group $\Gal(\ov \bQ_p/K)$  on it. We distinguish the two possible cases.\\
Case 1: $G$ splits over a cyclic field extension $E$ of $K$. Then take $\Gammapr=\Gal(E/K)$ and $u=\gammapr$ or $u=1$ according as the field extension is tamely ramified or wildly ramified, and choose $\cB$ to be the set of simple roots of $G$, if $G$ is adjoint, and the set of fundamental weights dual to the simple co-roots of $G$ (i.e. those weights pairing with one simple co-root to 1, and with all others to 0), if $G$ is simply connected. \\
Case 2: $G$ does not split over a cyclic field extension. Then $G$ has to be of type $D_4$ and split over a field extension $E$ of $K$ of degree six with $\Gal(E/K) \simeq S_3$, where $S_3$ is the symmetric group on three letters. In this case we observe (using that $G$ is simply connected or adjoint) that the orbits of the action of $\Gal(E/K)$ on $X$ are the same as the orbits of a subgroup $\bZ/3\bZ \subset \Gal(E/K) \simeq S_3$. Moreover, as $S_3$ does not contain a normal subgroup of order two, i.e. there does not exist a tamely ramified Galois extension of $K$ of degree three, this case can only occur if $p=3$, and we can choose $\Gammapr=\bZ/3\bZ$, $u$ the nontrivial element in $\Gal(\Et/K)\simeq \bZ/2\bZ$, and $\cB$ as in Case 1 to see that $G$ is \good{}.

\textbf{\ref{item-products}} In order to show part \ref{item-products}, suppose that $G_1, \hdots, G_k$ are \good{} groups with splitting fields $E_1, \hdots, E_k$ and corresponding cyclic groups $\Gammapr_1=\<\gammapr_1\>, \hdots, \Gammapr_k=\<\gammapr_k\>$
and generators $u_i \in \Gal(E_i \cap K^{\rm{tame}}/K), 1 \leq i \leq k$. Let $G=G_1 \times \hdots \times G_k$. Then $G$ splits over the composition field $E$ of $E_1, \hdots, E_k$, and $\abs{\Gal(E\cap K^{\rm{tame}}/K)}$ is the smallest common multiple of $\abs{\Gal(E_i\cap K^{\rm{tame}}/K)}, 1 \leq i \leq k$. Choose a generator $u$ of $\Gal(E\cap K^{\rm{tame}}/K)$. For $i \in [1,k]$, the image of $u$  in $\Gal(E_i \cap K^{\rm{tame}}/K)$ equals $u_i^{r_i}$ for some integer $r_i$ coprime to $\abs{\Gal(E_i \cap K^{\rm{tame}}/K)}$, which we assume to be coprime to $p$ by adding $\abs{\Gal(E_i \cap K^{\rm{tame}}/K)}$ if necessary. Hence $({\gammapr}_i)^{r_i}$ is a generator of $\Gammapr_i$, and we define $\gammapr=({\gammapr}_1)^{r_1} \times \hdots \times ({\gammapr}_k)^{r_k}$ and $\Gammapr=\<\gammapr\>$. Note that the order $\abs{\Gammapr}=p^sN$ of $\Gammapr$ is the smallest common multiple of $\abs{\Gammapr_i}, 1 \leq i \leq k$, and hence $N$ divides $\abs{\Gal(\Et/K)}$. 
By \ref{Def-good}\ref{item-root-action}  if $\alpha \in \Phi(G_i)$, then there exists $\ov u_{1,\alpha} \in \Gal(E_i/K)$ such that 
$$ \gammapr(\alpha)=(\gammapr_i)^{r_i}(\alpha)=\ov u_{1,\al} \alpha \, \text{ with } \ov u_{1,\al} \equiv u_i^{r_i}  \equiv u \text{ in } \Gal(E_i \cap K^{\rm{tame}}/K).$$
Let $u_{1,\al}$ be a preimage of $\ov u_{1,\al}$ in $\Gal(E/K)$. Using that 
\begin{eqnarray*}
	& & \abs{\Gal(E/E \cap E_i^{\rm{tame}})}\abs{\Gal(E\cap E_i^{\rm{tame}}/E_i)}\abs{\Gal(E_i/K)} \\
	&=& \abs{\Gal(E/K)}
	=\abs{\Gal(E/E \cap K^{\rm{tame}})}\abs{\Gal(E\cap K^{\rm{tame}}/E_i\cap K^{\rm{tame}})}\abs{\Gal(E_i\cap K^{\rm{tame}}/K)},
\end{eqnarray*}
we obtain by considering the factors prime to $p$ that $\abs{\Gal(E\cap E_i^{\rm{tame}}/E_i)}=\abs{\Gal(E \cap K^{\rm{tame}}/E_i \cap K^{\rm{tame}})}.$
 Moreover, the kernel of 
 $\Gal(E\cap E_i^{\rm{tame}}/E_i) \ra \Gal(E \cap K^{\rm{tame}}/E_i \cap K^{\rm{tame}})$, where the map arises from reduction mod $\Gal(E \cap E_i^{\rm{tame}}/E \cap K^{\rm{tame}})$,
  has order a power of $p$, hence is trivial; so we deduce that the map is an isomorphism. Thus we can choose an element $u_0 \in \Gal(E/E_i)\subset \Gal(E/K)$ such that $u_0 \equiv u^{\abs{\Gal(E_i \cap K^{\rm{tame}}/K)}} \mod \Gal(E/\Et)$, because $ u^{\abs{\Gal(E_i \cap K^{\rm{tame}}/K)}} \in \Gal(\Et/E_i \cap K^{\rm{tame}})$. Since $u_{1,\al} \equiv u \mod \Gal(E/E_i \cap K^{\rm{tame}})$ and $ u^{\abs{\Gal(E_i \cap K^{\rm{tame}}/K)}}$ is a generator of $\Gal(\Et/E_i \cap K^{\rm{tame}})$, by multiplying $u_{1,\al}$ with powers of $u_0 \in \Gal(E/E_i)$
 if necessary we can ensure that $u_{1,\al} \equiv u \mod \Gal(E/\Et)$. As $\Gal(E/E_i)$ fixes $\alpha$, we also have $\gammapr(\alpha)=u_{1,\al}(\al)$, and we conclude that $G$ satisfies Property \ref{item-root-action} of Definition \ref{Def-good} for all $\alpha \in \Phi(G)=\coprod_{i=1}^k\Phi(G_i)$.

  Choosing $\cB$ to be the union of the bases $\cB_i$ corresponding to the \good{} groups $G_i$ (by viewing $X_i$ embedded into $X:=X_1 \times \hdots \times X_k$), we conclude similarly that $G$ satisfies Property \ref{item-Et}. This proves that $G$ is a \good{} group and finishes part \ref{item-products}.

\textbf{\ref{item-restriction-of-scalars}} Let $G=\Res_{F/K} \wt G$ for $\wt G$ a \good{} group over $F$, $K \subset F \subset E$. Then there exists a corresponding $\Gal(E/K)$-stable decomposition $X=\bigoplus_{i=1}^{d} X_i$, where $d=[F:K]$, together with a decomposition of $\Phi$ as a disjoint union $\coprod\limits_{1 \leq i \leq f} \wt \Phi_i$ such that $\Gal(E/K)$ acts transitively on the set of subspaces $X_i$ with $\Stab_{\Gal(E/K)}(X_i) \simeq \Gal(E/F)$, and $(X_i, \wt \Phi_i, \check X_i, \check{\wt{\Phi_i}})$ is isomorphic to the root datum $R(\wt G)$ of $\wt G$ for $1 \leq i \leq f$. 
 We suppose without loss of generality that the fixed field of $\Stab_{\Gal(E/K)}(X_1)$ is $F$, i.e. $\Stab_{\Gal(E/K)}(X_1)=\Gal(E/F)$, and we write $d=d_p \cdot d_{p'}$, where $d_p$ is a power of $p$ and $d_{p'}$ is coprime to $p$.
 As $\wt G$ is \good{}, there exist a cyclic group $\wt \Gamma = \< \wt \gamma \>$ acting on $(X_1,\wt \Phi_1, \Delta_1)$ and a generator $\wt u$ of $\Gal(E \cap F^{\rm{tame}}/F)$ satisfying the conditions in Definition \ref{Def-good}. Fix a splitting $\Gal(E\cap F^{\rm{tame}}/F) \hookrightarrow \Gal(E/F)$, and let $\wt u_0$ be the image of $\wt u$ under the composition $\Gal(E\cap F^{\rm{tame}}/F)  \hookrightarrow \Gal(E/F) \hookrightarrow \Gal(E/K)$. Note that we have a commutative diagram (where $N'=\abs{\Gal(E\cap F^{\rm{tame}}/F)}$)
\begin{eqnarray*} 
	\xymatrix{
		\Gal(E\cap F^{\rm{tame}}/F)   \ar@{^{(}->}[r] \ar[d]_{\simeq} &  \Gal(E/F) \ar@{^{(}->}[r] \ar[d]_{\simeq} & \Gal(E/K)  \ar[d]_{\simeq}  \\
		\bZ/ N' \bZ  \ar@{^{(}->}[r]  & \bZ/ N' \bZ \ltimes \Gal(E/E\cap F^{\rm{tame}}) \ar@{^{(}->}[r] & \bZ/ (N' d_{p'}) \bZ \ltimes \Gal(E/E\cap K^{\rm{tame}})   \quad , 
	} 
\end{eqnarray*}
Hence  we can choose $u_0 \in \Gal(E/K)$ such that 
 $$u_0^{d} \equiv \wt u_0 \mod \Gal(E/E \cap K^{\rm{tame}}),$$ and $u := u_0 \mod \Gal(E/E \cap K^{\rm{tame}})$ is a generator of  $\Gal(E \cap K^{\rm{tame}}/K)$ (because $d=d_pd_{p'}$ with $d_p$ invertible in $\bZ/ (N' d_{p'}) \bZ $).
After renumbering the subspaces $X_i$ for $i>1$, if necessary, we can choose elements $\gamma_{t_2d_{p'}} \in 
 \Gal(E/K)$
   with $$\gamma_{t_2d_{p'}} \equiv u_0 = u \mod \Gal(E \cap K^{\rm{tame}}/K)$$
  for $1 \leq t_2 \leq d_p$ 
   such that if we set $\gamma_{t_1+t_2d_{p'}}= u_0$  for $1 \leq t_1 < d_{p'}, 0 \leq t_2 < d_p$ then $\gamma_i(X_i)=X_{i+1}$, $1 \leq i < d$ and $\gamma_d(X_d)=X_1$.  By multiplying $\gamma_d$ by an element in $\Gal(E/E\cap K^{\rm{tame}})$ if necessary, we can assume that $\gamma_d \circ \gamma_{d-1} \circ \hdots \circ \gamma_1 =\wt u_0$. Define $\gammapr \in \Aut(R(G), \Delta)$ by 
$$ X=\bigoplus_{i=1}^{d} X_i \ni (x_1, \hdots, x_d) \mapsto (\wt \gamma \circ \wt u_0^{-1} \circ \gamma_d x_d, \gamma_1x_1,\gamma_2x_2, \hdots, \gamma_{d-1}x_{d-1}).$$
Then the cyclic group $\Gammapr=\<\gammapr\>$ preserves $\Delta$, and we claim that $\Gammapr$ and $u$ satisfy the conditions for $G$ in Definition \ref{Def-good}. 

Property \ref{item-root-action} of Definition \ref{Def-good} is satisfied by the construction of $\gammapr$. 

In order to check Property \ref{item-Et}, 
  let $\wt \cB$ be a basis of $X_1 \subset X$ stabilized by $\Gal(E/E\cap F^{\rm{tame}})$ with a set of representatives $\{\wt B_1, \hdots, \wt B_{\wt n'}\}$ and $\wt v_{t_1,i} \in \Gal(E/F)$ with $(\wt \gamma)^{t_1}(B_i)=\wt v_{t_1,i}(B_i)$ ($1 \leq t_1 \leq p^sN/d$) satisfying all conditions of Property \ref{item-Et-prime} of Remark \ref{Rem-good} for $\wt G$. For $1 \leq i  \leq \wt n'$ and $1 \leq j \leq d_{p'}$, define
$$ B_{(i-1)d_{p'}+j}=u_0^{j-1}(\wt B_i)=\gamma_{j-1} \circ \cdots \circ \gamma_{1} (\wt B_i).$$
Note that $\langle{\gammapr}^{N}\rangle(X_1) = \coprod\limits_{0 \leq i < d_p} X_{1+id_{p'}}$, and hence, setting $n'=\wt n' \cdot d_{p'}$, the set
$$\cB= \bigcup_{1 \leq i \leq n'} \<{\gammapr}^{N}\>(\{B_i\})$$ 
forms a basis of $X$ (because $\gammapr^N$ has order $d_p$).
We will show that $\cB$ satisfies Property \ref{item-Et-prime} of Remark \ref{Rem-good} with set of orbit representatives $\{B_i\}_{1 \leq i \leq n'}$ (and hence satisfies Property \ref{item-Et} of Definition \ref{Def-good}).

  For $1 \leq t \leq p^sN, 1 \leq i  \leq \wt n', 1 \leq j \leq d_{p'}$, we define $v_{t,(i-1)d_{p'}+j} \in \Gal(E/K)$ by
\begin{equation*}
	v_{t,(i-1)d_{p'}+j}
	     =\left\{ \begin{array}{ll}
	     	\gamma_{j-1+t} \circ \cdots \circ \gamma_{j} & \text{ if } j+t \leq d\\
	     	\gamma_{t_2} \circ \cdots \gamma_1 \circ \wt v_{t_1,i} \circ \gamma_1^{-1} \circ \cdots \gamma_{j-1}^{-1} & \text{ if } j+t > d, t=dt_1+t_2-j+1
		\end{array}\right. .
\end{equation*}
Then using $(\gammapr)^d|_{X_1}=\wt \gamma$ and $\wt \gamma ^{t_1}(\wt B_i)=\wt v_{t_1,i}(\wt B_i) \in X_1$, 
we obtain 
$$(\gammapr)^t(B_i)=v_{t,i}(B_i) \, \text{ for all  }1 \leq t \leq p^sN \text{ and } 1 \leq i \leq n'.$$
Moreover, since
\begin{eqnarray*}
	 \wt v_{t_1,i} \equiv \wt u^{t_1} \mod \Gal(E/E\cap F^{\rm{tame}}) %\\
	&\Ra& \wt v_{t_1,i} \equiv \wt u_0^{t_1} \equiv  u_0^{dt_1} \equiv u^{dt_1} \mod \Gal(E/ E\cap K^{\rm{tame}})
\end{eqnarray*}
and $\gamma_k \equiv u \mod  \Gal(E\cap K^{\rm{tame}}/K)$ for all $1 \leq k < d$ by definition, we obtain
\begin{equation} \label{eqn-equivu} v_{t,i} \equiv u^t \mod  \Gal(E/E\cap K^{\rm{tame}}) \quad \text{ for all } 1 \leq t \leq p^sN \text{ and } 1 \leq i \leq n'.
\end{equation}

This shows that the action of ${(\gammapr)}^{t_1}$ on $B_i$ for $1 \leq t_1 \leq p^sN$ and $1 \leq i \leq n'$ is as required by Condition \ref{item-Et-prime} of Remark \ref{Rem-good}. It remains to show that $\cB$ is $\Gal(E/\Et)$-stable and that the $\Gal(E/\Et)$-orbits coincide with the $\langle{\gammapr}^N\rangle$-orbits.

In order to do so, note that Equation \eqref{eqn-equivu} implies in particular that for $1 \leq t_2 \leq d_p$, we have $v_{Nt_2,i} \equiv u^{N t_2} \mod \Gal(E/\Et)$, and hence $v_{Nt_2,i} \in \Gal(E/E_t)$ and 
\begin{equation}\label{eqn-inclusion}
	\<{\gammapr}^N\>(B_i) \subset \Gal(E/E_t)(B_i),
\end{equation}
	 where $E_t$ is the tamely ramified degree $N$ field extension of $K$ inside $E$. Let us denote by $\wt E_t$ the tamely ramified Galois extension of $F$ of degree $N/d_{p'}$ contained in $E$.  Note that $E_t$ is the maximal tamely ramified subextension of $\wt E_t$ over $K$, and $[\wt E_t:E_t]=d_p$. As $\wt G$ is \good{}, we obtain from Property \ref{item-Et} of Definition \ref{Def-good} and Lemma \ref{Lemma-good} \ref{item-old-ii} that $$\<{\gammapr}^{N{d_p}}\>(B_i)  = \<{\wt \gamma}^{N/{d_p'}}\>(B_i) = \Gal(E/E \cap F^{\rm{tame}})(B_i)= \Gal(E/\wt E_t)(B_i) . $$
Using $\<{\gammapr}^{N}\>(X_1) = \coprod\limits_{0 \leq i < d_p} X_{1+id_{p'}}$ and the inclusion \eqref{eqn-inclusion}, we deduce that $$\abs{\Gal(E/E_t)(B_i)} \geq \abs{\langle{\gammapr}^N\rangle(B_i)}=d_p \cdot \abs{\<{\gammapr}^{N{d_p}}\>(B_i)} =d_p \cdot \abs{\Gal(E/\wt E_t)(B_i)} \leq  \abs{\Gal(E/E_t)(B_i)},$$
which implies that $\langle{\gammapr}^N\rangle(B_i) = \Gal(E/E_t)(B_i) \supset \Gal(E/\Et)(B_i)$. In order to show that $\langle{\gammapr}^N\rangle(B_i) =  \Gal(E/\Et)(B_i)$, we observe that $\Gal(E/E\cap F^{\rm{tame}})$ is a subgroup of $\Gal(E/\Et)$ of index $d_p$ coprime to the index $N/d_{p'}$ of $\Gal(E/E\cap F^{\rm{tame}})$ inside $\Gal(E/F)$. Therefore $\Gal(E/\Et) \cap \Gal(E/F) = \Gal(E/E\cap F^{\rm{tame}})$ inside $\Gal(E/K)$. As $\Gal(E/F)$ is the stabilizer of $X_1$ in $\Gal(E/K)$, we deduce that there exist $d_p$ representatives in $\Gal(E/\Et)$ of the $d_p$ classes in $\Gal(E/\Et)/\Gal(E/E\cap F^{\rm{tame}})$ mapping $X_1$ to $d_p$ distinct components $X_i$ of $X$. In particular, we obtain that $$\abs{\Gal(E/\Et)(B_i)} \geq d_p \abs{\Gal(E/E \cap F^{\rm{tame}})(B_i)}=d_p \abs{\<{\gammapr}^{N{d_p}}\>(B_i)} =\abs{\<{\gammapr}^{N}\>(B_i) },$$
 and hence the $\Gal(E/\Et)$-orbits on $\cB$ agree with the $\langle{\gammapr}^{N}\rangle$-orbits on $\cB$. This finishes the proof that Property \ref{item-Et-prime} of Remark \ref{Rem-good} and hence Property \ref{item-Et} of Definition \ref{Def-good} is satisfied for our choice of $\Gammapr$ and $u$, and hence $G$ is \good{}.
\qed

From now on we assume that our group $G$ is \good{}.

\subsection{Construction of $G_q$} \label{section-Gq} 
In this section we define reductive groups $G_q$ over nonarchimedean local fields with arbitrary positive residue field characteristic $q$ whose Moy--Prasad filtration quotients are in a certain way (made precise in Theorem \ref{Thm-general}) the ``same'' as those of the given \good{} group $G$ over $K$.

For the rest of the paper, assume $x \in \sB(G,K)$ \index{notation}{$x$} is a rational point of order $m$. Here \index{definition}{rational point} \textit{rational} \label{page-rational} means that $\psi(x)$ is in $\bQ$ for all affine roots $\psi \in \Root$, and the \index{definition}{order} \textit{order of a point in the Bruhat--Tits building}  \label{page-order} $m$ \index{notation}{$m$} of $x$ is defined to be the smallest positive integer such that $\psi(x)\in \frac{1}{m} \bZ$ for all affine roots $\psi \in \Root$. 

Fix a prime number $q$\index{notation}{$q$}, and let $\Gammapr$ be the finite cyclic group acting on $R(G)$ as in Definition \ref{Def-good}. Let $\Fp$\index{notation}{$\Fp$} be a Galois extension of $K$ containing the splitting field of $(x^2-2)$ over $E$, such that 
\begin{itemize}%[label=$\cdot$]
	\item $M:=[F:K]$\index{notation}{$M$}  is divisible by the order $p^sN$ of the group $\Gammapr$, 
	\item $M$ is divisible by the order $m$ of the point $x \in \sB(G,K)$.
\end{itemize} 
This implies that the image of $x$ in $\sB(G_\Fp,\Fp)$ is hyperspecial, and by the last condition the set of valuations $\Gamma'_a$ (defined in Section \ref{section-affine-roots}) is contained in $\va(\Fp)$ for all $a \in \Phi_K$. In particular, $F$ satisfies all assumptions made in Section \ref{section-iota} in order to define $\iota_{K,F}$ and $\iota_{K,F,r}$. For later use, denote by $\varpi_F$\index{notation}{$\varpi_F$} a uniformizer of $F$ such that
 $\varpi_F^{[F:E]}\equiv \varpi_E \mod \varpi_F^{[F:E]+1}$, and let $\cO_F$ be the ring of integers of $F$.

 Let $K_q$\index{notation}{$K_q$} be the splitting field of $x^M-1$ over $\bQ_q^{ur}$, with ring of integers $\cO_q$\index{notation}{$\cO_q$} and uniformizer $\varpi_q$\index{notation}{$\varpi_q$}. Let $\Fq=K_q[x]/(x^M-\varpi_q)$\index{notation}{$\Fq$} with uniformizer $\varpi_{\Fq}$\index{notation}{$\varpi_{\Fq}$} satisfying $\varpi_{\Fq}^M=\varpi_{q}$ and ring of integers $\cO_{F_q}$\index{notation}{$\cO_{\Fq}$}. 
Recall that every reductive group over $K_q$ is quasi-split, and that there is a one to one correspondence between (quasi-split) reductive groups over $K_q$ with root datum $R(G)$ and elements of $\Hom(\Gal(\ov \bQ_q / K_q), \Aut(R(G),\Delta))/\mbox{Conjugation by} \Aut(R(G),\Delta)$, where $\Aut(R(G),\Delta)$ denotes the group of automorphisms of the root datum $R(G)$ that stabilize $\Delta$.
Thus we can define a reductive group $G_q$ \index{notation}{$G_q$} over $K_q$ by requiring that $G_q$ has root datum $R(G)$ and that the action of $\Gal(\ov \bQ_q/K_q)$ on $R(G)$ defining the $K_q$-structure factors through $\Gal(\Fq/K_q)$ and is given by
$$\Gal(\Fq/K_q)\simeq \bZ/M\bZ \xrightarrow{1 \mapsto \gammapr} \Gammapr \ra \Aut(R(G), \Delta),$$
where the last map is the action of $\Gammapr$ on $R(G)$ as in Definition \ref{Def-good}.
This means that $G_q$ is already split over $\Eq:=K_q[x]/(x^{p^sN}-\varpi_q)$\index{notation}{$\Eq$}. 
Note that by construction, Definition \ref{Def-good} and Lemma \ref{Lemma-good}, the restricted root data of $G_q$ and $G$ agree:
\begin{equation*}
	R_{K_q}(G_q)=R_K(G),
\end{equation*}
and we have for all $\alpha \in \Phi=\Phi(G)=\Phi(G_q)$
\begin{equation} \label{equation-Gal-orbits}
\abs{\Gal(E/K)\cdot \alpha}=\abs{\Gal(F_q/K_q)\cdot \alpha}.
\end{equation}
All objects introduced in Section \ref{section-basics} can also be constructed for $G_q$, and we will denote them by the same letter(s), but with a $G_q$ in parentheses to specify the group; e.g., we write $\Gamma_a'(G_q)$\index{notation}{$\Gamma_a'(G_q)$}. 

\subsection{Construction of $x_q$} \label{section-xq}
In order to compare the Moy--Prasad filtration quotients of $G_q$ with those of $G$ at $x$, we need to specify a point $x_q$ in the Bruhat--Tits building $\sB(G_q,K_q)$ of $G_q$. To do so, choose a maximal split torus $S_q$\index{notation}{$S_q$} in $G_q$ with centralizer denoted by $T_q$\index{notation}{$T_q$}, and fix a Chevalley--Steinberg system $\{\xFq_\alpha\}_{\alpha \in \Phi}$\index{notation}{$\xFq_\alpha$} for $G_q$ with respect to $T_q$. For later use, we choose the Chevalley--Steinberg system to have signs $\eps_{\alpha,\beta}$ as in Definition \ref{Def-signs}, i.e. \index{notation}{$m_\al^{\Fq}$}
$$m_\al^{\Fq}:= \xFq_\alpha(1)\xFq_{-\alpha}(\eps_{\al,\al})\xFq_{\alpha}(1) \in N_{G_q}(T_q)(\Fq),$$
where $N_{G_q}(T_q)$ denotes the normalizer of $T_q$ in $G_q$, and
$$ \Ad(m_\al^{\Fq})(\Lie(\xFq_\beta)(1))=\eps_{\alpha,\beta}\Lie(\xFq_{s_\alpha(\beta)})(1).$$ 
 Using the valuation constructed in Section \ref{Sec-rootgroups} attached to this Chevalley--Steinberg system, we obtain a point $x_{0,q}$\index{notation}{$x_{0,q}$} in the apartment $\sA_q$ of $\sB(G_q, K_q)$ corresponding to $S_q$.
Fixing an isomorphism $f_{S,q}: X_*(S)\ra X_*(S_q)$ that identifies $R_K(G)$ with $R_{K_q}(G_q)$, we define an isomorphism of affine spaces $f_{\sA,q}: \sA \ra \sA_q$ by 
\begin{equation} \label{equation-f-def}
	f_{\sA,q}(y)=x_{0,q}+f_{S,q}(y-x_0)- \frac{1}{4} \sum_{a \in \Phi_K^{+,\text{mul}}} \va(\lambda_a) \cdot \check a,
\end{equation}
where $\Phi_K^{+,\text{mul}}$\index{notatoin}{$\Phi_K^{+,\text{mul}}$} are the positive multipliable roots in $\Phi_K$, $\lambda_a \in (E_{\alpha})^1_{\text{max}}(G)$ for some $\alpha \in \Phi_a$, and $\check a$ is the coroot of $a$, so we have $\check a (a)=2$. 
We define $x_q:=f_{\sA,q}(x)$. \index{notation}{$x_q$} 

\begin{Lemma} \label{lemma-affine-roots}
	The isomorphism $f_{\sA,q}: \sA \ra \sA_q$ induces a bijection of affine roots $\Rootq(\sA_q) \ra \Root(\sA), \psi \mapsto \psi \circ f_{\sA,q}$.

	 Moreover, we have for all $a \in \Phi_K$ and $r\in \bR$ that $r-a(x-x_{0}) \in \Gamma_a'(G)$ if and only if $r-a(x_q-x_{0,q}) \in \Gamma_a'(G_q)$.
\end{Lemma}
\textbf{Proof.}
As the set of affine roots for $G$ on $\sA$ (and analogously for $G_q$ on $\sA_q$) is  
$$\Root = \Root(\sA) = \left\{y \mapsto a(y-x_{0})+\gamma' \, | \, a \in \Phi_K, \gamma' \in \Gamma_a' \right\},$$
we need to show that, for every $a\in \Phi_K=\Phi_K(G)=\Phi_{K_q}(G_q)$, we have
\begin{equation} \label{equation-gammas-equality} \Gamma'_a(G) = \Gamma'_{a}(G_q) -\frac{1}{4} \sum_{b \in \Phi_K^{+,\text{mul}}} \va(\lambda_b) \cdot \check b(a).
\end{equation}
Let us fix  $a \in \Phi_K$, and $\alpha \in \Phi_a \subset \Phi=\Phi(G)=\Phi(G_q)$. Recall that $E_\alpha(G)$ is the fixed subfield of $E$ under the action of $\Stab_{\Gal(E/K)}(\alpha)$. Using Equation (\ref{equation-Gal-orbits}) on page \pageref{equation-Gal-orbits}, we obtain
\begin{eqnarray*}
	[E_\alpha(G):K]&=& \frac{\abs{\Gal(E/K)}}{\abs{\Stab_{\Gal(E/K)}(\alpha)}}=\abs{\Gal(E/K)\cdot \alpha} = \abs{\Gal(F_q/K_q)\cdot \alpha} \\
	&=&\frac{\abs{\Gal(\Fq/K_q)}}{\abs{\Stab_{\Gal(\Fq/K_q)}(\alpha)}} = [E_\alpha(G_q):K_q],
\end{eqnarray*}
and hence
\begin{equation} \label{equation-valuation-Ealphas}
	\va(E_\alpha(G)-\{0\})= [E_\alpha(G)/K]^{-1} \cdot \bZ= [E_\alpha(G_q)/K_q]^{-1} \cdot \bZ=\va(E_\alpha(G_q)-\{0\}).
\end{equation}

 Note that the Dynkin diagram $\Dyn(G)$ of $\Phi(G)$ is a disjoint union of irreducible Dynkin diagrams, and  if $a$ is a multipliable root, then $\alpha$ is contained in the span of the simple roots of a Dynkin diagram of type $A_{2n}$. Thus by Equation (\ref{equation-valuation-Ealphas}) and the description of $\Gamma'_a$ as in Equation (\ref{equation-Ealpha-nonmult}) on page \pageref{equation-Ealpha-nonmult}, the Equality (\ref{equation-gammas-equality}) holds for $\alpha$ in the span of simple roots of an irreducible Dynkin diagram of any type other than $A_{2n}$, $n\in \bZ_{>0}$, or in the span of an irreducible Dynkin diagram of type $A_{2n}$ whose $2n$ simple roots lie in $2n$ distinct Galois orbits.  We are therefore left to prove the lemma in the case of $\Dyn(G)$ being a disjoint union of finitely many $A_{2n}$ whose simple roots form $n$ orbits under the action of $\Gal(E/K)$. An easy calculation (see the proof of Lemma \ref{Lemma-iota} for details) shows that, in this case, the positive multipliable roots of $\Phi_K$ form an orthogonal basis for the subspace of $X^*(S) \otimes \bR$ generated by $\Phi_K$, where by ``orthogonal'' we mean that $\check b(a)=0$ if $a$ and $b$ are distinct positive multipliable roots, and that, if $b \in \Phi_K$ and $b=\sum\limits_{a \in \Phi_K^{+,\text{mul}}} \kappa_a a$ is not multipliable, then $\sum\limits_{a \in \Phi_K^{+,\text{mul}}} \kappa_a  \in 2 \cdot \bZ$. Moreover, by the definition of $K_q$ and $\Fq$, it is easy to check that for $\lambda_q \in (E_{\alpha})^1_{\text{max}}(G_q)$, we have $\va(\lambda_q) \in 2 \cdot \va(E_\alpha-\{0\})$. Thus using the description of $\Gamma'_a$ as in Equation (\ref{equation-Ealpha-mult}) on page \pageref{equation-Ealpha-mult} and Equation (\ref{equation-Ealpha-div}) on page \pageref{equation-Ealpha-div}, we see that the desired Equation (\ref{equation-gammas-equality}) holds.
 
 The second claim of the lemma follows from combining Equation \eqref{equation-gammas-equality} and the definition of $x_q$ using the map in Equation \eqref{equation-f-def} on page \pageref{equation-f-def}.
 \qed
 
Note that Lemma \ref{lemma-affine-roots} implies in particular that $x_q$ is also a rational point of order $m$.
Let us denote the reductive quotient of $G_q$ at $x_q$ by $\RP_{x_q}$\index{notation}{$\RP_{x_q}$}; the corresponding Moy--Prasad filtration groups by $G_{x_q,r}$\index{notation}{$G_{x_q,r}$}, $r\geq 0$; the Lie algebra filtration by $\fg_{x_q,r}, r \in \bR$\index{notation}{$\fg_{x_q,r}$}; and the filtration quotients of the Lie algebra by $\PV_{x_q,r}, r \in \bR$\index{notation}{$\PV_{x_q,r}$}.
Then using Lemma \ref{lemma-reductive-quotient}, we obtain the following corollary to Lemma \ref{lemma-affine-roots}.
\begin{Cor} \label{Cor-root-data}
	The root data $R(\RP_x)$ and $R(\RP_{x_q})$ are isomorphic.
\end{Cor}

\subsection{Global Moy--Prasad filtration representation} \label{section-global-MP-filtration}
Since $R(\RP_x)=R(\RP_{x_q})$ (Corollary \ref{Cor-root-data}), we can define a split reductive group scheme $\ZG$\index{notation}{$\ZG$} over $\bZ$ by requiring that $R(\ZG)=R(\RP_x)$, and then $\ZG_{\fpb}\simeq \RP_x$ and $\ZG_{\fqb} \simeq \RP_{x_q}$; i.e., we can define the reductive quotient ``globally''. In this section we show that we can  define not only the reductive quotient globally, but also the action of the reductive quotient on the Moy--Prasad filtration quotients. More precisely, we will prove the following theorem, where $N$ is as in Definition \ref{Def-good}, i.e., in particular, $N$ is coprime to $p$.

\begin{Thm} \label{Thm-general} 
	Let $r$ be a real number, and keep the notation from Section \ref{section-Gq} and \ref{section-xq}, so
	 $G$ is a \good{} reductive group over $K$ and $x$ a rational point of $\sB(G,K)$. 
	Then there exists a split reductive group scheme $\ZG$ over $\ov \bZ[1/N]$ acting on a free $\ov \bZ[1/N]$-module $\ZV$
	satisfying the following. For every prime $q$ coprime to $N$, there exist isomorphisms $\ZG_{\ov \bF_q} \simeq \RP_{x_q}$ and $\ZV_{\ov \bF_q} \simeq \PV_{x_q,r}$ such that the induced representation of $\ZG_{\ov \bF_q}$ on $\ZV_{\ov \bF_q}$ corresponds to the usual adjoint representation of $\RP_{x_q}$ on $\PV_{x_q,r}$. Moreover, there are isomorphisms $\ZG_{\fpb} \simeq \RP_{x}$ and $\ZV_{\fpb} \simeq \PV_{x,r}$ such that the induced representation of $\ZG_{\fpb}$ on $\ZV_{\fpb}$ is the usual adjoint representation of $\RP_{x}$ on $\PV_{x,r}$. 
	In other words, we have commutative diagrams
		\begin{eqnarray*}  
			\xymatrix{
				\ZG_{\fpb} \times  \ZV_{\fpb} \ar[r] \ar[d]^{\simeq \times \simeq} &   \ZV_{\fpb} \ar[d]^{\simeq}
				& & \ZG_{\fqb} \times  \ZV_{\fqb} \ar[r] \ar[d]^{\simeq \times \simeq} &   \ZV_{\fqb} \ar[d]^{\simeq} \\
				\RP_x \times  \PV_{x,r}   \ar[r]    &  \PV_{x,r} 
				&  & \RP_{x_q} \times  \PV_{x_q,r}   \ar[r]    &  \PV_{x_q,r} \quad .
			} 
		\end{eqnarray*}
\end{Thm}

\begin{Rem}  \label{Remark-failure}
The existence of a reductive group $G_q$ over $K_q$ and a point $x_q \in \sB(G_q,K_q)$ satisfying the conditions of the above theorem fails for some reductive groups $G$ that are not good groups. For example, let $K$ be a maximal unramified extension of $\bQ_2$,  $E=K(\sqrt{-1})$, and $G$ the corresponding norm one torus, i.e. the kernel of the norm map from $\Res_{E/K}\bG_m$ to $\bG_m$. Then $\sB(G,K)$ consists of only one point $x$, the reductive quotient $\RP_x$ is trivial, and \linebreak % careful with this line break when chaging to a different latex style!!!!!
 $\PV_{x,r} \simeq \left\{\begin{array}{ll} \ov \bF_p & \text{ if } r \in \bZ  \\ \{ 0 \} & \text{ if } r \in \bR-\bZ\end{array} \right.$. However, for $q>2$, there does not exist a reductive group $G_q$ over a finite extension $K_q$ of $\bQ_q^{ur}$ and $x_q \in \sB(G_q,K_q)$ so that the above theorem holds. (Sketch of the argument: Assume such a group $G_q$ exists. Since the reductive quotient is trivial, $G_q$ has to be anisotropic. Since $\sum_{s \leq r <s+1} \dim \PV_{x_q,r}=\sum_{s \leq r <s+1} \dim \PV_{x,r}=1$ for any $s \in \bR$, the group $G_q$ has to be a one dimensional torus, hence $G_q$ has to be the norm one torus of a quadratic extension $E_q$ of $K_q$. However, this implies $\PV_{x_q,r} \simeq \left\{\begin{array}{ll} \ov \bF_q & \text{ if } r \in \frac{1}{2}+ \bZ  \\ \{ 0 \} & \text{ if } r \in \bR-(\frac{1}{2}+\bZ)\end{array}\right.$.)
\end{Rem}

We prove the theorem in two steps. In Section \ref{Section-global-reductive-quotient} we construct a morphism from $\ZG$ to an auxiliary split reductive group scheme $\ZZG$, and in Section \ref{Section-global-filtration-quotient} we construct $\ZV$ (largely) inside the Lie algebra  of $\ZZG$ and use the adjoint action of $\ZZG$ on its Lie algebra to define the action of $\ZG$ on $\ZV$.  

\subsubsection{Global reductive quotient} \label{Section-global-reductive-quotient}
Let $\ZZG$ \index{notation}{$\ZZG$} be a split reductive group scheme over $\bZ$ whose root datum is the root datum of $G$. In this section we construct a morphism $\iota: \ZG \ra \ZZG$ that lifts all the morphisms $\iota_{K,F}: G_{x,0}/G_{x,0+} \hookrightarrow \GF_{x,0}/\GF_{x,0+}$ and $\iota_{K,\Fq}: G_{x_q,0}/G_{x_q,0+} \hookrightarrow \GFq_{x_q,0}/\GFq_{x_q,0+}$ defined in Section \ref{section-iota}. In order to do so, let us first describe the image of $\iota_{K,F}$ more explicitly. In analogy to the root group parametrization $x_a$ defined in Section \ref{Sec-rootgroups}, and using the notation from that section, we define for $a \in \Phi_K(G)$ multipliable the more general map $X_a: F \times F \ra G(F)$ by
\begin{equation*}
	X_a(u,v)=\prod_{\beta \in [\Phi_a]}	\xE_\beta(u_\beta)\xE_{\beta + \wt \beta}(-v_\beta)\xE_{\wt \beta}(\sigma(u)_\beta),
\end{equation*}
where $\sigma$ denotes an element of $\Gal(F/E_{\alpha+\wt \alpha})$ that projects to the nontrivial element of $\Gal(E_\alpha/E_{\alpha+\wt \alpha})$ and where $u_\beta=\gamma(u)$ for some fixed choice of $\gamma \in \Gal(F/K)$ with $\gamma(\alpha)=\beta$. Note that $X_a|_{H_0(E_\alpha,E_{\alpha+\wt \alpha})}$ ($\alpha \in \Phi_a$) agrees with $x_a$. We then have the following lemma.

\begin{Lemma} \label{Lemma-xbar}
	Let $\chi: \fpb \ra \cO_{\bQ_p^{ur}}$ (if $\bQ_p \subset F$) or $\chi: \fpb \ra \cO_{\bF_p((t))^{ur}}$ (if $\bF_p((t)) \subset F$)  be the Teichmüller lift, and $\RU_a$ the root group of $\RP_x$ corresponding to the root $a \in \Phi(\RP_x) \subset \Phi_K(G)$. Define the map $y_a: \fpb \ra G_{x,0}^F$ by
	{ \small \begin{eqnarray*} 
		 u   & \mapsto  & \left\{ \begin{array}{ll} 
		 		{X_a(\sqrt{2}\chi(u) \cdot \varpi_F^{-a(x-x_0)\cdot M},\chi(u)\varpi_F^{-a(x-x_0)\cdot M}\sigma(\chi(u)\varpi_F^{-a(x-x_0)\cdot M}))} & \mbox{if $a$ is multipliable and $p \neq 2$} \\ 
		 		{X_a(0,\chi(u)\sigma(\chi(u))\varpi_F^{-2a(x-x_0)\cdot M}}) & \mbox{if $a$ is multipliable and $p=2$} \\
		 		X_a(0,\chi(u) \cdot \varpi_F^{-2a(x-x_0)\cdot M}) & \mbox{if $a$ is divisible} \\
		 		{x_a(\chi(u) \cdot \varpi_F^{-a(x-x_0)\cdot M})} & \mbox{otherwise.} 
		 		 \end{array} \right. 
	\end{eqnarray*} }

Then the composition $\ov y_a$ of $y_a$ with the quotient map $G_{x,0}^F \twoheadrightarrow G_{x,0}^F/G_{x,0+}^F$ is isomorphic to $\iota_{K,\Fp} \circ \Rx_a : \fpb \ra \iota_{K,\Fp}(\RU_a(\fpb)) \subset \RP_x^F(\fpb)$.
\end{Lemma}

\textbf{Proof.} 
If $p \neq 2$ or if $a$ is not multipliable, the lemma follows immediately from Lemma \ref{Lemma-Chevalley}.

In the case $p=2$, note that (using the notation from Lemma \ref{Lemma-Chevalley})
$$ \va\left(\chi(u)\varpi_\Fp^{s'}\sigma(\chi(u)\varpi_\Fp^{s'})\cdot \varpi_\Fp^{\va(\lambda)M}\right) < 2\va\left(\sqrt{1/\lambda_0}\chi(u) \varpi_\Fp^{{s'}}\right), $$
where ${s'}=-(a(x-x_0)+\va(\lambda)/2)M$, because $\va(\lambda)<0$. Moreover, $\sigma(\varpi_\Fp) \equiv \varpi_\Fp \mod \varpi_F^2$ in $\varpi_F\cO_F/\varpi_F^2 \cO_F$, and hence $\ov y_a(u)= \iota_{K,\Fp}(\Rx_a(u))$ by Lemma \ref{Lemma-Chevalley}. \qed

\begin{Rem} \label{Rem-xbar}
 An analogous statement holds for $\RP_{x_q}$. In the sequel 
 we denote the root group parameterizations constructed for $\RP_{x_q}$ analogously to Lemma \ref{Lemma-Chevalley} by ${\Rxq}_a: \bG_a \ra {\RUq}_a, a \in \Phi(\RP_{x_q})$\index{notation}{${\Rxq}_a$}.
\end{Rem}

Recall that $x$ is hyperspecial in $\sB(G_\Fp,\Fp)$, and hence the reductive quotient $\RPF_x$ of $G_\Fp$ at $x$ is a split reductive group over $\fpb$ with root datum $R(\RPF_x)=R(G)$. The analogous statement holds for $x_q$. Thus $\ZZG_{\fpb}$ is isomorphic to $\RPF_x$, and $\ZZG_{\fqb}$ is isomorphic to $\RPFq_{x_q}$. In order to construct explicit isomorphisms, let us fix a split maximal torus $\ZZT$\index{notation}{$\ZZT$} of $\ZZG$ and a Chevalley system $\{ \ZZx_\alpha: \bG_a \xra{\simeq} \ZZU_\alpha \subset \ZZG \}_{\alpha \in \Phi(\ZZG)=\Phi}$ for $(\ZZG, \ZZT)$\index{notation}{$\ZZx_\alpha$}\index{notation}{$\ZZU_\alpha$} with signs equal to $\eps_{\al,\beta}$ as in Definition \ref{Def-signs}; i.e., the Chevalley system $\{\ZZx_\al\}_{\al \in \Phi}$ for $(\ZZG,\ZZT)$ and the Chevalley--Steinberg system $\{x_\al\}_{\al \in \Phi}$ for $(G,T)$ have the same signs. 

Moreover, the split  maximal torus $T_\Fp \subset G_\Fp$ and the Chevalley system $\{ \xF_\alpha\}_{\alpha \in \Phi}$ yield a  split maximal torus $\RTF_x$ of $\RPF_x$ and a Chevalley system $\{ \ov \xF_\alpha: \bG_a \xra{\simeq} \RUF_\alpha \subset \RPF_x \}_{\alpha \in \Phi}$\index{notation}{$\ov \xF_\alpha$} for $(\RPF_x,\RTF_x)$ with signs $\eps_{\alpha,\beta}$. Similarly, we obtain a split maximal torus $\RTFq_{x_q}$ of $\RPFq_{x_q}$ and a Chevalley system $\{ \ov \xFq_\alpha: \bG_a \xra{\simeq} \RUFq_\alpha \subset \RPFq_{x_q} \}_{\alpha \in \Phi}$ for $(\RPFq_{x_q},\RTFq_{x_q})$ with signs $\eps_{\alpha,\beta}$. 
In addition, we denote by $\RT_x$ \index{notation}{$\RT_x$} and $\RT_{x_q}$ the maximal split tori of $\RP_x$ and $\RP_{x_q}$ corresponding to $S$ and $S_q$.

Moreover, we define constants $c_{\al,q} \in \cO_{\Fq}$ and $c_\alpha \in \cO_{\Fp}$ for $\alpha \in \Phi$ as follows. 
We choose $\gamma' \in \Gal(F/K)$ such that 
$$ \gamma' \mod \Gal(\Fp/ E \cap K^{\rm{tame}}) \equiv u \quad \in \Gal(E \cap K^{\rm{tame}}/K) $$
and $\zetaG \in \cO_K$ \index{notation}{$\zetaG$}
 satisfying
$$ \gamma'(\varpi_\Fp) \equiv \zetaG \varpi_\Fp \mod \varpi_\Fp^2 .$$
Similarly, let $\gamma_q \in \Gal(\Fq/K_q) \simeq \bZ/M\bZ$ correspond to $1\in \bZ/M\bZ$, i.e.
$$ \gamma_q \mod \Gal(\Fq/ \Eq) \equiv {\gammapr} \quad \in \Gal(\Eq/K)$$
 and $\zetaGq \in \cO_{K_q}$ \index{notation}{$\zetaGq$}
 such that
$$ \gamma_q(\varpi_\Fq) = \zetaGq \varpi_\Fq.$$

Let $C_1, \hdots, C_n$ be the representatives for the action of $\Gammapr=\<\gammapr\>$ on the connected components of $\Dyn(G)$ as given in Remark \ref{Rem-good}\ref{item-root-action-prime}, and recall that $\Phi_i$ denotes the roots that are a linear combination of simple roots corresponding to $C_i$.
For $\alpha \in \Phi$ there exists a unique triple $(i,\alpha_i, e_q(\al))$ with $i \in [1,n]$, $\alpha_i \in \Phi_i$ and $e_q(\alpha)$ minimal in $\bZ_{\geq 0}$ such that  $\gamma_q^{e_q(\al)}(\al_i)=\al$. Note that  $e_q(\al)$ is independent of the choice of prime number $q$. We also write $e(\al)=e_q(\al)$\index{notation}{$e(\al)$}. 
We define \index{notation}{$c_{\alpha,q}$}\index{notation}{$c_{\alpha}$}
 $$ c_{\alpha,q}:=\zetaGq^{e(\al) \cdot\alpha_i(x_q-x_{0,q})\cdot M}=\zetaGq^{e(\al) \cdot\alpha(x_q-x_{0,q})\cdot M} \quad \text{and} \quad  c_\alpha:=\zetaG^{e(\al) \cdot \alpha_i(x-x_0)\cdot M}=\zetaG^{e(\al) \cdot \alpha(x-x_0)\cdot M}.$$
 Note that $\alpha_i(x-x_0)\cdot M$ is an integer, as the order $m$ of $x$ divides $M$ and $\Gamma_a' \subset \va(F)=\frac{1}{M}\bZ$, where $a$ is the image of $\alpha$ in $\Phi_K$.

Finally, we denote by $\ov \zetaG$\index{notation}{$\ov \zetaG$} and $\ov \zetaGq$\index{notation}{$\ov \zetaGq$} the images of $\zetaG$ and $\zetaGq$  and by $\ov c_\alpha$\index{notation}{$\ov c_{\alpha}$} and $\ov c_{\al, q}$\index{notation}{$\ov c_{\alpha,q}$} the images of $c_\alpha$ and $c_{\al, q}$ under the surjections $\cO_{\Fp} \twoheadrightarrow \fpb$ and $\cO_{\Fq} \twoheadrightarrow \fqb$, respectively.

\begin{Rem} \label{Rem-zeta} 
The integers $e(\alpha)$ depend only on the connected component of $\Dyn(G)$ in whose span $\al$ lies.
\end{Rem}

The definitions of $\zetaG$, $\zetaGq$ and $e({\alpha})$ 
 are chosen so that the following lemma holds. 
\begin{Lemma} \label{Lemma-action-on-varpi} We keep the notation from above and let $r \in \bR$.
	\begin{enumerate}[(i)]
		\item \label{item-ealpha-i} If $\wt \gamma \in \Gal(\Fq/K_q)$ with $\wt \gamma(\alpha_i)=\alpha$ and $r':=r-\al(x_q-x_{0,q}) \in \Gamma_a'(G_q)$, then 
	$$\wt \gamma(\varpi_{\Fq}^{r'M})  \equiv  \zetaGq^{e(\al) \cdot(r-\alpha(x_q-x_{0,q}))M} \varpi_{\Fq}^{r'M}  \mod \varpi_{\Fq}^{r'M+1} .$$  
	\item \label{item-ealpha-ii} If $\wt \gamma \in \Gal(\Fp/K)$ with $\wt \gamma(\alpha_i)=\alpha$ and $r':=r-\al(x-x_{0}) \in \Gamma_a'(G)$, then 
	$$\wt \gamma(\varpi_{F}^{r'M}) \equiv  \zetaG^{e(\al) \cdot(r-\alpha(x-x_{0}))M} \varpi_{F}^{r'M}  \mod \varpi_{F}^{r'M+1} .$$
	\end{enumerate}
\end{Lemma}
\textbf{Proof.} 
If $\wt \gamma \in \Gal(\Fq/K_q)$ with $\wt \gamma(\alpha_i)=\alpha$, then $\wt \gamma=\gamma_q^{e(\alpha)+z\abs{\<\gammapr\>(\al_i)}}$ for some integer $z$. As $r'\in \Gamma_a'(G_q)=\frac{1}{\abs{\<\gammapr\>(\al_i)}}\bZ$, we have $\ov \zetaGq^{\abs{\<\gammapr\>(\al_i)}r'M}=1$ and 
$$\wt \gamma(\varpi_{\Fq}^{r'M})  \equiv \gamma_q^{e(\alpha)+z\abs{\<\gammapr\>(\al_i)}}(\varpi_{\Fq}^{r'M})  \equiv  \zetaGq^{e(\al) r'M} \varpi_{\Fq}^{r'M} \equiv \zetaGq^{e(\al) \cdot(r-\alpha(x_q-x_{0,q}))M} \varpi_{\Fq}^{r'M}  \mod \varpi_{\Fq}^{r'M+1}, $$
which shows part \ref{item-ealpha-i}.

In order to prove part \ref{item-ealpha-ii}, let $\wt \gamma \in \Gal(\Fp/K)$ with $\wt \gamma(\alpha_i)=\alpha$, and write $\wt \gamma =  {\gamma'}^{\wt e} \wt w$ for some integer $\wt e$ and $\wt w \in \Gal(F/E\cap K^{\rm{tame}})$. By Property \ref{item-root-action} of Definition \ref{Def-good} and the definition of $e(\al)$ there exists $w \in \Gal(F/E \cap K^{\rm{tame}})$ such that ${\gamma'}^{e(\alpha)} w (\al_i)=\al$, and hence $\wt w^{-1} {\gamma'}^{e(\alpha)-\wt e} w (\al_i)=\al_i$, and therefore ${(\gammapr)}^{e(\al)-\wt e}(\al_i)\in \Gal(F/\Et)(\al_i)$. On the other hand, as the $\Gammapr$-orbits on $\Phi$ agree with  the $\Gal(F/K)$-orbits on $\Phi$ and $X^{{\gammapr}^N}=X^{\Gal(F/E \cap K^{\rm{tame}})}$ (by Property \ref{item-Et} of Definition \ref{Def-good} and Lemma \ref{Lemma-good}), the $\Gal(F/E \cap K^{\rm{tame}})$-orbits on $\Gal(F/K)(\alpha_i)$ coincide with the $\langle{\gammapr}^{N}\rangle$-orbits, which are the same as the $\langle{\gammapr}^{N_i}\rangle$ orbits, where $N_i$ is coprime to $p$ such that $\abs{\Gal(F/K)(\alpha_i)}=p^{s_i}N_i$ for some integer $s_i$. Thus $e(\al)-\wt e \equiv 0 \mod N_i$. Note that $\ov \zetaG^{N_ir'M}=1$ in $\fpb$, because
$r' \in \Gamma'_a(G) = \frac{1}{p^{s_i}N_i} \bZ $ if $p\neq 2$ and $r' \in \Gamma'_a(G) \subset \frac{1}{2p^{s_i}N_i} \bZ $ if $p = 2$.
Moreover, for $g \in \Gal(F/\Et)$, $g(\varpi_F) \equiv \varpi_F \mod \varpi_F^2$ as all $p$-power roots of unity in $\fpb$ are trivial. Hence 
$$\wt \gamma(\varpi_{F}^{r'M})  \equiv {\gamma'}^{\wt e} (\varpi_{F}^{r'M}) \equiv \zetaG^{\wt e \cdot r'M} \varpi_{F}^{r'M}
 \equiv \zetaG^{e(\al) \cdot(r-\alpha(x_q-x_{0,q}))M} \varpi_{F}^{r'M}  \mod \varpi_{F}^{r'M+1}, $$
which proves part \ref{item-ealpha-ii}. \qed

Now let $f_T:\RTF_x \ra \ZZT_{\fpb}$ be an isomorphism that identifies the root data $R(\RPF_x)$ and $R(\ZZG)$. Then we can extend $f_T$ as follows.

\begin{Lemma} \label{Lemma-f}
	There exists an isomorphism $f: \RPF_x \ra \ZZG_{\fpb}$\index{notation}{$f$} extending $f_T$ such that for $\alpha \in \Phi$ and $u \in \bG_a(\fpb)$ we have 
	\begin{equation} \label{equation-rootgroupmap}
	 f(\ov \xF_\alpha(u))=\ZZx_\alpha(\ov c_\alpha \cdot u).
	 \end{equation}
\end{Lemma}
\textbf{Proof.} Note that there exists a unique isomorphism $f: \RPF_x \ra \ZZG_{\fpb}$ extending $f_T$ and satisfying Equation \eqref{equation-rootgroupmap} for all $\alpha \in \Delta$. 
So we need to show that this $f$ satisfies Equation \eqref{equation-rootgroupmap} for all $\alpha \in \Phi$. In order to do so, it suffices to show that the root group parameterizations $\{{\ZZx_\alpha}_{\fpb} \circ \ov c_\alpha\}_{\alpha \in \Phi}$ form a Chevalley system of $(\ZZG_{\fpb}, \ZZT_{\fpb})$ whose signs $\eps_{\alpha,\beta}'$  are equal to $\eps_{\alpha,\beta}$ ($\alpha, \beta \in \Phi$), i.e. to the signs of $\{{\ov \xF_\alpha}\}_{\alpha \in \Phi}$. 
 If $\alpha$ and $\beta$ are linear combinations of roots in different connected components of the Dynkin diagram of $\Phi$, then $\eps_{\alpha,\beta}'=1=\eps_{\alpha,\beta}$. Thus suppose $\alpha, \beta \in {\gamma'}(\Phi_1)$, and hence also $s_\alpha(\beta) \in {\gamma'}(\Phi_1)$, for some ${\gamma'} \in \Gal(\Fp/K)$. By Remark \ref{Rem-zeta} this implies that $\ov \zeta_{\gamma'}:=\ov \zetaG^{e(\alpha)}= \ov \zetaG^{e(\beta)} = \ov \zetaG^{e(s_\al(\beta))}$. 
 We obtain (using \cite[Cor.~5.1.9.2]{ConradSGA3} for the second equality) 
 \begin{eqnarray*}
 	& &\Ad\left(\ZZx_\alpha(\ov c_\alpha)\ZZx_{-\alpha}(\eps_{\al,\al}\ov c_{-\alpha})\ZZx_\alpha(\ov c_\alpha)\right)\left(\Lie({\ZZx_\beta}_{\fpb} \circ \ov c_\beta)(1)\right) \\
 	& =&\Ad\left(\ZZx_\alpha(\ov \zeta_{\gamma'}^{\alpha(x-x_0)\cdot M})\ZZx_{-\alpha}(\eps_{\al,\al}\ov \zeta_{\gamma'}^{-\alpha(x-x_0)\cdot M})\ZZx_\alpha(\ov \zeta_{\gamma'}^{\alpha(x-x_0)\cdot M})\right)\left(\ov \zeta_{\gamma'}^{\beta(x-x_0)\cdot M} \Lie({\ZZx_\beta}_{\fpb})(1)\right) \\
 	& =&\Ad\left(\check \alpha(\ov \zeta_{\gamma'}^{\alpha(x-x_0)\cdot M})\right)\Ad\left(\ZZx_\alpha(1)\ZZx_{-\alpha}(\eps_{\al,\al})\ZZx_\alpha(1)\right)\left(\ov \zeta_{\gamma'}^{\beta(x-x_0)\cdot M}\Lie({\ZZx_\beta}_{\fpb})(1)\right) \\
 	& =&\ov \zeta_{\gamma'}^{\beta(x-x_0)\cdot M} \cdot \Ad\left(\check \alpha(\ov \zeta_{\gamma'}^{\alpha(x-x_0)\cdot M})\right)\left(\eps_{\alpha,\beta}\Lie({\ZZx_{s_\alpha(\beta)}}_{\fpb})(1)\right) \\
  	& =&\ov \zeta_{\gamma'}^{\beta(x-x_0)\cdot M} \cdot (s_\alpha(\beta))(\check \alpha(\ov \zeta_{\gamma'}^{\alpha(x-x_0)\cdot M}))\cdot \eps_{\alpha,\beta}\Lie({\ZZx_{s_\alpha(\beta)}}_{\fpb})(1) \\
  	& =&\ov \zeta_{\gamma'}^{\beta(x-x_0)\cdot M} \cdot \ov \zeta_{\gamma'}^{\<\check \alpha, s_\alpha(\beta)\>\cdot\alpha(x-x_0)\cdot M}\cdot \eps_{\alpha,\beta}\Lie({\ZZx_{s_\alpha(\beta)}}_{\fpb})(1) \\
  	& =&\ov \zeta_{\gamma'}^{(s_\alpha(\beta))(x-x_0)\cdot M} \eps_{\alpha,\beta}\Lie({\ZZx_{s_\alpha(\beta)}}_{\fpb})(1) \\
  	& =& \eps_{\alpha,\beta}\left(\Lie({\ZZx_{s_\alpha(\beta)}}_{\fpb}\circ \ov c_{s_\alpha(\beta)})(1)\right). 
 \end{eqnarray*}
Thus the signs of the Chevalley system  $\{{\ZZx_\alpha}_{\fpb} \circ \ov c_\alpha \}_{\alpha \in \Phi}$  are $\eps_{\alpha,\beta}$ as desired. \qed

Similarly, for each prime $q$, let $f_{T,q}:\RTFq_{x_q} \ra \ZZT_{\fqb}$ be an isomorphism that identifies the root data $R(\RPF_{x_q})$ and $R(\ZZG)$. 
 Then we have the analogous statement.
\begin{Lemma} \label{Lemma-fq}
	There exists an isomorphism $f_q: \RPF_{x_q} \ra \ZZG_{\fqb}$ \index{notation}{$f_q$} extending $f_{T,q}$ such that for $\alpha \in \Phi$ and $u \in \bG_a(\fqb)$ we have 
	\begin{equation} \label{equation-rootgroupmapq}
		f_q(\ov \xFq_\alpha(u))=\ZZx_\alpha(\ov c_{\alpha,q} \cdot u).
	\end{equation}
\end{Lemma}

This allows us to define a map $\iota$ from $\ZG$ to $\ZZG$ as follows. 

Let $\ZS$ \index{notation}{$\ZS$} be a split maximal torus of $\ZG$. Then we have 
\begin{equation*}
	X_*(\ZS)=X_*(\RT_x)=X_*(S)=X_*(T)^{\Gal(\Fp/K)}\hookrightarrow X_*(T)=X_*(\ZZT),
\end{equation*}
where the first identification arises from $R(\ZG)=R(\RP_x)$, the second from Lemma \ref{lemma-reductive-quotient} and the fourth from $R(\ZZG)=R(G)$.
This yields a closed immersion $f_\ZS:\ZS \ra \ZZT$.
Note that $f_\ZS$ also corresponds to the injection
\begin{equation*}
	X_*(\ZS)=X_*(\RT_{x_q})=X_*(S_q)=X_*(T_q)^{\Gal(\Fq/K_q)}\hookrightarrow X_*(T_q)=X_*(\ZZT),
\end{equation*}
and we have commutative diagrams
\begin{eqnarray*} 
		\xymatrixcolsep{4pc}\xymatrix{
			\ZS_{\fpb}  \ar[r]^{f_\ZS} \ar@{-}[d]_{\simeq} &  \ZZT_{\fpb} & \ZS_{\fqb}  \ar[r]^{f_\ZS} \ar@{-}[d]_{\simeq} &  \ZZT_{\fqb}   \\
			\RT_x  \ar[r]_{\iota_{K,\Fp}}    & \RTF_{x}  \ar@{-}[u]_{\simeq} & \RT_{x_q}  \ar[r]_{\iota_{K_q,\Fq}}    & \RTFq_{x_q}  \ar@{-}[u]_{\simeq}\quad .
		} 
\end{eqnarray*}

To define $\iota$ on root groups, let $\{ \Zx_a\}_{a \in \Phi(\ZG)=\Phi(\RP_x)}$ \index{notation}{$\Zx_a$} be a Chevalley system for $(\ZG,\ZS)$ such that there exists an isomorphism $f_{H,q}: \ZG_{\fqb} \ra \RP_{x_q}$ mapping $\ZS_{\fqb}$ to  $\RT_{x_q}$ and identifying $(\Zx_a)_{\fqb}$ with ${\Rxq}_a$, or equivalently having the same signs as the Chevalley system $\{ {\Rxq}_a\}_{a \in \Phi_K}$, for some $q \neq 2$.

Moreover, note that for $a \in \Phi_K=\Phi(\ZG)$, there exists a unique integer in $[1,n]$, denoted by $n(a)$\index{notation}{$n(a)$}, such that $\Phi_a \cap \Phi_{n(a)} \neq \emptyset$ (see Remark \ref{Rem-good} for the definition of $\Phi_i, i \in [1,n]$). We label the elements in $\Phi_a \cap \Phi_{n(a)}$ by $\{\alpha_i\}_{1 \leq i \leq \abs{\Phi_a \cap \Phi_{n(a)}}}$ so that they satisfy the following two properties:
\begin{itemize}%[label=$\cdot$]
\item  If $a$ is a multipliable root, we assume that $\alpha_1 \in [\Phi_a]$, where $[\Phi_a]$ is as defined in Section \ref{Sec-rootgroups}. (Note that a priori we have either $\alpha_1$ or $\alpha_2$  in $[\Phi_a]$.) 
\item  Let $\gammapr$  be the generator of $\Gammapr$ as in Definition \ref{Def-good}, then for all $a \in \Phi_K$ with $\abs{\Phi_a \cap \Phi_{n(a)}}=3$, there exists a minimal integer $e'(a)$ such that ${\gammapr}^{e'(a)}$ preserves and acts nontrivially on $\Phi_a \cap \Phi_{n(a)}$, and we require that ${\gammapr}^{e'(a)}(\alpha_1)=\alpha_2$. (Note that this implies ${\gammapr}^{e'(a)}(\al_2)=\al_3$.)
\end{itemize}

We may (and do) assume that $[\Phi_a]$ is chosen to be $\{{\gammapr}^i(\alpha_1) \, | \, 0 \leq i \leq{\abs{\Phi_a}}-1 \}$.

\begin{DefProp} \label{Prop-iota} 
	There exists a unique group scheme homomorphism $\iota: \ZG_{\ov \bZ} \ra \ZZG_{\ov \bZ}$\index{notation}{$\iota$} extending $f_\ZS$ such that for all $\ov \bZ$-algebras $A$, $a \in \Phi(\ZG)=\Phi_K$ and $u \in \bG_a(A)$ we have
	\begin{flalign} 
		& \iota(\Zx_a(u))  =  \prod\limits_{i =1}^{\abs{\Gammapr/\Gammapr_{n(a)}}} \ZZx_{{\gammapr}^{(i-1)}(\alpha_1)}(\sqrt 2 u) \ZZx_{{\gammapr}^{(i-1)}(\alpha_1+\alpha_2)}(-(-1)^{-a(x-x_0)M}u^2)\ZZx_{{\gammapr}^{(i-1)}(\alpha_2)}((-1)^{-a(x-x_0)M} \sqrt 2 u) &  \label{eqn-DefProp1}
	\end{flalign}
	if $a$ is multipliable,
	\begin{flalign} &  \iota(\Zx_a(u))  =  \prod\limits_{i =1}^{\abs{\Gammapr/\Gammapr_{n(a)}}} \ZZx_{{\gammapr}^{(i-1)}(\alpha_1)}(-u) & &	\text{if $a$ is divisible, and} &  \label{eqn-DefProp2} \\
	& \iota(\Zx_a(u))  = \prod\limits_{i =1}^{\abs{\Gammapr/\Gammapr_{n(a)}}} \prod_{j=1}^{\abs{\Phi_a \cap \Phi_{n(a)}}} \ZZx_{{\gammapr}^{(i-1)}(\alpha_j)}(\zeta_{\abs{\Phi_a \cap \Phi_{n(a)}}}^{-a(x-x_0)M(j-1)} u)  & & \label{eqn-DefProp3} \text{otherwise,} 
	\end{flalign}
	
 where $\zeta_i$ is a primitive $i$-th root of unity, $i=1,2$ or $3$, and $\Gammapr_{n(a)}=\Stab_{\Gammapr}(\Phi_{n(a)})$. 
	
	Moreover, there exist unique isomorphisms $f_{\ZG}:\RP_x \xra{\simeq} \ZG_{\ov \bF_p}$\index{notation}{$f_{\ZG}$} and $f_{\ZG, q}:\RP_{x_q} \xra{\simeq} \ZG_{\ov \bF_q}$\index{notation}{$f_{\ZG, q}$} such that we have commutative diagrams
	\begin{eqnarray*} 
		\xymatrixcolsep{4pc}\xymatrix{
			\ZG_{\fpb}  \ar[r]^{\iota} \ar[d]_{\simeq}^{f_{\ZG}} &  \ZZG_{\fpb} & \ZG_{\fqb}  \ar[r]^{\iota} \ar[d]_{\simeq}^{f_{\ZG, q}} &  \ZZG_{\fqb}   \\
			\RP_x  \ar[r]_{\iota_{K,\Fp}}    & \RPF_{x}  \ar[u]_{\simeq}^{f} & \RP_{x_q}  \ar[r]_{\iota_{K_q,\Fq}}    & \RPFq_{x_q}  \ar[u]_{\simeq}^{f_q}\quad 
		} 
	\end{eqnarray*}
	for all primes $q$.
\end{DefProp}
\textbf{Proof.}
Combining Lemma \ref{Lemma-xbar} and Remark \ref{Rem-xbar} with Lemma \ref{Lemma-f} and Lemma \ref{Lemma-fq}, we observe in view of Property \ref{item-root-action-prime} of Remark \ref{Rem-good} and Lemma \ref{Lemma-action-on-varpi} that $f \circ \iota_{K,\Fp} \circ \Rx_a$ and $f_q \circ \iota_{K_q,\Fq} \circ {\Rxq}_a$ are described by the (reduction of the) right hand side of the three equations in the definition / proposition for all primes $q$. As $\iota_{K_q,\Fq} \circ {\Rxq}_a$ (and $\iota_{K,\Fp} \circ \Rx_a$ ) are isomorphisms from $\bG_a$ to $\iota_{K_q,\Fq} (\RUq_a)$ (and $\iota_{K,\Fp} (\RU_a)$) for $q \neq 2$ (and for $p \neq 2$),  the signs of the Chevalley systems $\{{\Rxq}_a\}_{a \in \Phi_K}$ coincide with those of $\{\Rx_a\}$ and of $\{\Zx_a\}$ for all $q$. (Note that $1 = -1$ in characteristic two, i.e. the previous statement is trivial in this case.) This implies for every prime $q$ the existence of a unique isomorphism $f_{\ZG, q}:\RP_{x_q} \xra{\simeq} \ZG_{\fqb}$ that identifies $\RT_{x_q}$ with $\ZS_{\fqb}$ and ${\Rxq}_a$ with $(\Zx_a)_{\fqb}$ for all $a \in \Phi_K$, and similarly for $\RP_x$.

Note that the Equations \eqref{eqn-DefProp1}, \eqref{eqn-DefProp2} and \eqref{eqn-DefProp3} in the definition / proposition define group scheme homomorphisms $f_a:\bG_a \ra \ZZG_{\ov \bZ}$ over $\ov \bZ$ for $a \in \Phi(\ZG)$. The maps $\{f_a\}_{a \in \Delta(\ZG)}$ and $f_\ZS$ together with the requirement that $\Zx_a(1)\Zx_{-a}(\eps_{a,a})\Zx_{a}(1) \mapsto f_a(1)f_{-a}(\eps_{a,a})f_a(1)$ for $a \in \Delta(\ZG)$ define by \cite[XXIII, Theorem~3.5.1]{SGA3III} a unique group scheme homomorphism $\iota: \ZG_{\ov \bZ} \ra \ZZG_{\ov \bZ}$. (The required relations asked for in \cite[XXIII, Theorem~3.5.1]{SGA3III} can be checked to be satisfied using that they hold in $\fqb$ for all primes $q$ by the existence of $\iota_{K_q,\Fq}$ (similar to the subsequent argument).)

We are left to check that the Equations \eqref{eqn-DefProp1}, \eqref{eqn-DefProp2} and \eqref{eqn-DefProp3} hold for $a \in \Phi-\Delta(\ZG)$. For this note that $\iota(\Zx_{s_b(a)}(\eps_{b,a}u))=\left(f_b(1)f_{-b}(\eps_{b,b})f_b(1)\right)\iota(\Zx_a(u))\left(f_b(1)f_{-b}(\eps_{b,b})f_b(1)\right)^{-1}$ for $a \in \Phi, b \in \Delta(\ZG)$, where $\{\eps_{a,b}\}_{a,b \in \Phi_K}$ are the signs of the Chevalley system $\{\Zx_a\}_{a \in \Phi_K}$. For $a, b \in \Delta(\ZG)$, the trueness of the equations in the proposition for $s_b(a)$ for all $u \in \bG_a(A)$ is therefore equivalent to the vanishing of a finite number of polynomials with coefficients in $\ov \bZ$. As the latter vanish mod $q$ for all primes $q$, these polynomials vanish also over $\ov \bZ$, and the equations are satisfied for $s_b(a)$ ($b, a \in \Delta(\ZG))$, and hence by repeating the argument for all roots $a \in \Phi$.
\qed

\begin{Rem}
	The morphism $\iota$ can be defined over $\bZ[x]/(x^3-1)=\bZ[\zeta_3]$ or even over $\bZ$ if none of the connected components of $\Dyn(G)$ is of type $D_4$ with vertices contained in only two orbits.
\end{Rem}

In order to provide a different construction of $\ZG$ in Section \ref{section-Vinberg}, we use the following Lemma.

\begin{Lemma} \label{Lemma-closed-immersion}
	 Let $\iota$ be as in Definition / Proposition \ref{Prop-iota}. Then $\iotaQ: \ZG_{\ov \bQ} \ra \ZZG_{\ov \bQ}$ is a closed immersion.
\end{Lemma}
\textbf{Proof.}  \\
In order to show that $\iotaQ$ is a closed immersion, it suffices to show that its kernel is trivial (\cite[Proposition~1.1.1]{ConradSGA3}). As $\ov \bQ$ is of characteristic zero, the kernel of $\iotaQ$ (a group scheme of finite type) is smooth. Hence we only need to show that $\iotaQ$ is injective on $\ov \bQ$-points. Let $g \in \ZG(\ov \bQ)$. Let $\dot{W}$ be a set of representatives of the Weyl group of $\ZG$ in the normalizer of $\ZS$. Without loss of generality, we assume that the elements of $\dot{W}$ are products of $\Zx_a(1)\Zx_{-a}(\eps_{a,a})\Zx_{a}(1)$, $a \in \Delta(\ZG)$, or the identity. Let $U$ be the unipotent radical of the Borel subgroup corresponding to $\Delta(\ZG)$, $U^-$ the one of the opposite Borel corresponding to $-\Delta(\ZG)$, and $U_w=U(\ov \bQ) \cap wU^-(\ov \bQ)w^{-1}$. By the Bruhat decomposition, we can write $g$ uniquely as $u_1wtu_2$ with $w \in \dot{W}$, $t \in \ZS(\ov \bQ)$, $u_1 \in U_w$ and $u_2\in U(\ov \bQ)$. By the uniqueness $1=\iota(g)=\iota(u_1)\iota(w)\iota(t)\iota(u_2)$ if and only if $1=\iota(u_1)=\iota(w)=\iota(t)=\iota(u_2)$. Note that $\iota(w)=1$ implies $w=1$ by our choice of $\dot{W}$, and $\iota(t)=1$ implies $t=1$. Choosing an order of the positive roots of $\Phi_K^+$, there is a unique way to write $u_2=\prod_{a \in \Phi_K^+}\Zx_a(u_a)$ with $u_a \in \ov \bQ$ for all $a \in \Phi_K^+$. By choosing a compatible ordering of the roots in $\Phi^+$ and the uniqueness of writing $\iota(u_2)=\prod_{\alpha \in \Phi^+}\ZZx_\alpha(u'_\alpha)$ with $u'_\alpha \in \ov \bQ$ together with the explicit description of $\iota$ on root groups given in Definition / Proposition \ref{Prop-iota}, we conclude that $u_a=0$ for all $a \in \Phi_K^+$, and hence $u_2=1$. Similarly, $u_1=1$, which shows that the map $\iota$ is injective as desired. \qed

\subsubsection{Global Moy--Prasad filtration quotients} \label{Section-global-filtration-quotient}
In this section we will also lift the injections $\iota_{K, \Fp, r}: \PV_{x,r} \ra \PVF_{x,r}$
and $\iota_{K_q, \Fq, r}: \PV_{x_q,r} \ra \PVFq_{x_q,r}$ in such a way that we get a lift of the commutative Diagram \eqref{diagram-iota}. Using these injections we view $\PV_{x,r}$ as a subspace of $\PVF_{x,r}$ and $\PV_{x_q,r}$ as a subspace of $\PVFq_{x_q,r}$. 
We will afterwards modify the global action slightly to also accomodate the case where $p=2$ and there exists $a \in \Phi_K^{\mathrm{mul}}$ with $a(x-x_0) \in \Gamma_a'$ such that $a(x-x_0)-r \in \Gamma_a'$ or such that there exists $b \in \Phi_K^{\rm{nm}}$ with $b(x-x_0)-r \in \Gamma_{b}'$ and $\<\check a, b\> \neq 0$.

We begin with the construction of an integral model for $\PV_{x_q,r}$.
Fix $r \in \va(\Fp)=\va(\Fq)$ (otherwise the Diagram \eqref{diagram-iota} would be trivial)
 and let $\zeta_M$ be a primitive $M$-th root of unity in $\ov \bZ$ compatible with $\zeta_3$ in Proposition \ref{Prop-iota}, i.e. if $3 \, | \, M$, then $\zeta_M^{M/3}=\zeta_3$. Let $\vartheta$\index{notation}{$\vartheta$} denote the composition of the action of $\gammapr$ on $\Lie(\ZZT)(\ov \bZ[1/N])$ induced from its action on $R(\ZZG)=R(G)$ (as given by Definition \ref{Def-good}), and multiplication by $\zeta_M^{rM}$, and define $\ZV_T$\index{notation}{$\ZV_T$} to be the free $\ov \bZ[1/N]$-submodule of $\Lie(\ZZT)(\ov \bZ[1/N])$ fixed by $\vartheta$.

Next consider $a \in \Phi_K$. We 
 recall that $\Gammapr_{n(a)}$ denotes the stabilizer of the component $C_{n(a)}$ of the Dynkin diagram $\Dyn(G)$ inside $\Gammapr$, and set $\ZZX_\alpha=\Lie(\ZZx_\alpha)(1) \in \Lie(\ZZG)(\ov \bZ[1/N])$ \index{notation}{$\ZZX_\alpha$} for $\al \in \Phi$. We define \index{notation}{$Y_a$}
\begin{equation} \label{eqn-DefYa}
	 Y_a
	 =\sum_{i=1}^{\abs{\Phi_a \cap \Phi_{n(a)}}}\sum_{j=1}^{\abs{\Gammapr/\Gammapr_{n(a)}}}{\zeta_M}^{e(\gammapr(\al_1))(j-1)rM}\zeta_{\abs{\Phi_a \cap \Phi_{n(a)}}}^{(-a(x_q-x_{0,q})+r)\abs{\Gammapr/\Gammapr_{n(a)}}\abs{\Phi_a \cap \Phi_{n(a)}}(i-1)}\ZZX_{{\gammapr}^{(j-1)}(\alpha_i)} 
\end{equation}
(note that $\zeta_{\abs{\Phi_a \cap \Phi_{n(a)}}}^{(-a(x_q-x_{0,q})+r)\abs{\Gammapr/\Gammapr_{n(a)}}\abs{\Phi_a \cap \Phi_{n(a)}}(i-1)} \in \{1, -1, \zeta_3, \zeta_3^2\}$)
and let $\ZVold$\index{notatoin}{$\ZVold$} be the free $\ov \bZ[1/N]$-submodule of $\Lie(\ZZG)(\ov \bZ[1/N])$ generated by $\ZV_T$ and $Y_a$ for all $a \in \Phi_K$ with $r-a(x_q-x_{0,q}) \in \Gamma_a'(G_q)$, or equivalently $r-a(x-x_0) \in \Gamma_a'(G)$ by Lemma \ref{lemma-affine-roots}.
Note that $\ZVold$ as a $\ov \bZ[1/N]$-module is a direct summand of the free $\ov \bZ[1/N]$-module $\Lie(\ZZG)(\ov \bZ[1/N])$.

Also note that the $\RPF_{x}$ representation $\PVF_{x,r}$ is isomorphic to the adjoint representation of $\RPF_{x}$ on $\Lie(\RPF_{x})$ and, similarly, the $\RPFq_{x_q}$ representation $\PVFq_{x_q,r}$ is isomorphic to the adjoint representation of $\RPFq_{x_q}$ on $\Lie(\RPFq_{x_q})$. Hence the isomorphisms $f:\RPF_x \xra{\simeq} \ZZG_{\fpb}$ and  $f_q:\RPFq_{x_q} \xra{\simeq} \ZZG_{\fqb}$ from Lemma \ref{Lemma-f} and \ref{Lemma-fq} yield isomorphisms $df:=\Lie(f):\PVF_{x,r} \simeq \Lie(\RPF_x)(\fpb) \xra{\simeq} \Lie(\ZZG)({\fpb})$ and  $df_q:=\Lie(f_q):\PVFq_{x_q,r} \xra{\simeq} \Lie(\ZZG)({\fqb)}$.

\begin{Prop} \label{Prop-generala}
	The adjoint action of $\ZZG_{\ov \bZ[1/N]}$ on $\Lie(\ZZG)(\ov \bZ[1/N])$ restricts to an action of $\ZG_{\ov \bZ[1/N]}$ on $\ZVold$. 
	
	Let $q$ coprime to $N$. Then we have $df(\PV_{x,r})=\ZVold_{\fpb}$ and $df_q(\PV_{x_q,r})=\ZVold_{\fqb}$. Moreover, the following diagrams commute
	\begin{eqnarray*}  
		\xymatrix{
			\ZG_{\fpb} \times  \ZVold_{\fpb} \ar[r] \ar[d]^{\simeq}_{f_{\ZG}^{-1} \times d{f}^{-1}} &   \ZVold_{\fpb} \ar[d]_{d{f}^{-1}}^{\simeq}
			& & \ZG_{\fqb} \times  \ZVold_{\fqb} \ar[r] \ar[d]^{\simeq}_{f_{\ZG, q}^{-1} \times d{f}_q^{-1}} &   \ZVold_{\fqb} \ar[d]^{\simeq}_{d{f}_q^{-1}} \\
			\RP_x \times  \PV_{x,r}   \ar[r]    &  \PV_{x,r} 
			&  & \RP_{x_q} \times  \PV_{x_q,r}   \ar[r]    &  \PV_{x_q,r} \quad ,
		} 
	\end{eqnarray*}
		unless $p$ or $q$ is 2 (for the left or right diagram, respectively) and there exists $a \in \Phi_K^{\rm{mul}}$ with $a(x-x_0) \in \Gamma_a'$  such that $a(x-x_0)-r \in \Gamma_a'$ or such that there exists $b \in \Phi_K^{\rm{nm}}$ with $b(x-x_0)-r \in \Gamma_{b}'$ and $\<\check a, b\> \neq 0$.
\end{Prop}
\textbf{Proof.}  
We first show that $df_q(\PV_{x_q,r})=\ZVold_{\fqb}$ for $q$ coprime to $N$ and $df(\PV_{x,r})=\ZVold_{\fpb}$ by considering the intersection of $\ZVold$ with the subspaces $\bigoplus_{\alpha \in \Phi(G)} \Lie(\ZZG)(\ov \bZ[1/N])_\alpha$ and $\Lie(\ZZT)(\bZ[1/N])$ of $\Lie(\ZZG)(\bZ[1/N])$ separately.

For $\al \in \Phi$, 
 denote by $\Gammapr_\alpha$ the stabilizer of $\alpha$ in $\Gammapr$, and let $\RX_\alpha=\Lie(\ov{\xFq}_\alpha)(1)$, $n_a=\abs{\Phi_a\cap\Phi_{n(a)}} \in \{1,2,3\}$ and $\zeta_{\gammapr}:=\ov\zetaGq^{e(\gammapr(\alpha_1))}=\ov\zetaGq^{e(\gammapr(\alpha_i))}$, $1 \leq i \leq n_a$.\\
  The image of $\left(\ZVold \cap \bigoplus_{\alpha \in \Phi(G)} \Lie(\ZZG)(\ov \bZ[1/N])_\alpha\right) \otimes_{\ov \bZ[1/N]} \fqb$ under $df_q^{-1}$ is then spanned by 
\begin{eqnarray*}
	\ov Y_{a} &= & \sum_{i=1}^{n_a}\sum_{j=1}^{\abs{\Gammapr/\Gammapr_{n(a)}}}{\zeta_{\gammapr}}^{(j-1)rM}\zeta_{n_a}^{(-a(x_q-x_{0,q})+r)\abs{\Gammapr/\Gammapr_{n(a)}}n_a(i-1)}\ov c_{{\gammapr}^{(j-1)}(\alpha_i),q}^{-1} \RX_{{\gammapr}^{(j-1)}(\alpha_i)}  \\
	& = & \sum_{i=1}^{n_a}\sum_{j=1}^{\abs{\Gammapr/\Gammapr_{n(a)}}}{\zeta_{\gammapr}}^{(j-1)rM}\zeta_{n_a}^{(-a(x_q-x_{0,q})+r)\abs{\Gammapr/\Gammapr_{n(a)}}n_a(i-1)} \zeta_{\gammapr}^{-\alpha(x_q-x_{0,q})M(j-1)} \RX_{{\gammapr}^{(j-1)}(\alpha_i)}  \\
	& = & \sum_{j=1}^{\abs{\Gammapr/\Gammapr_\alpha}}{\zeta_{\gammapr}}^{(j-1)(r-a(x_q-x_{0,q}))M} \RX_{{\gammapr}^{(j-1)}(\alpha_1)} \\
	&=& \sum_{j=1}^{\abs{\Gammapr/\Gammapr_\alpha}}{\gammapr}^{(j-1)} \left(\RX_{\alpha_1}\right)  
\end{eqnarray*}
 for $a \in \Phi_K$ with $r-a(x_q-x_{0,q}) \in \Gamma_a'(G_q)$
(where $\zeta_M$ gets send to $\ov \zetaGq$ under the surjection $\ov \bZ[1/N] \twoheadrightarrow \fqb$). Here the action of $\Gammapr$ on $\PVFq_{x_q,r}$ is the one induced from the action on $\fgFq_{x_q,r}$. Thus by definition of the Moy--Prasad filtration and the inclusion $\iota_{\Fq,K_q,r}$ constructed in the proof of Lemma \ref{Lemma-iota} we obtain the equality 
\begin{eqnarray}  df_q^{-1}(\ZVold_{\fqb}) \cap \bigoplus_{\alpha \in \Phi}\Lie(\RPFq_{x_q})(\fqb)_\alpha 
&= & df_q^{-1}((\ZVold \cap \bigoplus_{\alpha \in \Phi}\Lie(\ZZG)(\ov \bZ[1/N])_\alpha)\otimes \fqb) \notag \\
&=& \PV_{x_q,r} \cap \bigoplus_{\alpha \in \Phi}\Lie(\RPFq_{x_q})(\fqb)_\alpha \label{eqn-rootpart-q}
\end{eqnarray} 
inside $\PVFq_{x_q,r} \simeq \Lie(\RPFq_{x_q})(\fqb)$. 

In order to show the analogous statement for $V_{x,r}$, we claim that $\zeta_{n_a}^{(-a(x_q-x_{0,q})+r)\abs{\Gammapr/\Gammapr_{n(a)}}n_a(i-1)}=\zeta_{n_a}^{(-a(x-x_0)+r)\abs{\Gammapr/\Gammapr_{n(a)}}n_a(i-1)}$ in $\ov \bF_p$. This is obviously true for $p \neq 2$ as $a(x-x_0)=a(x_q-x_{0,q})$ in this case. If $p=2$, then $\zeta_{2}=-1=1$ in $\fpb$ and we only have to consider the case $n_a=\abs{\Phi_a\cap \Phi_{n(a)}}=3$. However, $n_a=3$ implies that the corresponding component $C_{n(a)}$ of the Dynkin diagram $\Dyn(G)$ is of type $D_4$, and hence $\check b(a)=0$ for all multipliable roots $b \in \Phi_K^{+,\rm{mul}}$. Thus $a(x-x_0)=a(x_q-x_{0,q})$ by definition, see Equation \eqref{equation-f-def}, and the claim $\zeta_{n_a}^{(-a(x_q-x_{0,q})+r)\abs{\Gammapr/\Gammapr_{n(a)}}n_a(i-1)}=\zeta_{n_a}^{(-a(x-x_{0,q})+r)\abs{\Gammapr/\Gammapr_{n(a)}}n_a(i-1)}$ follows. Let $\zeta_{\gammapr}=\ov\zetaG^{e(\gammapr(\alpha_1))}$, $\RX_\alpha=\Lie(\ov{\xF}_\alpha)(1)$, and use otherwise the same notation as above. Then there exists a set of representatives $\left[\Gal(\Fp/K)/\Stab_{\Gal(\Fp/K)}(\alpha)\right]$ of $\Gal(\Fp/K)/\Stab_{\Gal(\Fp/K)}(\alpha)$ such that the image of $\left(\ZVold \cap \bigoplus_{\alpha \in \Phi(G)} \Lie(\ZZG)(\ov \bZ[1/N])_\alpha\right)  \otimes_{\ov \bZ[1/N]} \fpb$ under $df^{-1}$ is spanned by 
\begin{eqnarray*}
	\ov Y_a &= & \sum_{i=1}^{n_a}\sum_{j=1}^{\abs{\Gammapr/\Gammapr_{n(a)}}}{\zeta_{\gammapr}}^{(j-1)rM}\zeta_{n_a}^{(-a(x_q-x_{0,q})+r)\abs{\Gammapr/\Gammapr_{n(a)}}n_a(i-1)}\ov c_{{\gammapr}^{(j-1)}(\alpha_i)}^{-1} \RX_{{\gammapr}^{(j-1)}(\alpha_i)}  \\
	& = & \sum_{i=1}^{n_a}\sum_{j=1}^{\abs{\Gammapr/\Gammapr_{n(a)}}}{\zeta_{\gammapr}}^{(j-1)rM}\zeta_{n_a}^{(-a(x-x_{0})+r)\abs{\Gammapr/\Gammapr_{n(a)}}n_a(i-1)} \zeta_{\gammapr}^{-\alpha(x-x_{0})M(j-1)} \RX_{{\gammapr}^{(j-1)}(\alpha_i)}  \\
	&= & \sum_{j=1}^{\abs{\Gammapr/\Gammapr_\alpha}}{\zeta_{\gammapr}}^{(j-1)(r-a(x-x_0))M} \RX_{{\gammapr}^{(j-1)}(\alpha_1)} = \sum_{{\gamma'} \in \left[\Gal(\Fp/K)/\Stab_{\Gal(\Fp/K)}(\alpha)\right]}{\gamma'} \left(\RX_{\alpha_1}\right)  ,
\end{eqnarray*}
where the last equality follows from 
  Lemma \ref{Lemma-action-on-varpi}. Thus we obtain
\begin{equation} \label{eqn-rootpart-p} df^{-1}(\ZVold_{\fpb}) \cap \bigoplus_{\alpha \in \Phi}\Lie(\RPF_{x})(\fpb)_\alpha= \PV_{x,r} \cap \bigoplus_{\alpha \in \Phi}\Lie(\RPF_{x})(\fpb)_\alpha 
\end{equation}  
inside $\PVF_{x,r} \simeq \Lie(\RPF_{x})(\fpb)$.

Let us consider $\ZV_T$. From the definition of the Moy--Prasad filtration $\ftE_{x,r}$ of the Lie algebra $\ft_{E_t}$ of the torus $T_{E_t}$ together with Lemma \ref{Lemma-good} and the observation that all $p$-power roots of unity in $\ov \bF_p$ are trivial,  we deduce (by sending $\zeta_M \otimes 1$ to $\ov \zetaG$ under the isomorphism $\ov \bZ[1/N] \otimes_{\ov \bZ[1/N]} \fpb \simeq \fpb$, as above) that 
$$df\left(\iota_{E_t,F,r}\left(\ftE_{x,r}/\ftE_{x,r+}\right)\right) = (\Lie(\ZZT)(\ov \bZ[1/N]))^{\vartheta^{N}} \otimes_{\ov \bZ[1/N]}  \fpb .$$ 
 Moreover, by combining 4.6.2.~Proposition and Proposition~4.6.1 from \cite[Section~4.6]{Yu2}, we have $\ft_{x,r}=(\ftE_{x,r})^{\Gal(E_t/K)}$ as $E_t$ is tamely ramified over $K$, and we obtain (using tameness of $E_t/K$) that 
\begin{eqnarray}
df\left(\ft_{x,r}/\ft_{x,r+}\right)&=&df\left((\ftE_{x,r}/\ftE_{x,r})^{\Gal(E_t/K)}\right)
=\left(\left(\Lie(\ZZT)(\ov \bZ[1/N])\right)^{\vartheta^{N}} \otimes_{\ov \bZ[1/N]} \fpb \right)^\vartheta \notag \\
&=&(\Lie(\ZZT)(\ov \bZ[1/N]))^{\vartheta}\otimes_{\ov \bZ[1/N]} \fpb=\ZV_T \otimes_{\ov \bZ[1/N]} \fpb . \label{eqn-toruspart-p}
\end{eqnarray}

For $q$ coprime to $N$, we denote by $E_{t,q}$ the tamely ramified extension of degree $N$ of $K_q$. Then we obtain by the same reasoning (substituting $E_t$ by $E_{t,q}$)
\begin{eqnarray} \label{eqn-toruspart-q}
	df_q\left(\ft_{x_q,r}/\ft_{x_q,r+}\right)&=& \ZV_T \otimes_{\ov \bZ[1/N]} \fqb.
\end{eqnarray}

Combing Equations \eqref{eqn-rootpart-p} and \eqref{eqn-toruspart-p}, and \eqref{eqn-rootpart-q} and \eqref{eqn-toruspart-q}, we obtain for $q$ coprime to $N$ that
$$ df(\PV_{x,r})=\ZVold_{\fpb} \quad \text{ and } df_q(\PV_{x_q,r})=\ZVold_{\fqb} .$$

In order to show that the adjoint action of $\ZZG_{\ov \bZ[1/N]}$ on $\Lie(\ZZG)(\ov \bZ[1/N])$ restricts to an action of $\ZG_{\ov \bZ[1/N]}$ on $\ZVold$, we observe that the following diagram commutes
	\begin{eqnarray*}  
		\xymatrix{
			 & \ZZG_{\fqb} \times  \Lie(\ZZG)(\ov \bZ[1/N])_{\fqb} \ar[r] \ar[d]^{\simeq}_{f_q^{-1}  \times d{f}_q^{-1}} &   \Lie(\ZZG)(\ov \bZ[1/N])_{\fqb} \ar[d]^{\simeq}_{d{f}_q^{-1}} \\
			  & \RPFq_{x_q} \times  \PVFq_{x_q,r}   \ar[r]    &  \PVFq_{x_q,r} \quad .
		} 
	\end{eqnarray*}
Since $\iota_{K_q,\Fq}(\RP_{x_q})$ preserves $\PV_{x_q,r}$ (Lemma \ref{Lemma-iota}), we deduce that the induced action of $\ZG_{\fqb}$ on $\Lie(\ZZG)(\fqb)$  preserves $\ZVold_{\fqb}$ for all $q$ coprime to $N$. Hence the induced action of $\ZG_{\ov \bZ[1/N]}$ on $\Lie(\ZZG)(\ov \bZ[1/N])$ preserves $\ZVold$, and by construction and Lemma \ref{Lemma-iota} and Definition / Proposition \ref{Prop-iota} the diagrams in the proposition commute (assuming the condition in the proposition in characteristic 2). \qed

In order to also obtain commutative diagrams in the case when $p$ or $q$ is 2 and there exists $a \in \Phi_K^{\rm{mul}}$ with $a(x-x_0) \in \Gamma_a'$, we define the $\ov\bZ[1/N]$-submodule $\ZV$\index{notation}{$\ZV$} of $\Lie(\ZZG)(\ov\bZ[1/N])$ to be the submodule generated by $\ZV_T$, $\ZVnm$ and $\ZVmul$, where $\ZVnm$ is the $\ov\bZ[1/N]$-submodule generated by $Y_a$ for all $a \in \Phi_K^{\mathrm{nm}}$ with $r-a(x-x_0) \in \Gamma_a'(G)$ and  $\ZVmul$ is the $\ov\bZ[1/N]$-submodule generated by $\sqrt{2}Y_a$ for all $a \in \Phi_K^{\mathrm{mul}}$ with $r-a(x-x_0) \in \Gamma_a'(G)$. Note that $\ZV$ is a finite index submodule of $\ZVold$, and the injection $\ZV \hookrightarrow \ZVold$ yields an isomorphism $\ZV  \otimes \ov \bZ[1/2N] \xrightarrow{\simeq} \ZVold  \otimes \ov \bZ[1/2N]$.

\begin{Lemma} \label{Lemma-V}
	Let $R$ be a $\ov\bZ[1/N]$-algebra. Then the image of $\ZV \otimes R$ in $\ZVold \otimes R$ is preserved by the action of $\ZG(R)$.
\end{Lemma}
\textbf{Proof.} To simplify notation, we assume $R=\ov\bZ[1/N]$, but the proof is the same for general $R$. We need to show that $\ZG(R)$ maps $\ZV_T \oplus \ZVnm$ to $\ZV_T \oplus \ZVnm \oplus \ZVmul$. Since $\ZS$ preserves $\ZV_T \oplus \ZVnm$ it suffices to consider the action of the root groups $\Zx_a(R)$ for ${a \in \Phi(\ZG)=\Phi(\RP_x) \subset \Phi_K(G)}$. Let $a \in \Phi(\ZG)\subset \Phi_K(G)$. If $\alpha \in \Phi_a$ is not contained in the span of roots of a connected component of the Dynkin diagram $\Dyn(G)$ that is of type $A_{2n}$ and on which $\Gal(E/K)$ acts non-trivially, then $\Zx_a(R)$ preserves  $\ZV_T \oplus \ZVnm$. By the same reasoning as in the proof of Lemma \ref{Lemma-iota}, if $a$ corresponds to a non-multipliable root in $\Phi_K(G)$, then $\Zx_a(R)$ preserves  $\ZV_T \oplus \ZVnm$ as well. Thus assume $a$ is multipliable. Hence, for $u\in \bG_a(R)$, we have by Definition / Proposition \ref{Prop-iota} that
$$ \iota(\Zx_a(u))  =  \prod\limits_{i =1}^{\abs{\Gammapr/\Gammapr_{n(a)}}} \ZZx_{{\gammapr}^{(i-1)}(\alpha_1)}(\sqrt 2 u) \ZZx_{{\gammapr}^{(i-1)}(\alpha_1+\alpha_2)}(-(-1)^{-a(x-x_0)M}u^2)\ZZx_{{\gammapr}^{(i-1)}(\alpha_2)}((-1)^{-a(x-x_0)M} \sqrt 2 u) . $$
Let $H \in \ZV_T$. Using that $\ZZx_{\al}(u)(H) =  H - \Lie(\al)(H)u\ZZX_{\al}$ for all $\alpha \in \Phi$, we observe that $\Zx_a(u)(H) =  \iota(\Zx_a(u))(H)$ is contained in $\ZV_T \oplus \ZVnm \oplus \ZVmul$.

It remains to consider the action of $\Zx_a(u)$ on $Y_b$ for $b \in \Phi_K^{\mathrm{nm}}$ with $r-b(x-x_0) \in \Gamma_b'(G)$. Let us assume (without influence on the arguments to follow) that $\alpha_1$ and $\alpha_2$ above are the simple roots $\alpha_1$ and $\beta_1$ of a Dynkin diagram of type $A_{2n}$ as in Figure \ref{fig-A2n1} on page \pageref{fig-A2n1}. Then $\Zx_a(u)(Y_b)=Y_b$ unless $b$ is the restriction of $\alpha_2 + \hdots + \alpha_t$ or of $-(\beta_1+\alpha_1 + \hdots + \alpha_t)$ for some $2 < t \leq n$ using the notation from Figure \ref{fig-A2n1}. In both cases we observe using the explicit formulas for $\iota(\Zx_a(u))$ and $Y_b$ that $\Zx_a(u)(Y_b) =  \iota(\Zx_a(u))(Y_b)$ is contained in $\ZV_T \oplus \ZVnm \oplus \ZVmul$. \qed

The lemma allows us to define an action of $\ZG$ on $\ZV$ by requiring that if $R$ is an $\ov\bZ[1/N]$-algebra in which $2 \neq 0$, then the action of $\ZG(R)$ on $\ZV_R$ is the restriction of the action of $\ZG(R)$ on $\ZVold_R$. Note that if $N$ is odd, then for $\ov g \in \ZG(\ov \bF_2)$ and $\ov v \in \ZV(\ov \bF_2)$ there exist $g \in \ZG(\ov \bZ[1/N])$ whose image in $\ZG(\ov \bF_2)$ is $\ov g$ (because this holds for the root groups and the torus) and $v \in  \ZV(\ov \bZ[1/N])$ whose image in $\ZV(\ov \bF_2)$ is $\ov v$, and $\ov g. \ov v$ is the image of $g.v \in \ZV(\ov \bZ[1/N])$ in $\ZV(\ov \bF_2)$ (which is independent on the choice of $g$ and $v$).

Note that the action of $\ZG$ on $\ZV  \otimes \ov \bZ[1/2N]$ corresponds to the action of  $\ZG$ on $\ZVold  \otimes \ov \bZ[1/2N]$ under the identification $\ZV  \otimes \ov \bZ[1/2N] \xrightarrow{\simeq} \ZVold  \otimes \ov \bZ[1/2N]$ above. In order to treat the special fiber over $\ov \bF_2$, we define isomorphisms $\fZV: \PV_{x,r} \ra \ZV_{\ov\bF_p}$ if $p=2$ and ${\fZV}_{, 2}: \PV_{x_2,r} \ra \ZV_{\ov\bF_2}$ if $2 \nmid N$ as follows.

Suppose $p=2$. Let $\fZV|_{\overline{\fg_{x,r} \cap (\ft \oplus \bigoplus_{a \in \Phi_K^{\mathrm{nm}}}\fg_a)}}: \overline{\fg_{x,r} \cap (\ft \oplus \bigoplus_{a \in \Phi_K^{\mathrm{nm}}}\fg_a)} \xra{\simeq} (\ZV_T)_{\ov \bF_2} \oplus (\ZVnm)_{\ov\bF_2}$ be given by the restriction of $df$, and let 
$\fZV(\overline{Y}_a)= \overline{\sqrt{\lambda_0} \sqrt{2}Y_a}$ for $a \in \Phi_K^{\mathrm{mul}}$ with $r-a(x-x_0) \in \Gamma_a'(G)$, where $\lambda_0$ is as defined in Lemma~\ref{Lemma-Chevalley},  $\overline{\sqrt{\lambda_0} \sqrt{2} Y_a}$ denotes the image of $ \sqrt{\lambda_0} \sqrt{2} Y_a \in \ZV$ under the surjection $\ZV \ra \ZV \otimes \ov\bF_2$ and $\overline{Y}_a$ is as introduced in the proof of Proposition \ref{Prop-generala}, i.e. 
 $\ov{Y}_a=\sum_{{\gamma'} \in \left[\Gal(\Fp/K)/\Stab_{\Gal(\Fp/K)}(\alpha)\right]}{\gamma'} \left(\RX_{\alpha_1}\right)$ with the above notation.

Define the isomorphism ${\fZV}_{, 2}: \PV_{x_2,r} \ra \ZV_{\ov\bF_2}$ analogously.

\begin{Prop}\label{Prop-generalb}
	For $q$ coprime to $N$, the following diagrams commute 
		\begin{eqnarray*}  
			\xymatrix{
				\ZG_{\fpb} \times  \ZV_{\fpb} \ar[r] \ar[d]^{\simeq}_{f_{\ZG}^{-1} \times {\fZV}^{-1}} &   \ZV_{\fpb} \ar[d]_{{\fZV}^{-1}}^{\simeq}
				& & \ZG_{\fqb} \times  \ZV_{\fqb} \ar[r] \ar[d]^{\simeq}_{f_{\ZG, q}^{-1}  \times {\fZV}_{,q}^{-1}} &   \ZV_{\fqb} \ar[d]^{\simeq}_{{\fZV}_{,q}^{-1}} \\
				\RP_x \times  \PV_{x,r}   \ar[r]    &  \PV_{x,r} 
				&  & \RP_{x_q} \times  \PV_{x_q,r}   \ar[r]    &  \PV_{x_q,r} \quad ,
			} 
		\end{eqnarray*}
		where $\fZV:=df$ if $p \neq 2$ and ${\fZV}_{, q}:=df_q$ for $q \neq 2$.
\end{Prop}
\textbf{Proof.} By Proposition \ref{Prop-generala} and the above observation that $\ZV  \otimes \ov \bZ[1/2N] \hookrightarrow \ZVold  \otimes \ov \bZ[1/2N]$ is an isomorphism of $\ZG_{\ov \bZ[1/2N]}$-modules, the right diagram commutes for $q \neq 2$ and the left diagram commutes if $p \neq 2$.

Let us consider the commutativity of the left diagram in the case that $p=2$. The commutativity of the right diagram for $q=2$ follows from the same arguments.

By construction, the action of the maximal torus $\RT_x$ on $\PV_{x,r}$ corresponds to the action of $\ZS_{\ov \bF_2}$ on $\ZV_{\ov \bF_2}$, and it remains to consider the action the root groups $\RU_a \subset \RP_x$ for $a \in \Phi(\RP_x) \subset \Phi_K$. 
We first consider the action on $\overline{\fg_{x,r} \cap (\ft \oplus \bigoplus_{a \in \Phi_K^{\mathrm{nm}}}\fg_a)} \simeq (\ZV_T)_{\ov \bF_2} \oplus (\ZVnm)_{\ov\bF}$.
In the proofs of Lemma \ref{Lemma-iota} and Lemma \ref{Lemma-V}, we have seen that if $a$ is non-multipliable, then $\RU_a=\Rx_a(\bG_m)$ preserves $\overline{\fg_{x,r} \cap (\ft \oplus \bigoplus_{a \in \Phi_K^{\mathrm{nm}}}\fg_a)}$ and $\Zx_a((\bG_m)_{\ov \bF_2})$ preserves $(\ZV_T)_{\ov \bF_2} \oplus (\ZVnm)_{\ov\bF}$, and hence, by construction, the actions agree under the isomorphisms $f_\ZG|_{\RU_a}$ and $\fZV|_{\overline{\fg_{x,r} \cap (\ft \oplus \bigoplus_{a \in \Phi_K^{\mathrm{nm}}}\fg_a)}}$.

So consider $a$ multipliable, and let $\overline u \in \ov\bF_2$.
Then 
$$\Rx_a(\ov u)(\ov X)=\overline{x_a\left(\sqrt{\frac1{\lambda_0}}\chi(u) \varpi_E^{s} \eps_1, \, \chi(u)\varpi_E^s\eps_1\sigma(\chi(u)\varpi_E^s\eps_1)\cdot \varpi_E^{\va(\lambda)e}\eps_0\right)\left(X\right)}$$
for $X \in \fg_{x,r}$, where we use the notation from Lemma \ref{Lemma-Chevalley} and $\ov{?}$ denotes the image of $?$ in $\fg_{x,r}/\fg_{x,r+}=\PV_{x,r}$.
On the other hand, if $u \in \ov \bZ[1/N]$ maps to $\ov u \in \ov\bF_2$, then
$$\Zx_a(\ov u)(\ov X)=\overline{\iota(\Zx_a(u))(X)}$$ 
for $X \in \ZV$ and where $\ov{?}$ denotes the image of $?$ in $\ZV_{\ov \bF_2}$.
Moreover, by Definition / Proposition \ref{Prop-iota}
$$ \iota(\Zx_a(u))  =  \prod\limits_{i =1}^{\abs{\Gammapr/\Gammapr_{n(a)}}} \ZZx_{{\gammapr}^{(i-1)}(\alpha_1)}(\sqrt 2 u) \ZZx_{{\gammapr}^{(i-1)}(\alpha_1+\alpha_2)}(-(-1)^{-a(x-x_0)M}u^2)\ZZx_{{\gammapr}^{(i-1)}(\alpha_2)}((-1)^{-a(x-x_0)M} \sqrt 2 u) .$$
Using these equations and the description of $x_a$ in Equation \eqref{equation-xa}, easy calculations show that
$\fZV\left(\Rx_a(\ov u)(H)\right)=\Zx_a(\ov u)(\fZV(H))=f_\ZG(\Rx_a(\ov u))(\fZV(H))$
 for $H \in \overline{\fg_{x,r} \cap \ft}$ and
$\fZV\left(\Rx_a(\ov u)(\ov{Y}_b)\right)=\Zx_a(\ov u)(\ov{Y_b})=f_\ZG(\Rx_a(\ov u))(\fZV(\ov{Y}_b))$  
for $b \in \Phi_K^{\mathrm{nm}}$ with $r-b(x-x_0) \in \Gamma'_b(G)$.

It remains to consider the action on $\overline{\fg_{x,r} \cap \bigoplus{a \in \Phi_K^{\mathrm{mul}} \fg_a}} \xra{\simeq} (\ZVmul)_{\ov \bF_2}$. By Lemma \ref{Lemma-iota} (and the definition of $\iota_{K,F,r}$ in its proof) and the definition of $\ZV$ and $\ZVmul$, the groups $\RP_x$ and $\ZG_{\ov \bF_2}$ preserve $\overline{\fg_{x,r} \cap \bigoplus_{a \in \Phi_K^{\mathrm{mul}} \fg_a}}$ and $(\ZVmul)_{\ov \bF_2}$, respectively. Hence, by the underlying constructions, their action agrees under the isomorphisms $f_\ZG$ and $\fZV|_{\overline{\fg_{x,r} \cap \bigoplus_{a \in \Phi_K^{\mathrm{mul}} \fg_a}}}$.
\qed

Now Theorem \ref{Thm-general} is an immediate consequence of Proposition \ref{Prop-generalb}.

\section{Moy--Prasad filtration representations and global Vinberg--Levy theory} \label{section-Vinberg}

In this section we will give a different description of the reductive group scheme $\ZG$ and its action on $\ZV$ from Theorem \ref{Thm-general} as a fixed-point group scheme of a larger split reductive scheme $\ZZG$ acting on a graded piece of $\Lie \ZZG$ (see Theorem \ref{Thm-tame}). This means we are in the setting of a global version of Vinberg--Levy theory and the special fibers correspond to (generalized) Vinberg--Levy representations for all primes $q$. 
In order to give such a description integrally (i.e. over $\ov \bZ[1/N]$), we will specialize to reductive groups $G$ that become split over a tamely-ramified field extension in Section \ref{section-tame}. Afterwards, in Section \ref{section-VinbergQ}, we will then show that such a description holds over $\ov \bQ$ for all \good{} groups. This will also allow us to study the existence of (semi)stable vectors in Section \ref{section-semistablestablevectors}.

\subsection{The case of $G$ splitting over a tamely ramified extension} \label{section-tame}
Let $S$ be a scheme, then we denote by $\mmu_{M,S}$\index{notation}{$\mmu_{M,S}$} the group scheme of $M$-th roots of unity over $S$. We will often omit $S$ if it can be deduced from the context. Given an $S$-group scheme  $\ZZG$, we denote by $\underline\Aut_{\mbox{\small$\ZZG$}/S}$
 its automorphisms functor, i.e. the functor that sends an $S$-scheme $S'$ to the group of automorphisms of $\ZZG_{S'}$ in the category of $S'$-group schemes, and by $\Aut_{\mbox{\small$\ZZG$}/S}$ 
 its representing group scheme if it exists.  We will often omit $S$ if it can be deduced from the context.  Given, in addition, a morphism $\theta: \mmu_{M,S} \ra \Aut_{\mbox{\small$\ZZG$}}$, we denote by $\ZZG^\theta$ \index{notation}{$\ZZG^\theta$} the scheme theoretic fixed locus of $\ZZG$ under the action of $\mmu_{M,S}$ via $\theta$, if it exists, i.e. $\ZZG^\theta$ represents the functor that sends an $S$-scheme $S'$ to the elements of $\ZZG(S')$ on which $\mmu_{M,S'}$ acts trivially. If $\ZZG^\theta$ is a smooth group scheme over $S$ of finite presentation, we denote by $\ZZG^{\theta,0}$ \index{notation}{$\ZZG^{\theta,0}$} its identity component. Similarly, if $\cF$ is  a quasi-coherent $\cO_S$-module, we denote by $\underline \Aut_{\mbox{\small$\cF$}/\cO_S}$ 
 its automorphism functor, and by $\Aut_{\mbox{\small$\cF$}/\cO_S}$ 
  (or simply $\Aut_{\mbox{\small$\cF$}}$) the group scheme representing $\underline \Aut_{\mbox{\small$\cF$}/\cO_S}$ if it exists.

\begin{Thm} \label{Thm-tame}
	Suppose that $G$ is a reductive group over $K$ that splits over a tamely ramified field extension $E$ of degree $e$ over $K$. Let $r=\frac{d}{M}$ for some nonnegative integer $d < M$, and let $\ZG$ be the split reductive group scheme over $\ov \bZ[\frac{1}{e}]$ acting on the free $\ov \bZ[\frac{1}{e}]$-module $\ZV$ as provided by Theorem \ref{Thm-general}, i.e. such that the special fibers each correspond to the action of a reductive quotient on a Moy--Prasad filtration quotient. 
	Then there exists a split reductive group scheme $\ZZG$ defined over $\ov\bZ[\frac{1}{e}]$ and morphisms 
	$$\theta: \mmu_{M} \ra \Aut_{\mbox{\small$\ZZG$}} \quad \text{ and } \quad d\theta:  \mmu_{M} \ra \Aut_{\mbox{\small$\Lie(\ZZG)$}}$$ that induces a $\bZ/{M} \bZ$-grading $\Lie(\ZZG)= \oplus_{i=1}^{M}\Lie(\ZZG)_i$ such that $\ZG$ is isomorphic to $\ZZG^{\theta,0}$, $\ZV$ is isomorphic to  $\Lie(\ZZG)_{M-d}(\ov\bZ[\frac{1}{e}])$ and the action of $\ZG$ on $\ZV$ corresponds to the restriction of the adjoint action of $\ZZG$ on $\Lie(\ZZG)(\ov\bZ[\frac{1}{e}])$ via these isomorphisms.

		In particular, this implies that for $q$ coprime to $e$ we have commutative diagrams
		\begin{eqnarray}  \label{eqn-diagram-thm-tame}
		\begin{gathered}
			\xymatrix{
				\ZZG^{\theta,0}_{\fpb} \times  \Lie(\ZZG)_{M-d}(\fpb) \ar[r] \ar[d]^{\simeq \times \simeq} &  \Lie(\ZZG)_{M-d}(\fpb) \ar[d]^{\simeq}
				&  \ZZG^{\theta,0}_{\fqb} \times \Lie(\ZZG)_{M-d}(\fqb) \ar[r] \ar[d]^{\simeq \times \simeq} &  \Lie(\ZZG)_{M-d}(\fqb) \ar[d]^{\simeq} \\
				\RP_x \times  \PV_{x,r}   \ar[r]    &  \PV_{x,r} 
				& \RP_{x_q} \times  \PV_{x_q,r}   \ar[r]    &  \PV_{x_q,r} \quad .
			} 
		\end{gathered}
		\end{eqnarray}
\end{Thm}

\begin{Rem}
	If $p$ is odd, not torsion for $G$ and does not divide $m$, 
	then, if we choose $M$ to be $m$, the left diagram in \eqref{eqn-diagram-thm-tame} is proven to exist and commute in \cite[Theorem~4.1]{ReederYu}. The proof given in \textit{Loc. cit.} does not work for all primes $p$, because it relies among others crucially on the assumption that $p$ does not divide $m$. 
\end{Rem} 

\textbf{Proof of Theorem \ref{Thm-tame}.}\\
Let $\me, \mf$ be integers such that $e \, | \, \me$, $M=\me \mf$, $\gcd(\me,\mf)=1$ and $\me$ is minimal satisfying these properties. 
	Let $E_{e'}$ be the splitting field of $(x^{e'}-1)$ over $E$, and let $\OEpr$ be the ring of integers in $E_{e'}$. 

We let $\ZZG$ be a split reductive group scheme over $\OEpr[\frac{1}{e}]\subset \ov \bZ[\frac{1}{e}]$ whose root datum $R(\ZZG)$ coincides with the root datum $R(G)$ of $G$, i.e. $\ZZG$ is as defined in Section \ref{Section-global-reductive-quotient} base changed to $\OEpr[\frac{1}{e}]$, and $\ZZT$ denotes a split maximal torus of $\ZZG$. 
Let $G_{ad}$ be the adjoint group of $G$ and $T'$ the subtorus of $T$ that consists of the images of the coroots of $G$. We have the usual map $G \ra G_{ad}$, and we denote the image of $T$ under this map by $T_{ad}$. Restricting the map to $T'$ induces an injection $X_*(T') \hookrightarrow X_*(T_{ad})$ that yields an isomorphism $X_*(T')\otimes_\bZ \bR \xra{\simeq} X_*(T_{ad})\otimes_\bZ \bR$, which we use to identify the two real vector spaces. This allows us to choose  $\lambda \in X_*(T_{ad}) \subset X_*(T')\otimes \bR \subset X_*(T)\otimes \bR$ such that $x=x_0+\frac{1}{M }\lambda$. Note that then, using the identification of $X_*(T)$ with $X_*(T_q)$, we have $x_q=x_{0,q}+\frac{1}{M}\lambda$. We also denote by $\lambda$ the corresponding element in $X_*(\ZZT_{ad}) \subset X_*(\ZZT) \otimes \bR$ under the identification of $X_*(T)$ with $X_*(\ZZT)$. Consider the action $\tl$ of $\mmu_{M}$ on $\ZZG$ given by composition of the closed immersion  $\mmu_{M} \ra \bG_m$ with $\lambda$ and the adjoint action of $\ZZG_{ad}$ on $\ZZG$, i.e.
$$ \theta_\lambda: \mmu_{M} \ra \bG_m \xra{\lambda} \ZZT_{ad} \hookrightarrow \ZZG_{ad} \xra{\Ad} \Aut_{\mbox{\small$\ZZG$}}.$$ 

Let $\out \in \Aut(R(G),\Delta)$ denote the action of $\gammapr \in \Gammapr \simeq \Gal(E/K)$ on $R(G)$ given in the Definition \ref{Def-good} of a \good{} group,
and denote by $\un{\bZ/e\bZ}_{{\Ze}}$ the constant group scheme over $\Spec \OEpr[\frac{1}{e}]$ corresponding to the group $\bZ/e\bZ$.
 Using the Chevalley system $\{ \ZZx_\alpha: \bG_a \ra \ZZU_\alpha \subset \ZZG \}_{\alpha \in \Phi(\ZZG)=\Phi}$ for $(\ZZG,\ZZT)$ (defined in Section \ref{Section-global-reductive-quotient}), the automorphism $\out$ defines a morphism of $\Spec \OEpr[\frac{1}{e}]$-schemes $\un{\bZ/e\bZ}_{{\Ze}} \ra \Aut_{\mbox{\small$\ZZG$}}$. Note that we have an isomorphism of $\Spec {\OEpr[\frac{1}{e}]}$-schemes $\mmu_\me \xra{\simeq} \un{\bZ/\me\bZ}_{{\Ze}}$ that yields the following morphism, which we again denote by $\out$,
\begin{equation*} \out: \mmu_\me \xra{\simeq} \un{\bZ/\me\bZ}_{{\Ze}} \xra{\cdot \frac{\me}{e}} \un{\bZ/e\bZ}_{{\Ze}} \ra \Aut_{\mbox{\small$\ZZG$}}. \end{equation*}
Fix an isomorphism $\mmu_M \simeq \mmu_{e'} \times \mmu_f$. This yields a projection map $p_{M,e'}: \mmu_M \ra \mmu_{e'}$ and allows us to define $\theta: \mmu_M \ra \Aut_{\mbox{\small$\ZZG$}}$ as follows
\begin{equation*} \theta : \mmu_M \xra{\diag} \mmu_M \times \mmu_M  \xra{p_{M,e'} \times \Id} \mmu_\me \times \mmu_M \xra{\out \times \tl} \Aut_{\mbox{\small$\ZZG$}} \times \Aut_{\mbox{\small$\ZZG$}} \xra{\text{mult.}} \Aut_{\mbox{\small$\ZZG$}}. \end{equation*}

By \cite[Propostion~A.8.10]{Pseudoreductive}, the fixed-point locus of $\ZZG$ under the action of $\theta$ is representable by a smooth closed ${\OEpr[\frac{1}{e}]}$-subscheme $\ZZG^{\theta}$ of $\ZZG$.  
Moreover, by \cite[Proposition~A.8.12]{Pseudoreductive}, the fiber $\ZZG^{\theta, 0}_{\ov s}$ is a reductive group for all geometric points $\ov s$ of $\Spec {\OEpr[\frac{1}{e}]}$.
Similarly, $\ZZT^{\theta,0}=\ZZT^{\out,0}$ is a smooth closed subscheme of $\ZZT$.
Hence $\ZZT^{\out, 0}$ is a split torus over $\Spec \OEpr[\frac{1}{e}]$.

Let us denote $\ZZG^{\theta, 0}$ by $\ZzG$. We claim that $\ZZT^{\out,0}$ is a maximal torus of $\ZzG$. In order to prove the claim for geometric fibers, we use a similar argument to one used in \cite[Section~4]{Haines}. 
Let $q$ be an arbitrary prime number coprime to $e$, $\ZzB$ the Borel of $\ZZG$ corresponding to the positive roots, and $\ZzU$ its unipotent radical. As $\ZzG_{\ov \bF_q}$ is a closed subgroup of $\ZZG_{\ov \bF_q}$, $\ZzG_{\ov \bF_q}/(\ZzB_{\ov \bF_q} \cap \ZzG_{\ov \bF_q})$ is proper in $\ZZG_{\ov \bF_q}/\ZzB_{\ov \bF_q}$, hence is proper. Thus $\ZzB_{\ov \bF_q} \cap \ZzG_{\ov \bF_q}$ is a solvable parabolic subgroup, i.e. a Borel subgroup, and $\ZzB_{\ov \bF_q}^{\theta,0}=\ZzB_{\ov \bF_q} \cap \ZzG_{\ov \bF_q}$.  According to \cite[8.2]{Steinberg}, $\ZzU_{\ov \bF_q}^\theta$ is connected, and hence $\ZzB_{\ov \bF_q}^{\theta,0}=\ZZT_{\ov \bF_q}^{\theta,0} \ltimes \ZzU_{\ov \bF_q}^\theta$. This means that $\ZZT_{\ov \bF_q}^{\out,0}=\ZZT_{\ov \bF_q}^{\theta,0} $ is a maximal torus of $\ZzG_{\ov \bF_q}$.  
Hence $\ZZT^{\out,0}_{\ov s}$ is a maximal torus in $\ZzG_{\ov s}$ for all geometric points $\ov s$ of $\Spec {\OEpr[\frac{1}{e}]}$, because the locus of the former points is open. This means that $\ZZT^{\out,0}$ is a maximal torus of $\ZzG$.

 In addition, $\Pic(\Spec {\ov \bZ[\frac{1}{e}]})$ is trivial (by the principal ideal theorem), and hence the root spaces for $(\ZZG_{\ov \bZ[1/e]},\ZZT_{\ov \bZ[1/e]})$ are free line bundles. Using that $\Spec {\ov \bZ[\frac{1}{e}]}$ is connected, we conclude that $\ZzG_{\ov \bZ[1/e]}$ is a split reductive group scheme. 

If $q$ is a large enough prime number, then by \cite[Theorem~4.1]{ReederYu} we have $\ZzG_{\ov \bF_q}\simeq \RP_{x_q}$. 
 Hence $R(\ZG')=R(\ZG'_{\fpb})=R(\RP_{x_q})=R(\ZG)$, and $\ZzG_{\ov \bZ[1/e]}$ is (abstractly) isomorphic to $\ZG$ as desired. 

In order to give a new construction of $\ZV$, let $d:\Aut_{\mbox{\small$\ZZG$}} \ra \Aut_{\mbox{\small$\Lie(\ZZG)$}}$ be the map defined as follows. For any $\OEpr[\frac{1}{e}]$-algebra $R$, and $g \in \Aut_{\mbox{\small$\ZZG$}}(R)$, define $dg:=\Lie(g) \in \Aut(\Lie(\ZZG)_R)$. 
Then the action $d\theta$ defines a $\bZ/M\bZ$-grading on $\Lie(\ZZG)$, which we write as $\Lie \ZZG=\oplus_{i=1}^{M} (\Lie \ZZG)_i$.

We define $\ZzV$ to be the free $\OEpr[\frac{1}{e}]$-module $\Lie(\ZZG)_{M-d}(\OEpr[\frac{1}{e}])$, and the action of $\ZzG:=\ZZG^{\theta,0}$ on $\ZzV$ should be given by the restriction of the adjoint action of $\ZZG$ on $\Lie(\ZZG)(\OEpr[\frac{1}{e}])$. 

In order to show that the $\ZG$-representation on $\ZV$ corresponds to the $\ZzG_{\ov \bZ[1/e]}$-representation on $\ZzV_{\ov \bZ[1/e]}$, we observe that $\ZzV_{\ov \bZ[1/e]}$ is the $M-d$ weight space of the action of $\out \cdot \Ad(\lambda(\zeta_M))$ for some primitive $M$-th root of unity $\zeta_M$ in $\ov \bZ[\frac{1}{e}]$. Using the notation introduced in Section \ref{Section-global-reductive-quotient} preceding Remark \ref{Rem-zeta}, we let $C_\alpha=\zeta_M^{e(\alpha)\cdot \alpha(x-x_0)\cdot M}$. By the same arguments as in the proof of Lemma \ref{Lemma-f}, we see that there exists an  automorphism $h$ of $\ZZG_{\ov \bZ[1/e]}$ that preserves $\ZZT_{\ov \bZ[1/e]}$ and sends $\ZZx_\alpha$ to $\ZZx_\alpha \circ C_\alpha$ for all $\al \in \Phi$.

Let $q$ be a large enough prime, to be more precise: odd, not torsion for $G$ and not dividing $M$. Then we deduce from the arguments used in \cite[Section~4]{ReederYu} that we have commutative diagrams: 
\begin{eqnarray} \label{eqn-splitproof-diragams1}
	\xymatrixcolsep{4pc}
	\xymatrix{
		\ZzG_{\fqb}  \ar@{^{(}->}[r]^{\subset} \ar[d]^{\simeq}_{f_q^{-1}\circ h|_{\ZzG_{\fqb}}} & \ar[d]_{\simeq}^{f_q^{-1}\circ h}  \ZZG_{\fqb}   
		& &		\ZzV_{\fqb}  \ar@{^{(}->}[r]^{\subset} \ar[d]^{\simeq}_{\Lie(f_q^{-1}\circ h)|_{\ZzV({\fqb})}} & \ar[d]_{\simeq}^{\Lie(f_q^{-1}\circ h)}  \Lie(\ZZG)({\fqb}) \\
		\RP_{x_q}  \ar@{^{(}->}[r]_{\iota_{K_q,\Fq}}    & \RPFq_{x_q}  
		& &	\PV_{x_q,r}  \ar@{^{(}->}[r]_{\iota_{K_q,\Fq,r}}    & \PVFq_{x_q,r}   \quad .	
	} 
\end{eqnarray}

Moreover, the diagram on the right hand side is compatible with the action by the groups of the diagram on the left hand side. 

Recall that we constructed in Section \ref{section-global-MP-filtration} a map $\iota:\ZG \ra \ZZG_{\ov \bZ[1/e]}$ and $\ZV$ as a free $\ov \bZ[\frac{1}{e}]$-submodule of $\Lie(\ZZG)(\ov \bZ[\frac{1}{e}])$ (because if $e$ is odd, then $\Phi_K$ does not contain multipliable roots, and hence the submodules $\ZV$ and $\ZVold$ agree in all cases) such that we have the following commutative diagrams for all primes $q$ coprime to $e$ 
\begin{eqnarray} \label{eqn-splitproof-diragams2} 
	\xymatrixcolsep{4pc}
	\xymatrix{
		\ZG_{\fqb}  \ar[r]^{\iota} \ar[d]^{\simeq} & \ar[d]_{\simeq}^{f_q^{-1}}  \ZZG_{\fqb}   
		& &		\ZV_{\fqb}  \ar@{^{(}->}[r]^{\subset} \ar[d]^{\simeq} & \ar[d]_{\simeq}^{\Lie(f_q^{-1})}  \Lie(\ZZG)({\fqb}) \\
		\RP_{x_q}  \ar@{^{(}->}[r]_{\iota_{K_q,\Fq}}    & \RPFq_{x_q}  
		& &	\PV_{x_q,r}  \ar@{^{(}->}[r]_{\iota_{K_q,\Fq,r}}    & \PVFq_{x_q,r}   \quad ,	
	} 
\end{eqnarray}
where the diagram on the right hand side is compatible with the action of the groups on the left hand side by Proposition \ref{Prop-generala}. Note that $\iota_{K_q,F_q}$ is a closed immersion as either $q$ is odd or $e$ is odd (see Section \ref{section-iota}).

Thus we conclude that $h^{-1}(\iota(\ZG_{\fqb}))=\ZzG_{\fqb}$ for large enough $q$.

Let $q$ now be any prime coprime to $e$, and let $\ov g \in \ZG(\fqb)$. As $\ZG(\ov \bZ[\frac{1}{e}])$ surjects onto $\ZG(\fqb)$ (because this holds for the root groups and the torus), we can choose $g \in \ZG(\ov \bZ[\frac{1}{e}])$ whose image in $\ZG(\fqb)$ is $\ov g$. By combining the Diagrams \eqref{eqn-splitproof-diragams1} and \eqref{eqn-splitproof-diragams2}, we see that the image of $h^{-1}\iota(g)$ in $\ZZG(\ov \bF_{q'})$ is actually contained in $\ZzG(\ov \bF_{q'})$ for all sufficiently large primes $q'$. 
Hence $h^{-1}\iota(g) \in \ZzG(\ov \bZ[\frac{1}{e}]) \subset \ZZG(\ov \bZ[\frac{1}{e}])$, and $h^{-1}\circ \iota(\ZG(\fqb)) \subset \ZzG(\fqb)$.  Since we observed that $\ZzG_{\fqb}$ is abstractly isomorphic to $\RP_{x_q} \simeq h^{-1} \circ f_q(\iota_{K_q,\Fq}(\RP_{x_q})) \simeq h^{-1} \circ \iota(\ZG_{\fqb})$, we conclude that 
\begin{equation} \label{eqn-fiberq}
	h^{-1} \circ \iota(\ZG_{\fqb}) = \ZzG_{\ov \bF_q}
\end{equation}
for all primes $q$ coprime to $e$. The same arguments show that 
\begin{equation} \label{eqn-fiberp}
h^{-1} \circ \iota(\ZG_{\fpb}) = \ZzG_{\ov \bF_p}.
\end{equation}

Moreover, we claim that $h^{-1} \circ \iota(\ZG_{\ov \bQ}) = \ZzG_{\ov \bQ}$. In order to prove the claim, note that $(\mmu_M)_{\ov \bQ} \simeq \un{\bZ/M\bZ}_{\ov \bQ}$, and hence the action of the group scheme $\mmu_M$ on $\ZZG_{\ov \bQ}$ corresponds to the action of the finite group $\bZ/M\bZ$ generated by $\vartheta \cdot \Conjug(\lambda(\zeta_M))$. Therefore, by the construction of $\iota:\ZG_{\ov \bZ[1/e]} \ra \ZZG_{\ov \bZ[1/e]}$ (see Proposition \ref{Prop-iota}) and the definition of $h:\ZZG_{\ov \bZ[1/e]} \ra \ZZG_{\ov \bZ[1/e]}$, we see that $h^{-1} \circ \iota(\ZG(\ov \bQ)) \subset \ZZG^{\theta}({\ov \bQ})$. As $\iota_{\ov \bQ}:\ZG_{\ov \bQ} \ra \ZZG_{\ov \bQ}$ is a closed immersion by Lemma \ref{Lemma-closed-immersion}, $h^{-1} \circ \iota(\ZG_{\ov \bQ})  \simeq \ZG_{\ov \bQ} \simeq \ZZG_{\ov \bQ}^{\theta,0}=\ZzG_{\ov \bQ}$, and we conclude that 
\begin{equation} \label{eqn-fiberQ}
h^{-1} \circ \iota(\ZG_{\ov \bQ}) = \ZzG_{\ov \bQ}.
\end{equation}
Thus, as $\ZzG_{\ov \bZ[1/e]}$ is smooth over $\Spec \ov \bZ[\frac{1}{e}]$, hence reduced, we deduce from the Nullstellensatz that $h^{-1} \circ \iota: \ZG \ra \ZZG_{\ov \bZ[1/e]}$ factors via the closed subscheme $\ZzG_{\ov \bZ[1/e]}$ of $\ZZG_{\ov \bZ[1/e]}$, i.e. we may write $h^{-1} \circ \iota: \ZG \ra \ZzG_{\ov \bZ[1/e]}$. 
As we proved that $(h^{-1} \circ \iota)_s:\ZG_s \ra (\ZzG_{\ov \bZ[1/e]})_s$ is an isomorphism for all $s \in \Spec \ov \bZ[\frac{1}{e}]$ (see Equation \eqref{eqn-fiberq}, \eqref{eqn-fiberp}, \eqref{eqn-fiberQ}), we conclude  that by \cite[17.9.5]{EGAIV4} the morphism $h^{-1} \circ \iota: \ZG \ra \ZzG_{\ov \bZ[1/e]}$ is an isomorphism. 

Moreover, as $\Lie(h)(\ZzV_{\fqb})=\ZV_{\fqb}$ for large enough primes $q$, we deduce that $\Lie(h):\Lie(\ZZG)(\ov \bZ[\frac{1}{e}])$ $ \ra \Lie(\ZZG)(\ov \bZ[\frac{1}{e}])$ yields an isomorphism of the direct $\ov \bZ[\frac{1}{e}]$-module summands $\ZzV_{\ov \bZ[1/e]}$ and $\ZV$.

As the action of $\ZG$ on $\ZV$ was defined via the adjoint action of $\ZZG_{\ov \bZ[1/e]} \supset \iota(\ZG)$ onto $\Lie(\ZZG_{\ov \bZ[1/e]})(\ov \bZ[\frac{1}{e}]) \supset \ZV$, the isomorphisms 
$$h^{-1}: \ZG \ra \ZzG_{\ov \bZ[1/e]}=\ZZG_{\ov \bZ[1/e]}^{\theta, 0}
\quad \text{ and } \quad 
 \Lie(h^{-1}): \ZV \ra \ZzV_{\ov \bZ[1/e]}=\Lie(\ZZG_{\ov \bZ[1/e]})_{M-d}\left(\ov \bZ[\frac{1}{e}]\right)$$ map the action of $\ZG$ onto $\ZV$ to the action of $(\ZZG_{\ov \bZ[1/e]})^{\theta, 0}$ on $\Lie(\ZZG_{\ov \bZ[1/e]})_{M-d}(\ov \bZ[\frac{1}{e}])$ which arises from the restriction of the adjoint action of $\ZZG_{\ov \bZ[1/e]}$ on $\Lie(\ZZG_{\ov \bZ[1/e]})(\ov \bZ\left[\frac{1}{e}\right])$. 

The commutative diagrams in the theorem now follow by applying Theorem \ref{Thm-general}.
\qed

\begin{Rem}
	Let $E_{e'}$ be as defined in the proof of Theorem \ref{Thm-tame}. 
	 Denote by $E_H$ the Hilbert class field of $E_{e'}$ and by $\OH$ the ring of integers in $E_H$. Then the group schemes $\ZG$ and $\ZZG$ and the action of $\ZG$ on $\ZV$ appearing in Theorem \ref{Thm-tame} can be defined over $\Spec \OH[\frac{1}{e}]$.

\end{Rem}

\subsection{Vinberg--Levy theory for all \good{} groups} \label{section-VinbergQ}
Even though the Moy--Prasad filtration representation of groups that do not split over a tamely ramified extension might not be described as in Vinberg--Levy theory, its lift to characteristic zero can be described using Vinberg theory, i.e. as the fixed-point subgroup of a finite order automorphism on a larger group acting on some eigenspace in the Lie algebra of the larger group. To be more precise, we have the following corollary of Theorem \ref{Thm-tame} combined with Theorem \ref{Thm-general}.

\begin{Cor} \label{Cor-Vinberg}
	Let $G$ be a \good{} reductive group over $K$, $r=\frac{d}{M}$ for some nonnegative integer $d<M$, and let the representation of $\ZG$ on $\ZV$ be as in Theorem \ref{Thm-general}. Then there exist a reductive group scheme $\ZZG_{\ov \bQ}$ over $\ov \bQ$ and morphisms 
	$$ \theta: \mmu_M \ra \Aut_{{\mbox{\small$\ZZG_{\ov \bQ}$}}/ \ov \bQ}  \quad  \text{ and } \quad d\theta:  \mmu_{M} \ra \Aut_{\mbox{\small$\Lie(\ZZG_{\ov \bQ})$}/\ov \bQ}$$
	such that $\ZG_{\ov \bQ} \simeq \ZZG_{\ov \bQ}^{\theta,0}$ and $\ZV_{\ov \bQ} \simeq \Lie(\ZZG_{\ov \bQ})_{M-d}(\ov \bQ)$, and the action of $\ZG_{\ov \bQ}$ on $\ZV_{\ov \bQ}$ corresponds via these isomorphisms to the restriction of the adjoint action of $\ZZG_{\ov \bQ}$ on $\Lie(\ZZG_{\ov \bQ})(\ov \bQ)$.
\end{Cor}
\textbf{Proof.}\\
Let $q$ be a prime larger than $p^sN$. Then, by construction, the representation over $\ov \bZ[\frac{1}{p^sN}]$ associated to $G_q$ via the proof of Theorem \ref{Thm-general} agrees with the representation of $\ZG_{\ov \bZ[1/(p^sN)]}$ on $\ZV_{\ov \bZ[1/(p^sN)]}$. As $G_q$ splits over a tamely ramified extension, Theorem \ref{Thm-tame} allows us to deduce the corollary. \qed

\section{Semistable and stable vectors} \label{section-semistablestablevectors}
In this section we apply our results of Section \ref{section-globalMP} and Section \ref{section-Vinberg} to prove that the existence of stable and semistable vectors in the Moy--Prasad filtration representations is independent of the characteristic of the residue field. Recall that a vector $v$ in a vector space $V$ over an algebraically closed field is \textit{stable} \index{definition}{stable} under the action of a reductive group $G_V$ on $V$ if the orbit $G_Vv$ is closed and the stabilizer $\Stab_{G_V}(v)$ of $v$ in $G_V$ is finite. A vector $v \in V$ is called \textit{semistable} \index{definition}{semistable} if the closure of the orbit $G_Vv$ does not contain zero.

\subsection{Semistable vectors} \label{section-semistablevectors}

The global version of the Moy--Prasad filtration representation as provided by Theorem \ref{Thm-general} allows us to show that the existence of semistable vectors is independent of the residual characteristic $p$ of $K$ as follows, where $N$ is the integer coprime to $p$ introduced in Definition \ref{Def-good}.
\begin{Thm} \label{Thm-semistable}
	We keep the notation used in Theorem \ref{Thm-general}, in particular $G$ is a \good{} reductive group over $K$ and $x \in \sB(G,K)$. Then the following are equivalent
	\begin{enumerate}[(i)]
	\item \label{semistable-i} $\PV_{x,r}$ has semistable vectors under the action of $\RP_{x}$. 
	\item \label{semistable-ii} $\PV_{x_q,r}$ has semistable vectors under the action of $\RP_{x_q}$ for some prime $q$ coprime to $N$.
	\item \label{semistable-iii} $\PV_{x_q,r}$ has semistable vectors under the action of $\RP_{x_q}$ for all primes $q$ coprime to $N$.
	\end{enumerate}
\end{Thm}

\textbf{Proof.}\\
We first show that \ref{semistable-ii} implies \ref{semistable-i}. Suppose that \ref{semistable-ii} holds, i.e. that $\PV_{x_q,r}$ contains semistable vectors under $\RP_{x_q}$ for some prime $q$ coprime to $N$. This implies by \cite[Proposition~4.3]{MP1} that $\ZV_{\ov \bQ_q}$ has semistable vectors under the action of $\ZG_{\ov \bQ_q}$, where $\ZG$ and $\ZV$ are as in Theorem \ref{Thm-general}.  By \cite[p.~41]{Mumford} (based on \cite[Definition~1.7 and Proposition~2.2]{Mumfordbook}) this means that there exists a $\ZG_{\ov \bQ_q}$-invariant non-constant homogeneous element $P_q$ in $\Sym \check \ZV_{\ov \bQ_q}$. Moreover, there exists $X \in \ZV_{\ov \bQ} \subset \ZV_{\ov \bQ_q}$ such that $P_q(X) \neq 0$, i.e. $X$ is semistable in $\ZV_{\ov \bQ_q}$ under the action of $\ZG_{\ov \bQ_q}$.
 Hence $X \neq 0$ is also semistable in $\ZV_{\ov \bQ}$ under the action of $\ZG_{\ov \bQ}$, which implies $(\Sym \check \ZV_{\ov \bQ})^{\ZG(\ov\bQ)} \neq \ov \bQ$. Thus, there does also exist a $\ZG(\ov \bZ)$-invariant non-constant homogeneous element $P$ in $\Sym \check \ZV_{\ov \bZ}$.  
  As $P$ is non-constant and homogeneous, we can assume without loss of generality that the image $\ov P$ of $P$ in $\Sym \check  \ZV_{\ov \bZ} \otimes {\fpb} \simeq \Sym \check \ZV_{\ov \bF_p}$ is non-constant. Note that $\ZG(\ov \bZ)$ surjects onto $\ZG({\fpb})$, which follows from the surjections 
  on all root groups and the split maximal torus.
Hence $\ov P$ is $\ZG({\fpb})\simeq \RP_x({\fpb})$-invariant and there exists $\ov X \in \ZV_{{\fpb}} \simeq \PV_{x,r}$ such that $\ov f(\ov X) \neq 0$, i.e. $\ov X$ is semistable by \cite[p.~41]{Mumford}. Thus \ref{semistable-i} is true.

The same arguments show that if $\RP_{x,r}$ has semistable vectors, then $\RP_{x_q,r}$ has semistable vectors for all primes $q$ coprime to $N$, i.e. \ref{semistable-i} implies \ref{semistable-iii}. As \ref{semistable-iii} implies \ref{semistable-ii}, we conclude that all three statements are equivalent.
 \qed

Note that the same holds for the linear duals $\check \PV_{x,r}$ and $\check \PV_{x_q,r}$ of $\PV_{x,r}$ and $\PV_{x_q,r}$ using $\check \ZV$ instead of $\ZV$ in the proof above: 
\begin{Cor} \label{Cor-semistable}
		We use the same notation as above. Then $\check \PV_{x,r}$ has semistable vectors under the action of $\RP_{x}$ if and only if $\check \PV_{x_q,r}$ has semistable vectors under the action of $\RP_{x_q}$ for some prime $q$ coprime to $N$ if and only if $\check \PV_{x_q,r}$ has semistable vectors under the action of $\RP_{x_q}$ for all primes $q$ coprime to $N$.
\end{Cor}

\begin{Rem}
	For semisimple groups $G$ that split over a tamely ramified extension and sufficiently large residue field characteristic $p$, Reeder and Yu classified  in \cite[Theorem~8.3]{ReederYu} those $x$ for which $\check \PV_{x,r}$ contains semistable vectors in terms of conditions that are independent of the prime $p$. Corollary \ref{Cor-semistable} allows us to conclude that these prime independent conditions also classify points $x$ such that $\PV_{x,r}$ contains semistable vectors	for all \good{} semisimple groups $G$ (without any restriction on the residue field characteristic). Note that the removal of the restriction on the residue field characteristic for absolutely simple split reductive groups $G$ is also contained in a joint paper of the author with Romano, see \cite{Splitcase}. For this result, it suffices to construct $\ZG$ acting on $\ZV$ over $\ov \bZ_p$.
\end{Rem}

\subsection{Stable vectors}
In this section we show an analogous result to the one in Section \ref{section-semistablevectors} for stable vectors. This allows us to generalize the criterion in \cite{ReederYu} for the existence of stable vectors in the dual of the first Moy--Prasad filtration quotient to arbitrary residual characteristics $p$ and all \good{} semisimple groups, which in turn produces new supercuspidal representations.

\begin{Thm} 
	We keep the notation used above, in particular $G$ is a \good{} reductive group over $K$ and $x \in \sB(G,K)$. 
	Then the following are equivalent
	\begin{enumerate}[(i)]
		\item \label{stable-i} $\PV_{x,r}$ has stable vectors under the action of $\RP_{x}$. 
		\item \label{stable-ii} $\PV_{x_q,r}$ has stable vectors under the action of $\RP_{x_q}$ for some prime $q$ coprime to $N$.
		\item \label{stable-iii} $\PV_{x_q,r}$ has stable vectors under the action of $\RP_{x_q}$ for all primes $q$ coprime to $N$.
	\end{enumerate}
\end{Thm}

Before we prove the theorem, we would like to mention that part of the following proof appears as well in \cite{Splitcase} in order to proof Corollary \ref{Cor-stablevectors} below in the case of $G$ being absolutely simple and split.

\textbf{Proof.} 
We suppose without loss of generality that $r=\frac{d}{M}$ for some nonnegative integer $d<M$. 

Assume that \ref{stable-ii} is satisfied, i.e. there exists a prime $q$ coprime to $N$ such that $\PV_{x_q,r}$ contains stable vectors under the action of $\RP_{x_q}$.

A slight variation of the proof by Moy and Prasad of \cite[Proposition~4.3]{MP1} (see \cite[Lemma~2]{Splitcase} for a detailed proof) shows that then $\ZV_{\ov \bQ_q}$ contains stable vectors under $\ZG_{\ov \bQ_q}$, where $\ZG$ and $\ZV$ are as in Theorem \ref{Thm-general}.

Recall that by Corollary \ref{Cor-Vinberg}, $\ZG_{\ov \bQ} \simeq \ZZG_{\ov \bQ}^{\theta,0}$ and $\ZV_{\ov \bQ} \simeq \Lie(\ZZG_{\ov \bQ})_{M-d}(\ov \bQ)$ such that the action of $\ZG_{\ov \bQ}$ on $\ZV_{\ov \bQ}$ corresponds via these isomorphisms to the restriction of the adjoint action of $\ZZG_{\ov \bQ}$ on $\Lie(\ZZG_{\ov \bQ})(\ov \bQ)$. 
Let $\zeta_M$ be a primitive $M$-th root of unity in $\ov \bQ$,  denote $\ZZG_{\ov \bQ}^{\theta({\zeta_M})^{M/(d,M)}, 0}$ by $\ZZG'$, its Weyl group by $W'$, and let $\vartheta$ be the action of $\theta(\zeta_M)$ on the root datum $R(\ZZG'_{\ov \bQ_q})$.
 Then by \cite[Corollary~14]{RLYG}, the existence of stable vectors in $\ZV_{\ov \bQ_q}$ is equivalent to the action of $\theta(\zeta_M)$ on $\ZZG_{\ov \bQ_q}'$ (or, equivalently, on $\ZZG'$) being principal and $\frac{M}{(d,M)}$ being the order of an elliptic $\bZ$-regular element of $W'\vartheta$. Hence we conclude by the same equivalence for the prime $p$ that there exist stable vectors in $\ZV_{\ov \bQ_p}$ under the action of $\ZG_{\ov \bQ_p}$.

Thus the set of stable vectors $( \ZV_{\ov \bQ_p})_s$ in $ \ZV_{\ov \bQ_p}$ is non-empty and open (see \cite[1.4, p.~37]{Mumfordbook}). Hence there exists a nonzero polynomial $P$ in the space global sections $\cO_{\ZV}( \ZV_{\ov \bQ_p})=\cO_{\ZV}( \ZV_{\ov \bZ_p})  \otimes_{\ov \bZ_p} \ov \bQ_p \simeq \ov \bZ_p[x_1, \hdots, x_n] \otimes_{\ov \bZ_p} \ov \bQ_p = \ov \bQ_p[x_1, \hdots, x_n] $ such that the $\ov \bQ_p$-points of the closed reduced subvariety $V(P)$ of $\ZV_{\ov \bQ_p}$ defined by the vanishing locus of $P$  contain $\left({ \ZV_{\ov \bQ_p}}-( \ZV_{\ov \bQ_p})_s\right) \ni 0$. We can assume without loss of generality that the coefficients of $P$ are in $\ov \bZ_p$, i.e. $P \in \cO_{\ZV}( \ZV_{\ov \bZ_p}) \subset \cO_{\ZV}(\ZV_{\ov \bQ_p})$, and that at least one coefficient of $P$ has p-adic valuation zero.  Let $\ov P$ be the image of $P$ under the reduction map $\cO_{\ZV}({ \ZV_{\ov \bZ_p}}) \simeq \ov \bZ_p[x_1, \hdots, x_n] \ra \cO_{\ZV}({ \ZV_{{\fpb}}}) \simeq \ov \bF_p[x_1, \hdots, x_n] $. Then $\ov P$ is not constant, because $P(0)=0$, and there exists $\ov X \in  \ZV_{{\fpb}} \simeq  \PV_{x,r}$ such that $\ov P(\ov X) \neq 0$.

 We claim that $\ov X$ is a stable vector under the action of $\RP_x$. We will prove the claim using the Hilbert--Mumford Criterion that states that a vector is stable if and only if it has positive and negative weights for every non-trivial one-parameter subgroup, see \cite[p.~41]{Mumford} based on \cite[Theorem~2.1]{Mumfordbook}. Let $\ov \lambda: \bG_m \ra \RP_x \simeq \ZG_{{\fpb}}$ be a non-trivial one parameter subgroup. Then $\ov \lambda$ is defined over some finite extension of $\bF_p$, and hence by \cite[IX, Corollaire~7.3]{SGA3II} there exists a lift $\lambda: \bG_m \ra \ZG_{\ov \bZ_p}$ of $\ov \lambda$. The composition of $\lambda$ with the action of $\ZG_{\ov \bZ_p}$ on $ \ZV_{\ov \bZ_p}$ yields an action of $\bG_m$ on $ \ZV_{\ov \bZ_p}$, and we obtain a weight decomposition  $ \ZV_{\ov \bZ_p} = \oplus_{m \in \mathbb Z}  \ZV_m$. Denote $\oplus_{m \in \bZ_{>0}}  \ZV_m$ by $ \ZV_+$ and $\oplus_{m \in \bZ_{<0}}  \ZV_m$ by $ \ZV_-$, i.e. $ \ZV_{\ov \bZ_p}= \ZV_- \oplus  \ZV_0 \oplus   \ZV_+$. Let $X \in  \ZV_{\ov \bZ_p}$ be a lift of $\ov X$, and write $X=X_- + X_0 + X_+$ with $X_- \in  \ZV_-, X_0 \in  \ZV_0, X_+ \in  \ZV_+$. Note that the weight decomposition of $ \ZV_{{\fpb}}$ under the action of $\bG_m$ via the composition of $\ov \lambda$ with the action of $\ZG_{{\fpb}}$ on $ \ZV_{{\fpb}}$ is the image of the decomposition $ \ZV_- \oplus  \ZV_0 \oplus  \ZV_+$, i.e. $(  \ZV_{{\fpb}})_-=\oplus_{m\in \bZ_{<0}} ( \ZV_{{\fpb}})_m=(\ZV_-)_{{\fpb}}, ( \ZV_{{\fpb}})_0=( \ZV_0)_{{\fpb}} $ and $( \ZV_{{\fpb}})_+=\oplus_{m\in \bZ_{>0}}( \ZV_{{\fpb}})_m=( \ZV_+)_{{\fpb}}$. Hence $\ov X=\ov X_- + \ov X_0 + \ov X_+$ (where an overline denotes the image after base change to ${\fpb}$) has positive and negative weights with respect to $\ov \lambda$ if and only if $\va(X_-)=0=\va(X_+)$.

Suppose that $\va(X_-)>0$. Then $P(X) \equiv P(X_0+X_+)$ modulo the maximal ideal of $\ov \bZ_p$. However, $X_0+X_+$ is not a stable vector, because it has no negative weights with respect to the non-trivial one parameter subgroup $\lambda \times_{\ov \bZ_p} \ov \bQ_p$, which implies $P(X_0+X_+)=0$. Hence $\ov P( \ov X)=0$ contradicting the choice of $\ov X$. The same contradiction arises if we assume that $\va(X_-)>0$. Thus, $\ov X$ has positive and negative weights for every nontrivial one parameter subgroup, i.e. $\ov X$ is stable by the Hilbert--Mumford criterion. Hence, statement \ref{stable-i} of the theorem holds.

The same arguments show that if $\RP_{x,r}$ has stable vectors, then $\RP_{x_q,r}$ has stable vectors for all $q$ coprime to $N$, i.e. \ref{stable-i} implies \ref{stable-iii}. As \ref{stable-iii} implies \ref{stable-ii}, the three statements are equivalent.
 \qed
 
 As in the semistable case, the same proof works for the linear duals of the Moy--Prasad filtration quotients:

\begin{Cor} \label{Cor-stable}
	We use the same notation as above. Then $\check \PV_{x,r}$ has stable vectors under the action of $\RP_{x}$ if and only if $\check \PV_{x_q,r}$ has stable vectors under the action of $\RP_{x_q}$ for some prime $q$ coprime to $N$ if and only if $\check \PV_{x_q,r}$ has stable vectors under the action of $\RP_{x_q}$ for all primes $q$ coprime to $N$.
\end{Cor}
Denote by  $r(x)$ the smallest positive real number such that $\PV_{x,r(x)} \neq \{ 0\}$, and let $\check \rho=\frac{1}{2}\sum\limits_{\al \in \Phi^+}\check \alpha$, where $\Phi^+$ are the positive roots of $\Phi=\Phi(G)$ (with respect to the fixed Borel $B$). Then Corollary \ref{Cor-stable} allows us to classify the existence of stable vectors in $\check \PV_{x,r(x)}$ for arbitrary primes $p$ and \good{} semisimple groups below. This generalizes the result of \cite[Corollary~5.1]{ReederYu} for large primes $p$ and semisimple groups that split over tamely ramified extensions.

\begin{Cor} \label{Cor-stablevectors} Let $G$ be a \good{} semisimple group and $x$ a rational point of order $m$ in $\sA(S,K) \subset \sB(G,K)$. Then $\check \PV_{x,r(x)}$ contains stable vectors under $\RP_{x}$ if and only if $x$ is conjugate under the affine Weyl group $\Waff$ of the restricted root system of $G$ to $x_0+\check \rho/m$, $r(x)=1/m$ and there exists an elliptic $\bZ$-regular element $w\gammapr$ of order $m$ in $W\gammapr$, where $W$ is the absolute Weil group of $G$ and $\gammapr$ is the  automorphism of $R(G)$ given in the definition of a \good{} group (Definition \ref{Def-good}). 
\end{Cor}
\textbf{Proof.} 
Note that by Lemma \ref{lemma-affine-roots}, the order of $x_q$ is $m$, and by Theorem \ref{Thm-general}, we have $r(x_q)=r(x)$. Let $q$ be sufficiently large, i.e. coprime to $M$, not torsion and odd. Then $G_q$ is a semisimple group that splits over a tamely ramified extension, and we deduce from the proof of \cite[Lemma~3.1]{ReederYu} that $\check \PV_{x_q,r(x_q)}$ can only admit stable vectors under $\RP_{x_q}$ if $x_q$ is a barycenter of some facet of $\sA_q=\sA(S_q,K_q)$, and hence $r(x_q)=1/m$. 
 Therefore, as $q$ is chosen sufficiently large, we obtain by \cite[Corollary~5.1]{ReederYu} that $\check \PV_{x_q,r(x_q)}$ has stable vectors if and only if $x_q$ is conjugate under the affine Weyl  group ${\Waff}_{q}$ of the restricted root system of $G_q$ to $x_{0,q}+\check \rho/m$, $r(x)=1/m$ and there exists an elliptic $\bZ$-regular element $w\gammapr$ of order $m$ in $W\gammapr$, because $W$ is isomorphic to the absolute Weil group of $G_q$. Note that 
 $$x_q \sim_{{\Waff}_q} x_{0,q}+\check \rho/m \quad \text{ if and only if } \quad x \sim_{\Waff} x_0+\frac{1}{4} \sum_{a \in \Phi_K^{+,\text{mul}}} \va(\lambda_a) \cdot \check a+\frac{\check \rho}{m}, $$
  and $x_0+\frac{1}{4} \sum\limits_{a \in \Phi_K^{+,\text{mul}}} \va(\lambda_a) \cdot \check a+\check \rho/m$ is conjugate to $x_0+\check \rho/m$ under the extended affine Weyl group of the restricted root system of $G$. However, by checking the tables for all possible points $x_q$ whose first Moy--Prasad filtration quotient $\check \PV_{x_q,r(x_q)}$ admits stable vectors in \cite{RLYG} and \cite{ReederYu}, we observe that the latter conjugacy can be replaced by conjugacy under the (unextended) affine Weyl group. Hence using Corollary \ref{Cor-stable}, we conclude that $\check \PV_{x,r(x)}$ contains stable vectors under the action of $\RP_x$ if and only if $x \sim_{\Waff} x_0 + \check \rho /m$, $r(x)=1/m$, and there exists an elliptic $\bZ$-regular element of order $m$ in $W\gammapr$. \qed

Recall that $k$ is a nonarchimedean local field with maximal unramified extension $K$.
\begin{Cor}
	Let $G$ be a good semisimple group, and suppose that $G$ is defined over $k$. Assume that $W\gammapr$ contains an elliptic $\bZ$-regular element. Then using the construction of \cite[Section~2.5]{ReederYu} we obtain supercuspidal (epipelagic) representations of $G(k')$ for some finite unramified field extension $k'$ of $k$.
	\end{Cor}
\textbf{Proof.} Let $m$ be the order of an elliptic $\bZ$-regular element of $W\gammapr$, and $x=x_0+\check\rho/m \in \sA(S,K)$. By Corollary \ref{Cor-stablevectors}, $\check \PV_{x,r(x)} $ contains stable vectors under the action of $\RP_x$. Since $x$ is fixed under the action of the Galois group $\Gal(K/k'')$ for some finite unramified extension $k''$ of $k$, the vector space $\check \PV_{x,r(x)}$ is defined over the residue field $\ff''$ of $k''$. Hence there exists a finite unramified field extension $k'$ of $k$ with residue field $\ff'$ such that $\check \PV_{x,r(x)}$ contains a stable vector defined over $\ff'$. Applying \cite[Proposition~2.4]{ReederYu} yields the desired result. 	 \qed

\section{Moy--Prasad filtration representations as Weyl modules} \label{section-reps}
In this section we describe the Moy--Prasad filtration representations in terms of Weyl modules. Recall that for $\lambda \in X^*(\ZS)$ a dominant weight, the \textit{Weyl module} \index{definition}{Weyl module}  $V(\lambda)$ \index{notation}{$V(\lambda)$} (over $\ov \bZ[1/N]$) is given by
\begin{equation*}
	V(\lambda)=\ind_{\ZB^-}^\ZG(-w_0\lambda)^\vee,
\end{equation*}
where $\ZB$ is the Borel subgroup of $\ZG$ corresponding to $\Delta(\ZG)$, $\ZB^-$ is the opposite Borel subgroup corresponding to $-\Delta(\ZG)$, $w_0$ is the longest element of the Weyl group of $\Phi(\ZG)$, and $(.)^\vee$ denotes the dual (\cite[II.8.9]{Jantzen}). We define \index{notation}{$\Phi_{x,r}$} \index{notation}{$\Phixr$}
\begin{eqnarray*}
	\Phi_{x,r} & = &  \left\{ a \in \Phi_K \, | \, r-a(x-x_0) \in \Gamma_a'(G) \right\} \\
	\Phixr & = & \left\{ a \in \Phi_{x,r} \, | \, a+b \not \in \Phi_{x,r} \text{ for all } b \in \Phi^+(\ZG)\subset\Phi_K \right\}.
\end{eqnarray*}

\subsection{The split case}
If $G$ is split over $K$, then 
\begin{equation*}
	\Phixr = \left\{ \alpha \in \Phi \, | \, r-\alpha(x-x_0) \in \bZ, \alpha+\beta \not \in \Phi \text{ for all } \beta \in \Phi^+(\ZG)\subset\Phi \right\}.
\end{equation*}

\begin{Thm} \label{Thm-Weyl-split}
	Let $G$ be a split reductive group over $K$, $r$ a real number and $x$ a rational point of $\sB(G,K)$. Let $\ZV$ be the corresponding global Moy--Prasad filtration representation of the split reductive group scheme $\ZG$ over $\ov \bZ$ (Theorem \ref{Thm-general}). Then 
	\begin{equation*}
		\ZV \simeq \left\{ \begin{array}{ll} \Lie(\ZG)(\ov \bZ) & \text{ if $r$ is an integer } \\
							\bigoplus_{\lambda \in \Phixr } V(\lambda) & \text{ otherwise } \, .
				\end{array} \right. 
	\end{equation*}
\end{Thm}	
\textbf{Proof.}
If $r$ is an integer, then we have by Theorem \ref{Thm-tame} that $\ZV\simeq \Lie(\ZZG)_{M}(\ov \bZ)=\Lie(\ZZG^\theta)(\ov \bZ)=\Lie(\ZG)(\ov \bZ)$.

Suppose $r$ is not an integer. Then $\ZV \subset \Lie(\ZZG)(\ov \bZ)$ is spanned by $\ZZX_\alpha=\Lie(\ZZx_\alpha)(1)$ for $\alpha \in \Phi_{x,r}$ (Section \ref{Section-global-filtration-quotient}). Thus the weights in $\Phixr$ are the highest weights of the representation of $\ZG$ on $\ZV$, and we have $\ZV_{\ov \bQ} \simeq \bigoplus_{\lambda \in \Phi_{x,r}^{\rm{max}} } V(\lambda)_{\ov \bQ}$. In order to show that  $\ZV \simeq \bigoplus_{\lambda \in \Phi_{x,r}^{\rm{max}} } V(\lambda)$, it suffices by \cite[II.8.3]{Jantzen} to prove that $\{\ZG(\ov \bZ)(\ZZX_\alpha)\}_{\alpha \in \Phixr}$ spans $\ZV$, i.e. that $\<\ZG(\ov \bZ)(\ZZX_\alpha)\>_{\alpha \in \Phixr}$ contains $\ZZX_\alpha$ for all $\alpha \in \Phi_{x,r}$. Let $\alpha \in \Phi_{x,r}\backslash \Phixr$. Then there exists $\beta \in \Phi^+(\ZG)$ such that $\alpha+\beta \in \Phi$. Let $N_{\alpha,\beta}>0$ be the maximal integer such that $\alpha+N_{\alpha,\beta} \beta \in \Phi$, and let $N^-_{\alpha,\beta}$ be the maximal integer such that $\alpha-N^-_{\alpha,\beta} \beta \in \Phi$. We claim that $\ZZX_{\alpha}+N_{\alpha,\beta} \beta \in \<\ZG(\ov \bZ)(\ZZX_\alpha)\>_{\alpha \in \Phixr}$ implies that $\ZZX_\alpha \in \<\ZG(\ov \bZ)(\ZZX_\alpha)\>_{\alpha \in \Phixr}$, which will imply the theorem by induction.

Suppose that $\ZZX_{\alpha}+N_{\alpha,\beta} \beta \in \<\ZG(\ov \bZ)(\ZZX_\alpha)\>_{\alpha \in \Phixr}$.
Note that $N_{\alpha,\beta}+N^-_{\alpha,\beta} \in \{1,2,3\}$, and recall that 
\begin{equation} \label{eqn-action}
	\ZZx_{-\beta}(u)(\ZZX_{\alpha+N_{\alpha,\beta} \beta})=\sum_{i=0}^{N_{\alpha,\beta}+N^-_{\alpha,\beta}} m_{\alpha,\beta,i} u^i \ZZX_{\alpha+(N_{\alpha,\beta}-i) \beta} \,
\text{ with } m_{\alpha,\beta,i}\in \left\{ \pm 1 \right\} , 
\end{equation}
for  $u\in \bG_a(\ov \bZ)$. By varying $u\in \bG_a(\ov \bZ)$ and taking linear combinations, we conclude that $\ZZX_{\alpha}$ is in the $\ov \bZ$-span of $\{\ZG(\ov \bZ)(\ZZX_\alpha)\}_{\alpha \in \Phixr}$.
 \qed

The following corollary follows immediately by combining Theorem \ref{Thm-Weyl-split} and Theorem \ref{Thm-general}.
\begin{Cor}
	Let $G$ be a split reductive group over $K$, $r$ a real number and $x$ a rational point of $\sB(G,K)$. Then the representation of $\RP_x$ on $\PV_{x,r}$ is given by
	\begin{equation*}
	\PV_{x,r} \simeq \left\{ \begin{array}{ll} \Lie(\RP_x)(\ov \bF_p) & \text{ if $r$ is an integer } \\
	\bigoplus_{\lambda \in \Phi_{x,r}^{\rm{max}} } V(\lambda)_{\fpb} & \text{ otherwise } \, .
	\end{array} \right. 
	\end{equation*}  
\end{Cor}	

\begin{Rem} 
	Note that, if $p$ is sufficiently large, then $V(\lambda)_{\ov \bF_p}$ is an irreducible representation of $\RP_x$ of highest weight $\lambda$. 
\end{Rem}

\subsection{The general case}
Let $a \in \Phixr$ and let $\sU_H$ be the unipotent radical of $\sB_H$. By Frobenius reciprocity, we have (\cite[Proof of Lemma~II.2.13a)]{Jantzen})
\begin{eqnarray*} 
	\Hom_{\ZG}\left(V(a),\Lie(\ZZG)(\ov \bZ[1/N])\right) & \simeq & \Hom_{\ZG}\left(\Lie(\ZZG)(\ov \bZ[1/N])^\vee, \ind_{\ZB^-}^\ZG(-w_0a)\right) \\
	& \simeq &  \Hom_{\ZB^-}\left(\Lie(\ZZG)(\ov \bZ[1/N])^\vee, -w_0a\right) \\
	& \simeq &  \Hom_{\ZB^-}\left(w_0a, \Lie(\ZZG)(\ov \bZ[1/N])\right) \simeq  \left(\left(\Lie(\ZZG)(\ov \bZ[1/N]\right)^{\sU_H}\right)_{a}.
\end{eqnarray*}

Using these isomorpisms, the element $Y_a \in \left(\left(\Lie(\ZZG)(\ov \bZ[1/N]\right)^{\sU_H}\right)_{a} \subset \Lie(\ZZG)(\ov \bZ[1/N])$ yields a morphism $V(a) \ra  \Lie(\ZZG)(\ov \bZ[1/N])$ of representations of $\ZG$. This morphism is an injection, and we will identify $V(a)$ with its image in $\Lie(\ZZG)(\ov \bZ[1/N])$.

\begin{Thm} \label{Thm-Weyl-general}
	Let $G$ be a good reductive group over $K$, $r$ a real number and $x$ a rational point of $\sB(G,K)$. Let
		$
		\Nprime = \left\{ \begin{array}{ll} 2N & \text{ if $\Phi_K$ contains multipliable roots } \\
		N & \text{ otherwise } \, .
		\end{array} \right. 
		$
		
		Then
	\begin{equation} \label{eqn-ZV-Weyl}
		\ZV_{\ov \bZ[1/\Nprime]}  \simeq  (\ZV_T)_{\ov \bZ[1/\Nprime]}  + \bigoplus_{\lambda \in \Phi_{x,r}^{\rm{max}} } V(\lambda)_{\ov \bZ[1/\Nprime]}  \subset \Lie(\ZZG)(\ov \bZ[1/\Nprime]) 
	\end{equation}
	as representations of $\ZG_{\ov \bZ[1/\Nprime]}$.
\end{Thm}	

\textbf{Proof.} 
By the definition of $\Nprime$ the subspace $\ZV_{\ov \bZ[1/\Nprime]} \subset \Lie(\ZZG)(\ov \bZ[1/\Nprime])$ is spanned by $\ZV_T$ and $Y_a$ for $a \in \Phi_{x,r}$ (Section \ref{Section-global-filtration-quotient}). Thus, analogously to the argument in the proof of Theorem \ref{Thm-Weyl-split}, it suffices to show that $\<\ZG(\ov \bZ[1/\Nprime])(Y_a), \ZV_T\>_{a \in \Phixr}$ contains $Y_b$ for all $b \in \Phi_{x,r}$. Let $a \in \Phixr \backslash \Phi_{x,r}$, $b \in \Phi^+(\ZG)$ with $a+b \in \Phi_{x,r}$, and $N_{a,b}>0$ the maximal integer such that $a+N_{a,b}b \in \Phi_{x,r}$. We need to show that $Y_{a+N_{a,b}b} \in \<\ZG(\ov \bZ[1/\Nprime])(Y_a), \ZV_T\>_{a \in \Phixr}$ implies $Y_a \in \<\ZG(\ov \bZ[1/\Nprime])(Y_a), \ZV_T\>_{a \in \Phixr}$.
 We assume $Y_{a+N_{a,b}b} \in \<\ZG(\ov \bZ[1/\Nprime])(Y_a), \ZV_T\>_{a \in \Phixr}$ and distinguish four cases.
\begin{description}
\item{Case 1:} $a \bR \neq b \bR$ and $b$ is not multipliable. In this case the result follows from the proof of the split case (Theorem \ref{Thm-Weyl-split}) and Equation \eqref{eqn-DefProp3} on page \pageref{eqn-DefProp3} and Equation \eqref{eqn-DefYa} on page \pageref{eqn-DefYa} (if $b$ is non-divisible) or Equation \eqref{eqn-DefProp2} on page \pageref{eqn-DefProp2} and Equation \eqref{eqn-DefYa} (if $b$ is divisible).

\item{Case 2:} $a \bR = b \bR$ and $b$ is not multipliable. In this case $a=-(a+N_{a,b}b)$, and the element $s_b$ in the Weyl group of $\ZG$ corresponding to reflection in direction of $b$ sends $Y_{a+N_{a,b}b}$ to $\pm Y_{-(a+N_{a,b}b)}=\pm Y_{a}$. Hence $Y_{a} \in \<\ZG(\ov \bZ[1/\Nprime])(Y_a), \ZV_T\>_{a \in \Phixr}$.

\item{Case 3:} $a \bR \neq b \bR$ and $b$ is multipliable. By taking Galois orbits over different connected component and using Equation \eqref{eqn-DefProp1} on page \ref{eqn-DefProp1} and Equation \eqref{eqn-DefYa} on page \pageref{eqn-DefYa}, it suffices to consider the case that $\Dyn(G)=A_{2n}$ with non-trivial Galois action. We label the simple roots by $\alpha_{n}, \alpha_{n-1}, \hdots, \al_2, \alpha_1, \beta_1, \beta_2, \hdots, \beta_{n}$ as in Figure \ref{fig-A2n1} on page \pageref{fig-A2n1}. Then $b$ is the image of $\alpha_1 + \hdots + \alpha_s$ for some $1 \leq s \leq n$, and, as $\<\check b, a+N_{a,b}b\> >0$, the root $a+N_{a,b}b$ is the image of
\begin{itemize}[label={}]
	 \item $-(\alpha_{s+1} +\hdots + \alpha_{s_1})$ for some $s < s_1 \leq n$, or
	 \item $\alpha_{s_2} + \hdots + \alpha_s$ for some $1 < s_2 \leq s$, or
	 \item $\alpha_{1} + \hdots + \alpha_s + \beta_1 + \hdots + \beta_{s_3}$ for some $1 \leq s_3 < s$ or $s < s_3 \leq n$.
\end{itemize}

To simplify notation, we will prove the claim for the case that $b$ is the image of $\alpha_1$ and $a+N_{a,b}b$ is the image of $-\alpha_2$. All other cases work analogously. Combining Equation \eqref{eqn-DefProp1} on page \pageref{eqn-DefProp1}, Equation \eqref{eqn-DefYa} on page \pageref{eqn-DefYa} and Equation \eqref{eqn-action} on page \pageref{eqn-action}, and using that $\ZG_{\ov \bZ[1/\Nprime]}$ preserves the subspace $\ZV_{\ov \bZ[1/\Nprime]}$ of $\Lie(\ZZG)(\ov \bZ[1/\Nprime])$, we obtain that
\begin{eqnarray*}
 \Zx_{-b}(u)(Y_{a+N_{a,b} b}) & = & \left(\ZZx_{-\beta_1}(\sqrt{2}u)\ZZx_{-(\al_1+\beta_1)}(-(-1)^{b(x-x_0)M}u^2) \ZZx_{-\al_1}((-1)^{b(x-x_0)M}\sqrt{2}u) \right) \\
 & &  \left(\ZZX_{-\beta_2}+(-1)^{(-(a+N_{a,b}b)(x_q-x_{0,q})+r)\cdot 2} \ZZX_{-\alpha_2}\right)	\\
 & = & Y_{a+N_{a,b}b} + m_{a,b,1}' \sqrt{2} u Y_{a+(N_{a,b}-1)b} + m_{a,b,2}' u^2 Y_{a+(N_{a,b}-2)b} \, 
\end{eqnarray*}
with $m_{a,b,1}', m_{a,b,2}' \in \{\pm1\}$, for all $u \in \bG_a(\ov \bZ[1/\Nprime])$.  Since $2 \mid \Nprime$, taking $\ov\bZ[1/\Nprime]$-linear combinations of $Y_{a+N_{a,b}b} + m_{a,b,1}' \sqrt{2} u Y_{a+(N_{a,b}-1)b} + m_{a,b,2}' u^2 Y_{a+(N_{a,b}-2)b}$ for different $u$ implies that $Y_{a+(N_{a,b}-1)b}$ and $Y_{a+(N_{a,b}-2)b}$ are contained in $\<\ZG(\ov \bZ[1/\Nprime])(Y_a), \ZV_T\>_{a \in \Phixr}$, so $Y_a \in \<\ZG(\ov \bZ[1/\Nprime])(Y_a), \ZV_T\>_{a \in \Phixr}$.

\item{Case 4:} $a \bR = b \bR$ and $b$ is multipliable. As in Case 3, we can restrict to the case that $\Dyn(G)=A_{2n}$, and we may assume that $b$ is the image of $\al_1$. Then $a+N_{a,b}b$ is the image of $\alpha_1$ or the image of $\alpha_1+\beta_1$. If $N^-_{a,b}$ denotes the largest integer such that $a-N^-_{a,b}b \in \Phi_{x,r}$, then $Y_{a-N^-_{a,b}b}$ is conjugate to $\pm Y_{a+N_{a,b}b}$ under the Weyl group. Hence $Y_{a-N^-_{a,b}b} \in \<\ZG(\ov \bZ[1/\Nprime])(Y_a), \ZV_T\>_{a \in \Phixr}$. If $a+N_{a,b}b$ is the image of $\alpha_1$, then $N^-_{a,b}=0$, and we are done. Thus, suppose that $a+N_{a,b}b$ is the image of $\alpha_1+\beta_1$.

 Recall that for $\alpha \in \Phi$ and $H_\alpha:=\Lie(\check \alpha)(1)$, we have (\cite[Corollary~5.1.12]{ConradSGA3})
\begin{eqnarray*}
	\ZZx_{-\al}(u)(\ZZX_\al)&=&\ZZX_\al+\eps_{\al,\al}uH_{-\alpha}-\eps_{\al,\al}u^2\ZZX_{-\al} \\
	\ZZx_{-\al}(u)(H) & = & H + \Lie(\al)(H)u\ZZX_{-\al}
\end{eqnarray*}
for all $u \in \bG_a(\ov \bZ[1/\Nprime])$ and all $H \in \Lie(\ZZT)(\ov \bZ[1/\Nprime])$. Using these identities, we obtain
\begin{eqnarray} 
	\Zx_{-b}(u)(Y_{a+N_{a,b}b})&=&\left(\ZZx_{-\beta_1}(\sqrt{2}u)\ZZx_{-\al_1-\beta_1}(-(-1)^{b(x-x_0)M}u^2)\ZZx_{-\al_1}((-1)^{b(x-x_0)M}\sqrt{2}u) \right) \notag \\
	& &  \left( \ZZX_{\al_1+\beta_1} \right) \notag \\
	&= & Y_{a+N_{a,b}b}+m_{a,1}''\sqrt{2}uY_{a+(N_{a,b}-1)b}+H+m_{a,3}''\sqrt{2}u^3Y_{a+(N_{a,b}-3)b} \notag \\
	&  &  +m_{a,4}''u^4Y_{a+(N_{a,b}-4)b}, \notag 
\end{eqnarray}
with $m_{a,1}'', m_{a,3}'' \in \{\pm1\}$, $m_{a,4}'' \in \{\pm1,\pm3\}$ and $H\in \ZV_T$. As $Y_{a+(N_{a,b}-4)b}=Y_{a-N^-_{a,b}b}$ and $H$ are in $\<\ZG(\ov \bZ[1/\Nprime])(Y_a), \ZV_T\>_{a \in \Phixr}$, and since $2 \mid \Nprime$, we also obtain that $Y_{a+(N_{a,b}-1)b}$ and $Y_{a+(N_{a,b}-3)b}$ are contained in $\<\ZG(\ov \bZ[1/\Nprime])(Y_a), \ZV_T\>_{a \in \Phixr}$. \qed
\end{description}

\begin{Cor}
	Let $G$ be a good reductive group, $r \notin \frac{1}{p^sN}\bZ$ a real number, and $x$ a rational point of $\sB(G,K)$. 
	 Then
	\begin{equation*}
		\PV_{x,r} \simeq	\bigoplus_{\lambda \in \Phi_{x,r}^{\rm{max}} } V(\lambda)_{\fpb} \, .
	\end{equation*}
\end{Cor}	
\textbf{Proof.} If $r \notin \frac{1}{p^sN}\bZ$, then $\ZV_T=\{0\}$. Hence, if $p \neq 2$, the claim follows by combining Theorem \ref{Thm-Weyl-general} and Theorem \ref{Thm-general}. The proof in the case $p=2$ works completely analogous to the proof of Theorem \ref{Thm-Weyl-general} using that $\ZV$ is spanned by $\ZV_T$, $Y_a$ for all $a \in \Phi_K^{\mathrm{nm}}$ with $r-a(x-x_0) \in \Gamma_a'(G)$ and $\sqrt{2}Y_a$ for all $a \in \Phi_K^{\mathrm{mul}}$ with $r-a(x-x_0) \in \Gamma_a'(G)$. \qed

\printindex{definition}{Selected definitions}

\printindex{notation}{Selected notation}

\bibliographystyle{amsxport}

\bibliography{Jessisbib}

\end{document}